\input elliptic.def

%

\NoLabels

\MakeToc{toc}

\topmatter

\title
On deformation types of real elliptic surfaces
\endtitle


\dedicatory
\bigskip
\rightline{\epsfile{iching}}
\medskip
\endgraf\leftskip0pt plus.4\hsize \rightskip0pt
  \parfillskip0pt \parindent0pt \spaceskip.3333em \xspaceskip.5em
  \pretolerance9999\tolerance9999 \exhyphenpenalty10000
  \hyphenpenalty10000
\rm
Le Yi Jing n'est pas un livre, un texte qu'on lit du d\'{e}but \`{a} la fin,
mais un ouvrage que l'on consulte quand on en a besoin. Lorsqu'on
h\'{e}site sur une voie \`{a} suivre, une attitude \`{a} prendre, un choix \`{a} faire,
un dilemme \`{a} r\'{e}soudre, on peut alors s'en servir pour ce
qu'il est dans la pratique~: un manuel d'aide \`{a} la d\'{e}cision.
\endgraf\smallskip
Cyrille Javary, {\it Les Rouages du Yi Jing}, Ed. Phillipe
Picquier, 2001\mnote{Kh: to move right the epigraph}
\enddedicatory

\author
Alex Degtyarev, Ilia Itenberg, Viatcheslav Kharlamov
\endauthor

\address
Bilkent University\endgraf\nobreak
06800 Ankara, Turkey
\endaddress

\email
degt\@fen.bilkent.edu.tr
\endemail

\address
Universit\'{e} Louis Pasteur et IRMA (CNRS)\endgraf\nobreak
7 rue Ren\'{e} Descartes 67084 Strasbourg Cedex, France
\endaddress

\email
itenberg\@math.u-strasbg.fr
\endemail

\address
Universit\'{e} Louis Pasteur et IRMA (CNRS)\endgraf\nobreak
7 rue Ren\'{e} Descartes 67084 Strasbourg Cedex, France
\endaddress

\email
kharlam\@math.u-strasbg.fr
\endemail

\thanks
The second and third authors were partially funded by the ANR-05-0053-01
grant of Agence Nationale de la Recherche and a grant of
Universit\'{e} Louis Pasteur, Strasbourg.
\endthanks

\keywords
Real elliptic surface, \emph{dessins d'enfants}, Tate-Shafarevich group
\endkeywords

\subjclass
14J27, 14P25, 05C90
\endsubjclass

\contents
 \nosubsubsection
\endcontents

\abstract
We study real elliptic surfaces and trigonal curves
(over a base of an arbitrary genus)
and their equivariant deformations. We calculate
the real Tate-Shafarevich group and reduce
the deformation classification to
the combinatorics of a real version
of Grothendieck's \emph{dessins d'enfants}.
As a consequence,
we obtain an explicit description
of the deformation classes of $M$- and
$(M-1)$- (\ie, maximal and submaximal in the sense
of the Smith inequality) curves and surfaces.
\endabstract

\endtopmatter

\begin
\document

\begingroup
\def\include#1{\input #1\relax}
\let\end\endinput


\ifx\defloaded\undefined\input elliptic.def \fi

\begin

\section{Introduction}

\subsection{Motivation and historical remarks}
In geometry of nonsingular algebraic surfaces, over the reals as
well as over the complex numbers, there are two major equivalence
relations: the first one, called \emph{deformation equivalence},
is up to isomorphism and deformation (of the complex structure),
and the second one, called \emph{topological equivalence}, is up
to diffeomorphism (ignoring the complex structure). Certainly,
deformation equivalence implies topological equivalence, and
one of the principal questions in the subject is to what extent
the converse holds, \ie, to what extent is the deformation
class of a surface controlled by its topology.
Since we regard a real variety as a complex variety
equipped with a real structure (which is an
anti-holomorphic involution), by
a deformation of real varieties we mean
an equivariant Kodaira-Spencer deformation,
and by a diffeomorphism
between two
real varieties we mean an equivariant diffeomorphism. Therefore,
the $\operatorname{Dif}=\operatorname{Def}$ question above
stated over the reals would involve
the same question for the underlying complex varieties.
Luckily, due to Donaldson's and Seiberg-Witten's revolution in
four dimensional topology
(as well as the Enriques-Kodaira classification of algebraic
surfaces),
one does have
an advanced level of control over the discrepancy between the deformation
class of a compact complex surface and its diffeomorphism class.
This fact
makes it reasonable to fix a deformation class of compact
complex surfaces beforehand and to concentrate on the topology
and deformations of the real
structures that can appear
on (some of) the surfaces in question.

The problem of enumerating the equivariant deformation classes of
real structures within a fixed complex deformation class goes
back at least to
F.~Klein \cite{Klein},
who studied the nonsingular real cubic
surfaces in~$\Cp3$
(\ie, Del Pezzo surfaces of degree~$3$) from a similar point of view.
He proved that
the equivariant deformation class
of such a surface is already determined by the topology of its real part,
which is the real projective plane with up to three handles or up to
one sphere (Schl\"{a}fli's famous five `species' of nonsingular
cubics).
Further important steps in this direction were made
by A.~Commessatti~\cite{Co-rat}, \cite{Co-abel}, who found a
classification of all
real abelian
surfaces and all $\R$-minimal real rational surfaces, thus
extending (at least implicitly) Klein's result
to these special classes.

In general, we call a deformation class of complex varieties
\emph{quasi-simple} if a real variety
within the complex class is determined up to equivariant
deformation by the diffeomorphism type of the real structure.
For
curves, the problem
was settled by F.~Klein and G.~Weichold
(see, \eg, the survey~\cite{Nat}) who proved that the family of compact curves of
any given genus is indeed quasi-simple.
Note that the equivariant
deformation class of a
real curve
is no
longer determined by its genus and real part; in addition, one should take
into account the so called \emph{type} of the curve, \ie, whether
the real part does or does not divide the complexification.
However, the type is certainly a topological invariant of the real
structure.

Further advance in the study of quasi-simplicity called for
appropriate tools in complex algebraic geometry. Their development
took half a century, and it was not until the late 70s that the
study was resumed.
Now, due to the results
obtained in~\cite{Nikulin}, \cite{DK-rat}, \cite{DIK1},
\cite{Welschinger}, \cite{CF}, \cite{duke}, we know that
quasi-simplicity holds for any special (in the sense of the
Enriques-Kodaira classification) class of $\C$-minimal complex
surfaces except elliptic. (For the surfaces of general type there
are counter-examples, see, \eg,~\cite{KK}.)
A slightly different but related
\emph{finiteness} statement, \ie, finiteness of the number of
equivariant deformation classes of real structures within a given
deformation class of complex varieties,
is known to hold for all surfaces except elliptic
or ruled with irrational base. For ruled surfaces, the statement
is probably true and its proof should not be difficult, \cf., \eg,
\cite{DK-rat}, but it does not seem to appear in the literature.
Thus, elliptic surfaces are essentially
the last special class of surfaces
for which the quasi-simplicity and finiteness questions are
still open.

It is worth mentioning that, in spite of noticeable activity in
the theory of complex elliptic surfaces,
literature
dealing with the real case
is scant. Among the few works that we know
are~\cite{Mn-Akriche}, \cite{Ballico}, \cite{Mn-Bihan}, \cite{DK-empty},
 \cite{Frediani}, \cite{Gross-Wilson},
 \cite{Kharlamov}, \cite{Mangolte}, \cite{Silhol}, and~\cite{Wall}.


\subsection{Subject of the paper}
In this paper, our goal is to study
relatively minimal real elliptic surfaces without multiple fibers and, in
particular,
to understand
the extent to which
the
equivariant deformation class of
such a surface is controlled by the
topology of its real structure.
Recall that the complex deformation class of an elliptic surface
as above is determined by the genus~$g$ of the base curve and the
Euler characteristic~$\chi$ of the surface, provided that~$\chi$ is
positive. (The case $g=0$ is treated in A.~Kas~\cite{Kas}, and the
general case, in
W.~Seiler~\cite{Sei}, see also~\cite{FM}.) Note that, if $\chi$ is
small (for a given genus), one
deformation class may consist of several irreducible components:
the principal component formed by the non-isotrivial
surfaces
may be accompanied by few others, formed by the isotrivial ones.
Each isotrivial surface can
be deformed to a surface that perturbs to a non-isotrivial one.
However, from the known constructions it is not immediately obvious
that the deformation can be chosen real. For this reason,
we confine
ourselves to the more topological study of
non-isotrivial surfaces,
leaving the algebro-geometric aspects to subsequent papers.

An elliptic surface comes equipped with an elliptic fibration.
Moreover, for most surfaces, in particular, for all
elliptic surfaces
of Kodaira dimension~$1$, the elliptic fibration is unique.
(In the case of relatively minimal surfaces without multiple
fibers,
the Kodaira dimension is known to be equal to~$1$ whenever
$g>0$, as well as when $g=0$ and the Euler characteristic~$\chi$,
which is divisible by~$12$, is $>24$.)
Thus, the elliptic fibration is an important part of the structure,
and we include it into the setting of the
problem,
considering equivariant deformations of real elliptic fibrations
(with no confluence of singular fibers allowed) on the one hand,
and equivariant diffeo-/homeomorphisms on the other
hand.
Furthermore, as any non-isotrivial surface can be perturbed to an
\emph{almost generic} one, \ie, a surface with simplest singular
fibers only, we consider solely deformations of almost generic surfaces.
Here, `almost generic' can be thought of as `topologically
generic', as opposed to `generic', or `algebraically generic',
where one requires in addition that the fibers with nontrivial
complex multiplication should also be simple. We use the latter
assumption when treating an individual surface \via\
algebro-geometric tools.

Note that during the deformation we never assume the base curve
fixed; it is also subject to a deformation.
The classification of real elliptic surfaces over a fixed base does
not seem feasible; in general it may not even be possible to
perturb a given surface to an almost generic one.

\subsection{Tools and results}
As in the complex case, the study of real elliptic surfaces
is based upon two major tools:
the real version of the
Tate-Shafarevich group, which enumerates all real surfaces with a given
Jacobian, and a real version of the techniques of
\emph{dessins d'enfants}, which reduces the deformation
classification of non-isotrivial Jacobian elliptic surfaces (or,
more generally, trigonal curves on ruled surfaces) to a
combinatorial problem. We develop the two tools and,
as a first application,
obtain a rather explicit classification of the so called $M$- and
$(M-1)$-surfaces (\ie, those maximal and submaximal in the sense
of the Smith inequality) and $M$-
and $(M-1)$-curves. The principal results of the paper are stated
in~\ref{class-M}--\ref{S.M-1}.

As a straightforward consequence of the description of deformation
classes in terms of groups and graphs, the whole number of
equivariant deformation classes of real elliptic fibrations
with given numeric invariants
is finite. This
settles the finiteness problem stated
above for
non-isotrivial elliptic
surfaces without multiple fibers.

The real Tate-Shafarevich group~$\Rsha(J)$ is defined as the set
(with a certain group operation) of the isomorphism classes of all
real elliptic fibrations with a given Jacobian~$J$. Contrary to
the complex case, $\Rsha(J)$ is usually disconnected, and we
describe, in purely topological terms, its discrete part
$\Rsha\top(J)$, which enumerates the deformation classes of real
fibrations whose Jacobian is~$J$. This description
gives an explicit list of
all modifications that a fibration may undergo, and
the result shows
that they can all be seen in the real part. As a
consequence,
we prove
that, up to deformation,
an elliptic surface is Jacobian if and only if the real part of
the fibration
admits a topological section (Proposition~\ref{section->Jacobian}),
each Betti
number of an elliptic surface is bounded by the corresponding
Betti number of its Jacobian, and each $M$- or $(M-1)$-surface is
Jacobian up to deformation
(Proposition~\ref{M->Jacobian}).

The real version of \emph{dessins d'enfants} was first introduced by
S.~Orevkov~\cite{Orevkov}, who used
it to study real trigonal curves on $\C$-minimal rational
surfaces~$\Sigma_d$.
(A similar object was considered independently
in~\cite{SV} and~\cite{NSV}).
The curves considered by Orevkov do not
intersect the `exceptional' section of the surface; for even
values of~$d$, these are the branch curves of the Weierstra{\ss}
models of Jacobian elliptic surfaces. (In the Weierstra{\ss} model,
the elliptic surface appears as the double covering of~$\Sigma_d$
branched at the union of the exceptional section and a trigonal
curve.) Using the dessin techniques, Orevkov invented a kind of
Viro-LEGO$^\circledR$ game: he introduced a few elementary pieces
(which are the dessins of cubic curves), defined the operation of
connecting `free ends' of two pieces, and used this procedure to
construct bigger curves, thus proving a number of existence
statements.
Orevkov also noticed
(private communication~\cite{Orevkov-private}; \cf. similar
observations in V.~Zvonilov~\cite{Zvonilov}) that, as long as
almost generic $M$-curves over a rational base are concerned, this
procedure is universal, \ie, there is a unique way to break any
$M$-curve into elementary pieces. Clearly, this construction gives
a deformation classification of such $M$-curves.

We extend Orevkov's approach to trigonal curves over a base of
an arbitrary genus and obtain similar results for $M$- and
$(M-1)$-curves. We show that, as in the rational case, any $M$- or
$(M-1)$- curve breaks into certain elementary pieces.
(An essential ingredient here is Theorem~\ref{6k<=}, which states
that unbreakable curves must be sufficiently `small'.
Another decomposability statement, Theorem~\ref{reduction},
is used to handle large pieces of $(M-1)$-curves.)
In the
$M$-case, this procedure is unique; in the $(M-1)$-case it is
unique up to a few moves that are described explicitly. As a
consequence, we obtain a deformation classification of $M$- and
$(M-1)$-curves and, when combined with the results on~$\Rsha\top$,
that of $M$- and $(M-1)$-surfaces (see~\ref{th.M-surfaces}
and~\ref{M.g=0} for the $M$-case and \ref{M-1.surface}
and~\ref{M-1.g=0} for the $(M-1)$-case). A surprising by-product
of the classification is the fact that, essentially, $M$- and
$(M-1)$-surfaces and curves exist only over a base of genus
$g\le1$. (Here, `essentially' means that certain `trivial' handles
should be ignored. Without this convention the genus can be made
arbitrary large.)

\subsection{Contents of the paper}
Sections~\ref{irs} and~\ref{res} are introductory.
In Section~\ref{irs} we remind the reader a few basic facts concerning
topology of involutions, and
in Section~\ref{res} we discuss certain complex and real aspects
of the theory of trigonal curves, elliptic surfaces,
their Jacobians and Weierstra{\ss} models.
In Section~\ref{Tate-Shafarevich} we introduce a real version
of the Tate-Shafarevich
group, express it in cohomological terms, and study its discrete
part. The main results here are Theorems~\ref{RSha.BR}
and~\ref{Rsha(JR)}, as well as their corollaries.
Section~\ref{rtc} plays
a central r\^{o}le in the paper. Here we develop Orevkov's results
on real \emph{dessins d'enfants}. After a brief introduction, we
concentrate on a special class of dessins that represent
meromorphic functions having generic branching behavior, \ie,
$j$-invariants of generic trigonal curves. The principal results
of Section~\ref{rtc} are the decomposability
theorems~\ref{reduction} and~\ref{6k<=}, which assert that, under
certain assumptions, a dessin breaks into simple pieces. Finally,
in Section~\ref{applications} we apply the results obtained to the
case of $M$- and $(M-1)$-curves and surfaces. We prove the
structure theorems, derive a few simple consequences, and discuss
further generalizations and open problems.

\subsection{Acknowledgements}
Our thanks go to Stepan Orevkov, who enthusiastically
shared his observations with us, motivating
our interest in \emph{dessins d'enfants}.
We would also like to thank Victoria Degtyareva
for courageously reading and polishing a preliminary version
of the text.

We are grateful to
the
{\it Max-Planck-Institut f\"{u}r Mathematik\/}
and to the
{\it Mathematisches Forschungsinstitut Oberwolfach\/} and its RiP program
for their hospitality and excellent working
conditions which helped us to
complete
this project.
An essential part of the work was done during the
first
author's visits to \emph{Universit\'{e} Louis Pasteur}, Strasbourg.

\end


\ifx\defloaded\undefined\input elliptic.def \fi

\begin

\setcounter\section1

\section{Involutions and real structures}\label{irs}

In this section we recall basic results concerning topology of
involutions. Proofs and further details can be found in
the monograph~\cite{Bre1},
which deals with
general theory of compact transformation groups.
A survey of sheaf theory, cohomology, and spectral sequences
is found in~\cite{Bre2}.
For a self-contained exposition specially tailored
for the needs of topology of real algebraic varieties,
we refer to~\cite{DIK1}.

\subsection{Real structures and real sheaves}
Throughout this section
all topological spaces are assumed
paracompact
and Hausdorff.

\paragraph\label{real.structure}
A \emph{real structure} on a complex
variety~$X$
(not necessarily
connected or
nonsingular) is an anti-holomorphic involution $c_X\:X\to X$.
Clearly, any two real structures differ by an automorphism of~$X$.
By the Riemann extension theorem, an
\hbox{(anti-)}holomorphic
endomorphism~$f$ of the smooth part of~$X$
extends to an
\hbox{(anti-)}\penalty100 holomorphic
endomorphism of~$X$
if
and only if
$f$ admits a continuous extension.

A \emph{real variety}
is a complex variety~$X$ equipped with a real
structure~$c_X$. (Sometimes it is convenient to refer to the pair $(X,c_X)$
as a \emph{real form} of~$X$.) The fixed point set $\Fix c_X$ is called the
\emph{real part} of~$X$ and is denoted $X_{\R}$. A holomorphic map
$f\:X\to Y$ between two real varieties $(X,c_X)$ and $(Y,c_Y)$ is
called \emph{real} if it commutes with the real structures:
$c_Y\circ f=f\circ c_X$.

Recall that for any
continuous involution~$c_X$ on
a finite dimensional
topological
space~$X$ with
finitely generated
total cohomology group
$H^*(X;\Z_2)$ the following
\emph{Smith inequality} holds:
$$
\dim H^*(\Fix c_X;\Z_2)\le \dim H^*(X;\Z_2).
$$
Furthermore, the difference
$\dim H^*(X;\Z_2)-\dim H^*(\Fix c_X;\Z_2)$ is even.
If the difference is~$2d$, the
involution~$c_X$ is called an \emph{$(M-d)$-involution}. If $c_X$
is the real structure of a real variety~$X$, then $X$ itself is
called an \emph{$(M-d)$-variety}.

\paragraph\label{Z2cohomology}
Given an abelian group~$A$ with an involution $c\:A\to A$, we
define the cohomology groups $H^*(\Z_2;A)$ to be the cohomology of
the complex
$$
0@>>>A@>1-c>>A@>1+c>>A@>1-c>>\ldots
$$
(the leftmost copy of~$A$ being of degree zero). Similarly, given
a sheaf~$\CA$ with an involutive automorphism $c\:\CA\to\CA$, we
define the cohomology sheaves $\CH^*(\Z_2;\CA)$ to be the
cohomology of the complex
$$
0@>>>\CA@>1-c>>\CA@>1+c>>\CA@>1-c>>\ldots.
$$
(Certainly, the former
is nothing
but
a specialization
of
the
general definition of the cohomology of a discrete group~$G$
with coefficients in a $G$-module, see,
\eg,~\cite{Bro},
to the group $G=\Z_2$ and the simplest invariant cell decomposition of
the space $S^\infty=E\Z_2$.
The
latter
is a straightforward sheaf version of the
former.)

\paragraph\label{rsheaf}
Let $X$ be a topological space with an involution $c_X\:X\to X$.
Denote by $\pi\:X\to X/c_X$ the projection. Given a sheaf~$\CA$
on~$X$, any morphism $c\:\CA\to c^*_X\CA$ (over the identity
of~$X$) descends to a morphism $\pi_*c\:\pi_*\CA\to\pi_*\CA$. By a
certain abuse of the language, $c$~is called an \emph{involution}
if $\pi_*c$ is an involution. This condition is equivalent to the
requirement $c\circ c_X^*c=\id$, where $c_X^*c\:c^*_X\CA\to\CA$ is
the pull-back of~$c$. (Certainly, this definition is merely an
attempt to refer to involutive lifts $\CA\to\CA$ of~$c_X$ to~$\CA$
in terms of sheaf morphisms identical on the base.)

The constant sheaf $G_X$ (for any abelian group~$G$) has a
canonical involution, which is the identity
$G_X=c_X^*G_X$. As a lift of~$c_X$ it is given by
$s\mapsto s\circ c_X$.

\paragraph\label{rcanon}
Now, let~$X$ be a complex manifold and let $c_X\:X\to X$ be a real
structure. Then the structure sheaf~$\CO_X$
has a canonical
involution, called the \emph{canonical real structure} (defined
by~$c_X$); it is given by the complex conjugation
$\CO_X=\overline{c^*_X\CO_X\!}\,$, or,
as a lift of~$c_X$, by
$s\mapsto\overline{s\circ c_X}$. If
$\CA$ is a coherent sheaf on~$X$, the pull-back $c^*_X\CA$ is a
(coherent, in a sense) sheaf of $c^*_X\CO_X$-modules. An
involution $c\:\CA\to c^*_X\CA$ is called a \emph{real structure}
on~$\CA$ if it is compatible with the module structure (\via\
the canonical real structure on~$\CO_X$). A typical example is the
canonical real structure on the
sheaf~$\CF$ of sections of a \emph{`Real'}
vector bundle~$F$ (\ie, a holomorphic
vector bundle on~$X$ supplied with an involution~$c_F$
covering~$c_X$ and anti-linear on the fibers, so that it is a real
structure on the total space); it is given by
$s\mapsto c_F\circ s\circ c_X$. This formula applies as well in a
more general situation, when $F\to X$ is a holomorphic fibration
with
abelian groups
as fibers
(so that $\CF$ is a sheaf of abelian
groups) and $c_F\:F\to F$ is a fiberwise additive real structure
covering~$c_X$. Although, in general, $\CF$ is not a
coherent sheaf, we will still refer to the result as the
\emph{canonical real structure} on~$\CF$.

The canonical real structure on~$\CO_X$ defines involutions on the
other two members of the exponential sequence
$$
0@>\phantom{2\pi i}>>\Z_X@>2\pi i>>\CO_X
 @>\phantom{2\pi i}>>\CO_X^*@>\phantom{2\pi i}>>0.
$$
Note that the resulting involution on the constant sheaf~$\Z_X$
differs from the canonical involution above by~$(-1)$.
In order to emphasize this nonstandard real structure, we will use
the notation $\Z_X^-$ (and, more generally,~$G_X^-$).

\paragraph
Let $\CA$ be a sheaf on~$X$ with an involution
$c\:\CA\to c_X^*\CA$. Denote by $\pi\:X\to
X/c_X$ the projection
and consider the complex
$$
\pi_*\CA^*:\quad0@>\phantom{1-c}>>\pi_*\CA@>1-c>>\pi_*\CA@>1+c>>\pi_*\CA@>1-c>>\ldots
$$
of sheaves on $X/c_X$ (the leftmost copy of~$\pi_*\CA$ being of
degree
zero; for simplicity, we use the same notation~$c$ for the automorphism
$\pi_*c\:\pi_*\CA\to\pi_*\CA$).
We will refer to the hypercohomology $\HH^*(X/c_X;\pi_*\CA^*)$
as the hypercohomology of~$(\CA,c)$ and denote it $\HH^*(\CA,c)$
(or just $\HH^*(\CA)$, when $c$ is understood). For the constant
sheaf~$G_X$ with its canonical real structure we will also use the
notation $\HH^*(X;G)$.

Recall that there are natural spectral sequences
$$
H^q(X/c_X;\CH^p(\Z_2;\CA))\Longrightarrow\HH^{p+q}(\CA),\eqtag\label{GaloisI}
$$
where $\CH^p(\Z_2;\CA)$ stand for the cohomology sheaves
of~$\pi_*\CA^*$, and
$$
H^p(\Z_2;H^q(X/c_X;\pi_*\CA))=H^p(\Z_2;H^q(X;\CA))
 \Longrightarrow\HH^{p+q}(\CA).\eqtag\label{GaloisII}
$$
(Since $\pi$ is finite-to-one, the higher direct images
$\RR{i}\pi$,
$i>0$, vanish and one has $H^q(X/c_X;\pi_*\CA)=H^q(X;\CA)$.)
Furthermore, since $\pi$ is finite-to-one,
one can calculate $\HH^*(\CA)$ using
$c_X$-invariant \v{C}ech resolutions of~$\CA$. More precisely,
given a $c_X$-invariant open covering $\CU=\{U_i\}$ of~$X$, one
can consider the bi-complex
$$
(\check C^{p,*}_{\CU},d_1,d_2)=
 \bigoplus_{p\ge0}(\check C_{\CU}^*(\CA),d_2),\qquad
d_1=1-(-1)^pc\:C^{p,*}_{\CU}\to C^{p+1,*}_{\CU}
$$
(direct sum of copies of the ordinary \v{C}ech complex with the
first differential given above). Then $\HH^n(\CA)$ is the limit,
over all coverings, of the cohomology
$H^n(\check C^{*,*}_{\CU})$.

\subsection{Kalinin's spectral sequence}\label{Kalinin}
Let~$X$ be a
a finite dimensional
paracompact
Hausdorff topological space with an involution $c_X\:X\to X$.

\paragraph\label{Borel}
The \emph{Borel construction} over $(X,c_X)$ is the twisted
product
$$
X_c=X\times_cS^\infty=(X\times S^\infty)/(x,\AMSbold r)\sim(c_X(x),-\AMSbold r).
$$
The cohomology groups $H^*_c(X;G)=H^*(X_c;G)$
are called the
\emph{equivariant cohomology} of~$X$ (with coefficients in an
abelian group~$G$).
Note that the
subscript~$c$ stands for
the involution $c=c_X$;
as we never use cohomology with compact supports, we hope that
this notation will not lead
to a confusion.
The Leray spectral sequence of the fibration
$X_c\to\Rpi=S^\infty/\pm\id$
with fiber~$X$ is called the \emph{Borel-Serre
spectral sequence} of $(X,c_X)$:
$$
\RE2^{pq}(X;G)=H^p(\Z_2;H^q(X;G))\Longrightarrow H^{p+q}_c(X;G).
$$
Sometimes it is convenient to start the sequence at the term
$\RE1^{pq}=H^q(X;G)$ with the differential
$\Rd1^{p*}=1-(-1)^pc^*_X$.

There is a canonical isomorphism $H^p_c(X;G)=\HH^p(X;G)$, and
the Borel-Serre spectral sequence is isomorphic to the spectral
sequence~\eqref{GaloisII} for the constant sheaf $\CA=G_X$. If $G$
is a commutative ring, then the Borel-Serre spectral sequence is a
spectral sequence of $H^*(\Rpi;G)$-algebras.

\paragraph\label{K.Z2}
Let~$G=\Z_2$, and let $\hh\in H^1(\Rpi;\Z_2)=\Z_2$ be the generator.
Assume, in addition,
that $X$ is a {\sl CW}-complex of finite
dimension. Then the \emph{stabilization homomorphisms}
$$
\scup\hh\:\rE^{pq}(X;\Z_2)\to\rE^{p+1,q}(X;\Z_2),\quad
\scup\hh\:H^n_c(X;\Z_2)\to H^{n+1}_c(X;\Z_2)
$$
are isomorphisms for $p\gg0$ and one has
$$
\lim_{n\to\infty}H^n_c(X;\Z_2)
 =H^{n\gg0}(\Fix c_X\times\Rpi;\Z_2)
 =H^*(\Fix c_X;\Z_2).
$$
Thus, one obtains a $\Z$-graded spectral sequence
$\rH^q(X;\Z_2)=\lim_{p\to\infty}\rE^{pq}(X;\Z_2)$,
called
\emph{Kalinin's spectral sequence} of~$(X,c_X)$. As above, it is
convenient to start the sequence at the term
$\RH1^*(X;\Z_2)=H^*(X;\Z_2)$ with the differential $\Rd1^*=1+c^*$.

Kalinin's spectral sequence converges to
$H^*(\Fix c_X;\Z_2)$. More precisely, there is an increasing
filtration~$\{\CF_q\}=\{\CF_q(X;\Z_2)\}$ on~$H^*(\Fix c_X;\Z_2)$, called
\emph{Kalinin's filtration}, and homomorphisms
$\bv^q\:\iH^q(X;\Z_2)\to H^*(\Fix c_X;\Z_2)/\CF_{q-1}$, called
\emph{Viro homomorphisms}, which
establish isomorphisms of the graded groups.
In general,  the convergence does {\bf not} respect
the ordinary
grading of $H^*(\Fix c_X;\Z_2)$.

The Smith inequality in~\ref{real.structure} can be derived from
Kalinin's spectral sequence, and $c_X$ is an $M$-involution if and
only if the sequence degenerates at~$\RH1$. If the sequence degenerates
at~$\RH2$, the involution (real variety,
\etc.) is called \emph{$\Z_2$-Galois maximal}.

A similar construction applies to the homology, producing a
$\Z$-graded spectral sequence $\rH_q(X;\Z_2)$ starting from
$H_q(X;\Z_2)$ and converging to $H_*(\Fix c_X;\Z_2)$. The
corresponding decreasing filtration on $H_*(\Fix c_X;\Z_2)$ and
Viro homomorphisms are denoted by $\{\CF^q(X;\Z_2)\}$ and
$\bv_q\:\CF^q\to\iH_q$, respectively.

\paragraph\label{K.mult}
The cup-products in $H^*(X;\Z_2)$ descend to a multiplicative structure
in $\rH^*(X;\Z_2)$, so that $\rH^*$ is a $\Z_2$-algebra and the
differentials~$\rd^*$ are differentiations for all $r\ge2$, \ie,
$\rd^*(x\scup y)=\rd^*x\scup y+x\scup\rd^*y$. The
filtration $\CF_*$ and Viro homomorphisms $\bv^*$ are
multiplicative, \ie, $\CF_p\scup\CF_q\subset\CF_{p+q}$ and
$\bv^*(x\scup y)=\bv^*x\scup\bv^*y$.

\paragraph\label{K.PoincareZ2}
If $X$ is a closed connected $n$-manifold and
$\Fix c_X\ne\varnothing$, Kalinin's spectral sequence inherits
Poincar\'{e} duality: for each $r\le\infty$ one has
$\rH^n(X;\Z_2)=\Z_2$, the cup-product
$\rH^p(X;\Z_2)\otimes\rH^{n-p}(X;\Z_2)\to\Z_2$ is a perfect
pairing, and the differentials $\rd^p$ and $\rd^{n-p-r+1}$,
$1\le r<\infty$, are dual to each other.

The last member~$\CF^n$ of the homological filtration is the
group~$\Z_2$ 
spanned
by the class $w^{-1}(\nu)\scap[X_R]$, where $\nu$ is the normal
bundle of~$X_\R$ in~$X$ and $w(\nu)$ is its total Stiefel-Whitney
class.
(Recall that, if $X$ is a complex
manifold and $c$ is a real structure, the normal bundle~$\nu$ is
canonically isomorphic to the tangent bundle~$\tau$ of~$X_\R$; the
isomorphism is given by the multiplication by~$i$.) Hence, in
terms of the cohomology of $\Fix c$ the Poincar\'{e} duality above
can be stated as follows: {\proclaimfont
the pairing
$(x,y)\mapsto\langle x\scup y\scup w^{-1}(\nu),[X_\R]\rangle\in\Z_2$ is
perfect and, with respect to this pairing, one has
$\CF_{n-q-1}=\CF_q^\perp$.}

\paragraph\label{K.Z}
Now, let~$G=\Z$, and let $h\in H^2(\Rpi;\Z)=\Z_2$ be the
generator.
Assume, as above, that $X$ is a {\sl CW}-complex of finite
dimension. Then the \emph{stabilization homomorphisms}
$$
\scup h\:\rE^{pq}(X;\Z)\to\rE^{p+2,q}(X;\Z),\quad
\scup h\:H^n_c(X;\Z)\to H^{n+2}_c(X;\Z)
$$
are isomorphisms for $p\gg0$ and one has
$$
\gather
\lim_{k\to\infty}H^{2k}_c(X;\Z)
 =H^{2k\gg0}(\Fix c_X\times\Rpi;\Z)
 =H\even(\Fix c_X;\Z_2),\\
\lim_{k\to\infty}H^{2k+1}_c(X;\Z)
 =H^{2k+1\gg0}(\Fix c_X\times\Rpi;\Z)
 =H\odd(\Fix c_X;\Z_2).
\endgather
$$
(We use the notation $H^{p\bmod2}=\bigoplus_{i=p\bmod2}H^i$,
$H\even=H^{0\bmod2}$, $H\odd=H^{1\bmod2}$.)
Thus, one
obtains a $(\Z_2\times\Z)$-graded spectral sequence
$$
\rH^{pq}(X;\Z)=\lim_{k\to\infty}\rE^{2k+p,q}(X;\Z),\quad p\in\Z_2,
$$
which is also called \emph{Kalinin's spectral sequence}
of~$(X,c_X)$ (with coefficients in~$\Z$). It converges to
$H\even(\Fix c_X;\Z_2)\oplus H\odd(\Fix c_X;\Z_2)$, \ie, there are
increasing filtrations $\{\CF^p_q\}=\{\CF^p_q(X;\Z)\}$ on
$H^{p\bmod2}(\Fix c_X;\Z_2)$, $p\in\Z_2$,
called \emph{Kalinin's filtration}, and
homomorphisms
$\bv^{pq}\:\iH^{pq}(X;\Z)\to H^{p\bmod2}(\Fix c_X;\Z_2)/\CF^p_{q-1}$,
called \emph{Viro homomorphisms}, which
establish isomorphisms of the graded groups.

As in~\ref{K.Z2}, one can start the sequence at
the term $\RH1^{pq}(X;\Z)=H^q(X;\Z)$ with differential
$\Rd1^{pq}=1-(-1)^{p+q}c^*$.
If the sequence $\rH^{**}(X;\Z)$ degenerates
at~$\RH2$, the involution (real variety,
\etc.) is called \emph{$\Z$-Galois maximal}.

\paragraph\label{K.mult.Z}
Kalinin's spectral sequence $\rH^{**}(X;\Z)$ is multiplicative (in the same
sense as in~\ref{K.mult}), the multiplicative structure inducing
the product
$$
x\otimes y\mapsto x\scup y+\Sq^1x\scup\Sq^1y
$$
in the limit term $H\even(\Fix c_X;\Z_2)\oplus H\odd(\Fix c_X;\Z_2)$.
(Here $\Sq^1\:H^p(\,\cdot\,;\Z_2)\to H^{p+1}(\,\cdot\,;\Z_2)$
stands for the Bockstein homomorphism.)

\paragraph\label{K.ZtoZ2}
Reduction modulo~$2$ induces a homomorphism
$$
\rH^{0,q}(X;\Z)\oplus\rH^{1,q}(X;\Z)\to\rH^q(X;\Z_2)
$$
of $\Z$-graded spectral sequences, which is compatible with the
isomorphism
$$
H\even(\Fix c_X;\Z_2)\oplus H\odd(\Fix c_X;\Z_2)=H^*(\Fix c_X;\Z_2)
 @>1+\Sq^1>>H^*(\Fix c_X;\Z_2)
$$
of their limit terms. If $H^*(X;\Z)$ is free of $2$-torsion, reduction
modulo~$2$ is an isomorphism starting from the term $\RH2$.

\end


\ifx\defloaded\undefined\input elliptic.def \fi

\begin

\setcounter\section2

\def\pos#1{^{#1}}
\def\neg#1{^{-#1}}

\section{Real elliptic surfaces}\label{res}

In what follows, a \emph{surface} is a {\bf
nonsingular}
complex
manifold of complex dimension two. In the few cases when singular
surfaces are considered, it is specified explicitly.
Proofs of most statements in this section are omitted.
We refer the reader to the excellent founding paper
by K.~Kodaira~\cite{Kod}, or to
the more recent
monographs~\cite{FM}
and~\cite{BPV}.

\subsection{Elliptic surfaces}\label{es}

\subsubsection{}\label{def}
An \emph{elliptic surface} is a surface~$E$ equipped with an
\emph{elliptic fibration}, \ie, a proper holomorphic map $p\:E\to
B$ to a
nonsingular curve~$B$ (called the \emph{base} of the
fibration) such that
for all but finitely many points $b\in B$
the fiber $p^{-1}(b)$ is a
nonsingular curve of genus~$1$.
We will use the notation $p\:E|_U\to U$ (or just $E|_U$) for the
restriction of the fibration to a subset $U\subset B$. The
restriction to the subset~$B^\sharp$ of the regular values of~$p$ is
denoted by $p^\sharp\:E^\sharp\to B^\sharp$. (In other words,
$E^\sharp$ is formed by the nonsingular fibers of~$p$.)

Two fibrations $p\:E\to B$ and $p'\:E\to B'$ on the same surface~$E$
are considered \emph
{identical} if there is an isomorphism $b\:B\to B'$ such that
$p'=p\circ b$. A \emph{morphism} between two elliptic fibrations
$p\:E\to B$ and $p'\:E'\to B'$ is a pair of proper holomorphic
maps $\tilde f\:E\to E'$ and $f\:B\to B'$ such that
$p'\circ\tilde f=f\circ p$. Two fibrations $p\:E\to B$ and
$p'\:E'\to B'$ are \emph{isomorphic} if there is a pair of
bi-holomorphic maps $\tilde f\:E\to E'$ and $f\:B\to B'$ such
that $p'\circ\tilde f=f\circ p$.

A compact surface~$E$ of positive Kodaira dimension~$\kappa(E)$
admits at most one elliptic fibration. All compact surfaces of
Kodaira dimension~$1$ are elliptic; they are called \emph{properly
elliptic}.

An \emph{elementary deformation} of elliptic fibrations consists
of a nonsingular $3$-fold~$X$, a nonsingular surface~$S$, a proper
holomorphic map $p\:X\to S$, and a
deformation (in the sense of Kodaira-Spencer)
$\pi\:S\to D=\{z\in\C\mid\mathopen|z\mathclose|<1\}$ such that
$p\circ\pi$ is a submersion and each restriction $p_t\:X_t\to S_t$
of~$p$ to the slices $S_t=\pi^{-1}(t)$ and $X_t=p^{-1}S_t$, $t\in
D$, is an elliptic fibration. The restrictions~$p_t$ are said to
be \emph{connected by an elementary deformation}.
\emph{Deformation equivalence} of elliptic fibrations is the
equivalence relation generated by isomorphisms and elementary
deformations. 
Any deformation of a
properly elliptic surface~$X_0$ admits a (unique) structure of
deformation of elliptic fibrations.

\paragraph\label{rdef}
An elliptic
fibration $p\:E\to B$ is called \emph{real} if
both
$E$ and $B$
are equipped with real structures
$c_E\:E\to E$ and $c_B\:B\to B$ so that $c_B\circ p=p\circ c_E$.
(When it does not lead to a confusion, we will omit the subscripts
in the notation for the real structure.) If
$E$ is compact and
$\kappa(E)>0$,
then, due to the uniqueness of the elliptic
fibration, any real structure $c\:E\to E$ descends to~$B$.

The notions of morphism, isomorphism, deformation, \etc\ extend to
the real case in the obvious way: one requires that all manifolds
involved should be equipped with real structures that are
respected by all maps. (For elementary deformations, the real
structure on the unit disk $D\subset\C$ is that induced from the
complex conjugation.)

\paragraph\label{relmin}
In this paper, we only consider \emph{relatively minimal} elliptic
fibrations, \ie, those without $(-1)$-curves in the fibers. For a
compact elliptic surface~$E$ this is equivalent to the condition
$K_E^2=0$. As is known, a fiber of an elliptic fibration (as well
as any fibration whose generic fiber is a curve of positive genus)
can not contain intersecting $(-1)$-curves. This implies that each
elliptic fibration admits a unique relatively minimal model and,
in particular, the relatively minimal model is real whenever the
original fibration is. Moreover, by Kodaira's results on the
stability of exceptional curves, the deformation study of elliptic
fibrations (both complex and real) is reduced to the deformation
study of their relatively minimal models. In particular, the
elliptic fibrations deformation equivalent to a relatively minimal
elliptic fibration are relatively minimal.

Any relatively minimal elliptic fibration~$p:E\to B$ is
\emph{strongly relatively minimal}, \ie, all fiber-to-fiber
bi-meromorphic maps $E\to E$ (and, hence, fiber-to-fiber
bi-antimeromorphic maps $E\to E$) are regular. In particular, the
fibration is uniquely determined by its restriction
$p^\sharp\:E^\sharp\to B^\sharp$ (see \ref{def}), and any real
structure on $p^\sharp\:E^\sharp\to B^\sharp$ extends uniquely to
a real structure on $p\:E\to B$.

\subsection{Jacobian surfaces}\label{ws}
From now on, we consider only relatively minimal elliptic
fibrations {\bf without muptiple fibers}.

\paragraph\label{inv}
To each elliptic fibration $p\:E\to B$ one can associate its
\emph{functional} (or $j$-) \emph{invariant} $j\:B\to\Cp1$ and its
\emph{homological invariant} $\RR1p\Z_E$.

The functional invariant is the extension to~$B$ of the
meromorphic function $B^\sharp\to\C$ sending each nonsingular
fiber to its $j$-invariant; following Kodaira, we divide the
$j$-invariant by $12^3$, so that its `special' values are $j=0$
and~$1$: a nonsingular elliptic curve~$C$ with $j(C)=0$ or~$1$ has
a complex multiplication of order~$6$ or~$4$, respectively. Since
reversing the complex structure on a nonsingular elliptic
curve~$C$ transforms $j(C)$ to $\overline{j(C)}$, the functional
invariant of a real elliptic fibration is real, $j\circ
c_B=\bar\jmath$. (When speaking about a real structure on the
functional invariant $j\:B\to\Cp1$ we always assume that the real
structure on the Riemann sphere $\Cp1=\C\cup\{\infty\}$ is
standard, so that the points~$0$, $1$, and~$\infty$ are real.)

An elliptic fibration with $j=\const$ (respectively, $j\ne\const$)
is called \emph{isotrivial} (respectively,
\emph{non-isotrivial}). In this paper we deal mainly with
non-isotrivial fibrations. Note that, unless $j=0$ or~$1$, an
isotrivial fibration has no singular fibers.

Non-isotrivial fibrations have a \emph{strong extension
property}: given
a nonsingular curve~$B$, a point $b_0\in B$, and an elliptic
fibration over $B\sminus\{b_0\}$,
there is a
unique (relatively minimal) elliptic fibration
over~$B$ whose restriction to $B\sminus\{b_0\}$ is the given one.

\paragraph\label{hominv}
The homological invariant (see \ref{inv}) is often defined as the
monodromy in the $1$-cohomology of the nonsingular fibers, \ie, as
a local system~$\CM$ on~$B^\sharp$ with fiber $\Z\oplus\Z$, and
as such it is just the restriction of $\RR1p\Z_E$ to~$B^\sharp$.
Then $\RR1p\Z_E=i_*\CM$, $i\:B^\sharp\to B$ standing for the
inclusion. The homological invariant of a real elliptic fibration
inherits a real structure from that on~$\Z_E$.

The homological invariant of an elliptic fibration is closely
related to its $j$-invariant. Representing
(by means of the elliptic modular function
$j(z)=j(\C/L_z)$ where $L_z=\Z +z\cdot\Z$) the complex line
$\Cp1\smallsetminus\{\infty\}$ as the quotient of the upper half
plane by the modular group
$\operatorname{\text{\sl PSL}}(2,\Z)=
 \operatorname{\text{\sl SL}}(2,\Z)/\{\pm1\}$, one
equips the space $\Cp1\smallsetminus\{0,1,\infty\}$ with
a $\operatorname{\text{\sl
PSL}}(2,\Z)$-principal bundle $\Cal P$.
(Here
$0$ and~$1$ are the images of the
two
unstable points
of the action.)
A local system~$\CM$ as above
is said to \emph{belong to} a holomorphic map $j\:B\to\Cp1$ if
the principal
$\operatorname{\text{\sl PSL}}(2,\Z)$-bundle
associated with~$\CM$ is $j^*\Cal P$.
The homological invariant of an
elliptic fibration belongs to its functional invariant.

Two real
structures on a holomorphic map $j\:B\to\Cp1$ and a
$(\Z\times\Z)$-local system~$\CM$ belonging to~$j$ are called
\emph{concordant} if
they are both
lifts of the same real structure $c_B$ on $B$.
This is the case if $j$ and~$\CM$ are,
respectively, the functional and homological invariant of a real
elliptic fibration.

In both complex and real cases,
the passage from~$j$ to~$\CM$ involves a choice of one of the two
lifts over each loop $\gamma\subset B^\sharp$; the lifts
differ by the multiplication by $-\id\in\SL(2,\Z)$.
Next statement asserts that the elliptic fibrations obtained
from distinct lifts differ topologically.

\lemma\label{Aconj-A}
A matrix $A\in\SL(2,\Z)$ is never conjugate to~$-A$.
\endlemma

\proof
In fact, the statement holds for the bigger group
$\SL(2,\R)$. If a $(2\times2)$-matrix
$A$ is similar to~$-A$
and $\det A=1$, one can easily see that the eigenvalues of~$A$
must be $\pm i$, \ie, $A$ is a complex structure
(or the rotation through $\pm\pi/2$). Then $-A$ is the rotation in
the opposite direction; it is
not conjugate to~$A$ by an orientation preserving transformation.
\endproof

\paragraph\label{jac}
Among all elliptic fibrations with
given functional and homological invariants there is a unique,
up to isomorphism,
elliptic fibration
$p\:J\to B$
with a section. We equip~$J$
with a distinguished
section $s\:B\to J$ and
call the pair $(J,s)$ the \emph{Jacobian elliptic fibration}
associated to $p\:E\to B$ (or to the given pair of functional and
homological invariants).
The group of automorphisms of~$J$ preserving~$s$ is a cyclic group
of order~$2$, $4$, or~$6$, the last two cases occurring only if
$j=\const$. In all cases the element
of order~$2$ acts in
each nonsingular fiber as the multiplication by~$(-1)$.

An elementary deformation\label{jacobdeform}
$$
p\:X\to S,\quad  \pi\: S\to D
$$
of elliptic fibrations is called
\emph{Jacobian},
or
a \emph{deformation through Jacobian fibrations}, if it is equipped
with a section $s\:D\to S$ of $p\circ \pi$. This notion extends to
the real case in the usual way: one requires that the deformation
and the section should be real.

In order to construct the Jacobian elliptic fibration $J=J(E)$
associated to $p\:E\to B$, one can start from
$p^\sharp\:E^\sharp\to B^\sharp$ and replace each
fiber~$F_b=p^{-1}(b)$ by its Jacobian $\Pic^0(F_b)$. Then it
remains to complete the elliptic fibration
$\Pic^{0}_{B^\sharp}=\bigcup_{b\in B^\sharp}\Pic^0(F_b)\to
B^\sharp$ and to take its relatively minimal model.
(Note that the
completion step requires, in fact, a thorough understanding of
singular fibers and their local models. It turns out that the
Jacobian fibration has singular fibers of the same types as the
original one.) The strong relative minimality implies uniqueness.
Moreover, it shows that the construction is functorial, \ie, any
(anti-)isomorphism $E\to E'$ of elliptic fibrations induces an
(anti-)isomorphism $J(E)\to J(E')$ of their Jacobians respecting
the distinguished sections. In particular, the Jacobian elliptic
fibration associated to a real elliptic fibration inherits an
\emph{associated Jacobian real structure}, namely, the structure
$c_J\:J(E)\to J(E)$
induced by the action of $c_E$ on the Jacobians of the nonsingular
fibers. Unless $j=\const$, the only other real structure
preserving the section is $-c_J$, \ie,
the composition of~$c_J$
and the fiberwise multiplication by~$(-1)$.
The real structures~$c_J$ and $-c_J$ are called \emph{opposite} to
each other.

A deformation of elliptic surfaces
gives rise to a
natural
deformation of the associated Jacobian surfaces.
In particular, if the
original deformation is real,
so is the resulting Jacobian deformation.

\paragraph\label{counthominv}
As it follows from the strong extension property (see \ref{inv}),
for each pair $(j,\CM)$ consisting of a nonconstant holomorphic
map $j\:B\to\Cp1$ and belonging to it $(\Z\times\Z)$-local
system~$\CM$ on $B^\sharp=j^{-1}(\{0,1,\infty\})$ there exist a
Jacobian fibration $J(j,\CM)$
whose functional and homological
invariants are $j$ and $\CM$. The set of local systems belonging
to a given map~$j$ is an affine space over $H^1(B^\sharp;\Z_2)$.
In particular, their number and, thus, the number of Jacobian
fibrations with given functional invariant is $2^r$, where
$r={b_1(B^\sharp)}$. We recall that the monodromy
of $\CM$ along a small
loop about a point $b_0\in B\sminus B^\sharp$ determines and is
determined by the topology of the singular fiber at~$b_0$. Thus,
if the types of the singular fibers are fixed (\eg, under the
assumption that the fibration is \emph{almost generic}, \ie, does not have
singular fibers other than those of type~$\I_1$), the admissible
homological invariants form an affine space over $H^1(B;\Z_2)$.
(Note, though, that the latter assumption imposes certain
restrictions to the functional invariant, see~\ref{j.generic} below.)


\paragraph\label{realhominv}
Let $j$ and $\CM$ be as in~\ref{counthominv}.
Assume that $B_\R \ne \varnothing$. If~$j$ and~$\CM$ are equipped
with concordant real structures, then, as it follows from the
strong relative minimality and the local Torelli theorem for
elliptic curves with a marked point, the real structures lift to a
real structure on $J(j,\CM)$
that makes it a real Jacobian fibration
with given real invariants~$j$ and~$\CM$.
When nonempty,
the set of real local systems
belonging to and concordant with a given real map $j\:B\to\Cp1$ is
an affine space over $\HH^1(B^\sharp;\Z_2)$ (see~\ref{HHgen};
the existence question is discussed in~\ref{tg->realinv}).
Therefore, their number and, thus, the number of real Jacobian
fibrations with a given real functional invariant~$j$ is equal to
$2^r$, where
$r=\dim(H^1(B^\sharp;\Z_2))^c+1$,
see~\ref{HHcount}. (Note that the term
$1=\dim H^1(\Z_2;H^0(B^\sharp;\Z_2))$ accounts for the two
distinct real structures on a given complex Jacobian fibration,
see~\ref{jac}, and $2^{r-1}$ is the number of local systems
admitting a concordant real structure.) As in~\ref{counthominv},
if the fibration is assumed almost generic, this number reduces to $2^s$,
where $s=\dim(H^1(B;\Z_2))^c+1$.

\paragraph\label{sheafs}
Removing from a Jacobian fibration~$J$ the singular points
(including multiple components) of its singular fibers, one
obtains an analytic family $J\ab\to B$ of (not necessarily
connected) abelian Lie groups
with~$s$ as the zero section. Any
section of~$J$ is contained in $J\ab$, and any two sections differ
by a \emph{translation}, \ie, an automorphism of $p\:J\to B$ which
is a translation in each nonsingular fiber. Furthermore, one can
introduce two sheafs of abelian groups: the sheaf~$\tJac$ of germs
of holomorphic sections of~$J\ab$ and the sheaf~$\Jac$ of germs of
sections of~$J\ab$ intersecting each fiber at the same component
as~$s$. The two sheafs differ by a skyscraper~$\CS$ having finite
fibers and concentrated at the points of~$B$ corresponding to
reducible singular fibers with at least two simple components,
$$
0@>>>\Jac@>>>\tilde\Jac@>>>\CS@>>>0.\eqtag\label{stJac}
$$
(In fact, the order of the stalk $\CS_b$ at a point $b\in B$ is
exactly the number of simple components in the fiber $p^{-1}(b)$.)
From the exponential sequence over~$J$ one can also obtain the
following short exact sequence for~$\Jac$:
$$
0@>\phantom{2\pi i}>>\RR1p\Z_J 
@>2\pi i>>\RR1p\CO_J 
 @>\phantom{2\pi i}>>\Jac@>\phantom{2\pi i}>>0.\eqtag\label{sJac}
$$
Recall that, for any elliptic fibration $p\:E\to B$, the sheaf
$\RR1p\CO_E$ is of the form
$\CO_B(L\neg1)$, where
$L$ is a certain line bundle on~$B$ determined solely by the
functional and homological
invariants of the fibration. If the fibration has a
section $s\:B\to E$, then
$L\neg1$
is the normal bundle of~$s$; its
degree (\ie, the self-intersection of~$s$) is negative unless
$j=\const$.

\paragraph\label{realsheafs}
Let $p\:E\to B$ be an elliptic fibration and $J\to B$ its
Jacobian. Any \hbox{(anti-)}\penalty0automorphism $g\:E\to E$
induces an \hbox{(anti-)}\penalty0automorphism $J(g)\:J\to J$
preserving the distinguished section (see \ref{jac}). The kernel
$\Aut_0E$ of the map $g\mapsto J(g)$ is formed by the translations
of~$E$. Obviously, there is a canonical isomorphism
$\Aut_0E=\Aut_0J$ and both of the groups are isomorphic to the
\emph{Mordel-Weil group} $\Gamma(B;\tJac)$ (the group of sections
of $J\to B$). Furthermore, if $E$ has a section itself, it is
isomorphic to~$J$, and the set of isomorphisms $\Gf\:E\to J$
identical on the Jacobian is an affine space
over~$\Gamma(B;\tJac)$: one has $\Gf+t=\Gf\circ t_E=t_J\circ\Gf$,
where $t\in\Gamma(B;\tJac)$ is a section and $t_E$, $t_J$ are the
corresponding translations of~$E$ and~$J$, respectively.

Recall that the Mordel-Weil group of any compact non-isotrivial
elliptic fibration is discrete. This follows, \eg, from the fact
that the line bundle~$L$ has no
sections (as it has negative degree).

Now, let $c_E$ be a real structure on~$E$ and $c_J$ the Jacobian
real structure on~$J$.
The latter induces a real structure $c$
on~$\tJac$, $s\mapsto c_J\circ s\circ c_B$, \cf.~\ref{rsheaf},
which is compatible with~\eqref{stJac} and~\eqref{sJac}.
The set
of all real structures on~$E$ whose Jacobian is~$c_J$ is
an affine
space over the subgroup $\Ker(1+c)\subset\Gamma(B;\tJac)$, the
affine
action
being $c_E+t=t_E\circ c_E$.
(Here,
the condition
$(1+c)t=0$ is necessary and sufficient for
the composition $t_E\circ c_E$ to be an involution.) The
shift~$t_E$ by a section $t\in\Gamma(B;\tJac)$ is real (\ie,
commutes with~$c_E$) if and only if $(1-c)t=0$. Note that the
condition does not depend on~$c_E$; thus, $t_E$ commutes with any
real structure whose Jacobian is~$c_J$. More generally, for any
such real structure~$c_E$ on~$E$ one has $t_E^{-1}\circ c_E\circ
t_E=c_E-(1-c)t$. In particular, the set of all real structures
on~$E$ with the given Jacobian~$c_J$ is an affine space over
$H^1(\Z_2;\Gamma(B;\tJac))$.



\subsection{Trigonal curves and Weierstra{\ss} models}\label{trigonal-Weierstrass}
Traditionally,
elliptic curves with a rational point
are described \via\ the so called Weierstra{\ss} equation.
In the case of elliptic surfaces, this approach leads
to trigonal curves on ruled surfaces.

\subsubsection{Trigonal curves}\label{trigonal-curves}
Let $q\:\Sigma\to B$ be a geometrically ruled surface with a
distinguished section~$s$. A \emph{trigonal curve} on~$\Sigma$ is
a reduced curve $C\subset\Sigma$ disjoint from~$s$ and such that the
restriction $q\:C\to B$ is of degree three. Given a trigonal curve
$C\subset\Sigma$, the fiberwise center of gravity of the three
points of~$C$ (regarded as points in the affine fiber of
$\Sigma\sminus s$)
defines an additional section~$0$ of~$\Sigma$; thus,
the $2$-bundle
whose projectivization is~$\Sigma$ splits and, after a renormalization,
can be chosen in the form $1\oplus Y$. We choose the normalization
so that the projectivization of the $Y$ summand is the zero
section.

Any trigonal curve can be given by a \emph{Weierstra{\ss} equation};
in appropriate affine charts it has the form
$$
x^3+\p\, x+\q=0,\eqtag\label{trigeq}
$$
where $g_2$ and~$g_3$ are certain sections of~$Y^2$ and~$Y^3$,
respectively, and $x$ is a coordinate such that $x=0$ is the zero
section and $x=\infty$ is the distinguished section~$s$. The
sections~$g_2$, $g_3$ are determined by the curve uniquely up to
the transformation
$$
(\p,\q)\mapsto (t^2\p, t^3\q),\quad t\in\Gamma(B,\CO_B^*).
$$

If both $(\Sigma,s)$ and~$C$ are real, then $Y$ is a real line bundle
and the sections $g_2$, $g_3$ can also be chosen real; they are
defined uniquely up to the above transformation with a real
section~$t$.

The \emph{$j$-invariant} of a trigonal curve $C\subset\Sigma$
is the function
$j\:B\to\Cp1$ given by
$$
j=\frac{4\p^3}\Delta,\quad \Delta=4\p^3+27\q^2.\eqtag\label{jinv}
$$
Geometrically, the value of~$j$ at a generic point $b\in B$ is the
usual $j$-invariant of the quadruple of points cut by the union
$C\cup s$ in the projective line $q^{-1}(b)$.

As the equation suggests, in general the $j$-invariant does not
change continuously when the curve is deformed; even the degree
of~$j$ can change.

From now on, by a \emph{deformation} of a trigonal curve
$C\subset\Sigma$ we mean a deformation of the quadruple
$(B,q,s,C)$ (\ie, neither~$\Sigma$ nor $B$ are assumed fixed). As
usual, the \emph{deformation equivalence} of trigonal curves is
the equivalence relation generated by deformations and
isomorphisms (of the quadruples as above).

A trigonal curve $C \subset \Sigma$ is called \emph{almost generic}
if it is nonsingular and has no vertical flexes.
If this is the case,
the $j$-invariant $j\:B\to\Cp1$ has
degree $\deg j=6d$, where $d=\deg Y$,
the point $\infty\in\Cp1$
is a regular value of~$j$, its $6d$ pull-backs corresponding to
the vertical tangents of the curve,
and all pull-backs of~$0$ and~$1$ have ramification
index $0\bmod3$ (respectively, $0\bmod2$).
By an arbitrary small deformation (including a change
of the complex structure of the base)
one can achieve
that the $j$-invariant have so called
\emph{generic branching behavior} (which, in fact, is highly
non-generic for a function $B\to\Cp1$):
in addition to the above conditions one requires that the
ramification index of each pull-back of~$0$ (respectively,~$1$)
should be
exactly~$3$ (respectively,~$2$).
(Note that other critical values, which $j$
unavoidably has, are irrelevant.) A trigonal curve whose
$j$-invariant has generic branching behavior is called
\emph{generic}.

\subsubsection{Topology}\label{topology}
The real part~$\Sigma_\R$ of a real geometrically ruled surface
$q\:\Sigma\to B$ consists of
several connected components~$\Sigma_i$, one over each component~$B_i$
of~$B_\R$.
Each~$\Sigma_i$ is either a torus or a Klein bottle.
If $\Sigma$ has the form $\Cp{}(1\oplus Y)$ as above, then a
component $\Sigma_i$ is orientable if and only if the
restriction~$Y_i$ of the real part~$Y_\R$ of~$Y$ to the
corresponding component~$B_i$ is topologically trivial. In other
words, if $Y$ is given by a real divisor~$D$, then a
component~$\Sigma_i$ is orientable if and only if the degree
$\deg(D\cap B_i)$ is even.
The latter remark shows also that $\sum_i\deg Y_i=\deg Y\bmod2$.
Hence, $\Sigma_\R$ is necessarily nonorientable whenever $\deg Y$ is
odd.

Let $C\subset\Sigma$ be a real trigonal curve, and let
$q_\R\:C_\R\to B_\R$ be the projection. The real part $C_\R$
splits into groups of components $C_i=q_\R^{-1}(B_i)$. Each
restriction $q_\R\:C_i\to B_i$ is onto. A component~$B_i$ of~$B_\R$
(and the corresponding group~$C_i$)
is called \emph{hyperbolic} if the restriction $q_\R\:C_i\to B_i$
is generically three-to-one; otherwise,
it is called
\emph{non-hyperbolic}.
A curve~$C$ with non-empty real part
is called \emph{hyperbolic} if all its groups $C_i$
are hyperbolic; otherwise, it is called \emph{non-hyperbolic}.

Over a hyperbolic
component~$B_i$, the group~$C_i$ of an almost generic curve
consists of a `central' component which projects to~$B_i$
homeomorphically and two (if $\Sigma_i$ is orientable) or one (if
$\Sigma_i$ is nonorientable) additional components; the
restriction of~$q_\R$ to the union of the additional components is
a double covering, trivial in the former case and nontrivial in
the latter case.

\midinsert
\eps{S3Fig1}
\figure\label{fig.topology}
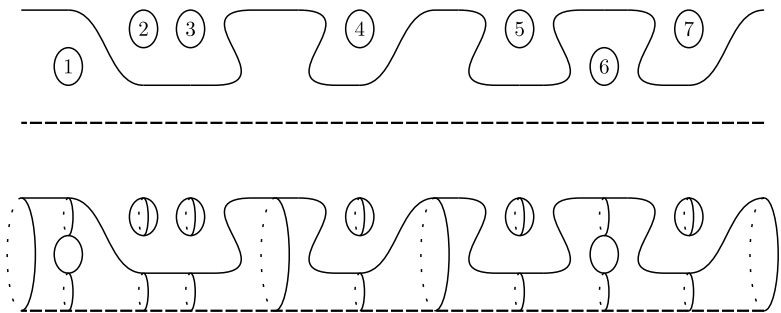
A typical non-hyperbolic trigonal curve (top) and a corresponding Jacobian
surface (bottom); the horizontal dotted lines represent
the distinguished sections~$s$ of the
surfaces.
\endfigure
\endinsert

The group~$C_i$ of an almost generic curve over a non-hyperbolic
component~$B_i$ looks as shown in Figure~\ref{fig.topology}. More
precisely, $C_i$ has a `long' component mapped onto $B_i$
and several contractible components, commonly called ovals;
the long component may contain a few `zigzags',
which are also preserved by fiberwise isotopies.
For the purpose of this paper,
we define \emph{ovals} and \emph{zigzags}
as the connected components
of the set $\{b \in B_i | \#q^{-1}_\R \ge 2\}$; ovals are those
whose pull-back is disconnected. The set of all ovals
within a non-hyperbolic component~$B_i$ inherits from~$B_i$
a pair of opposite cyclic orders.

Pick a non-hyperbolic component~$B_i$.
Fix a section $\sigma\:B_i\to\Sigma_i$ disjoint from~$s$
and taking each oval inside the corresponding contractible
component;
such a section is unique up to homotopy.
A set of consecutive ovals in $B_i$
is called a \emph{chain} if between any two neighboring
ovals of the set the section $\sigma$ intersects the long component
an even number of times. For example, in Figure~\ref{fig.topology}
the maximal chains are $[1]$, $[2,3,4,5]$, $[6]$, and $[7]$
(assuming that $\Sigma_i$ is orientable; otherwise,
ovals 7 and 1 form a single chain).
A chain of ovals is called \emph{complete}
if it contains all ovals in a single component~$B_i$.

The notions of oval and chain extend to all nonsingular
trigonal curves. Note that a non-hyperbolic nonsingular curve cannot
be isotrivial; hence, it can be perturbed to an almost generic
one.

\subsubsection{Weierstra{\ss} models}\label{wmodels}
The \emph{Weierstra{\ss} model} of a Jacobian elliptic surface
$p\:J\to B$
is obtained
from~$J$ by contracting all components of the fibers of~$p$ that
do not intersect the distinguished section~$s$.
The result is a
proper map $p\:J\w\to B$, where $J\w$ has at worst
simple singular points, and a section $s\:B\to J\w$ not
passing through the singular points of~$J\w$.
The original Jacobian surface~$J$ is recovered from~$J\w$ by
resolving its singularities.

The
quotient of~$J\w$ by the fiberwise multiplication by~$(-1)$ is
a geometrically ruled surface~$\Sigma$ over~$B$, the section~$s$
mapping to a section of~$\Sigma$. The projection
$J\w\to\Sigma$ is the double covering defined by the (fiberwise)
linear system
$\mathopen|2s\mathclose|$ on~$J\w$; its branch curve is the
disjoint union of the exceptional section~$s$
and a certain trigonal curve~$C$ on~$\Sigma$.
In particular, $\Sigma$ has the form $\Cp{}(1\oplus Y)$ and
$Y=L^2$,
where $L$ is the conormal bundle of~$s$ in~$J$
(\cf.~\ref{sheafs}).

The sections~$\p$, $\q$ defining~$C$, see~\eqref{trigeq},
must satisfy the following
conditions:
\roster
\item\local{WDelta}
the discriminant $\Delta=4\p^3+27\q^2$ is not identically zero,
and
\item\local{Worder}
at each point $b\in B$ one has $\min(3\ord_b(\p),2\ord_b(\q))<12$.
\endroster
(The former condition ensures that generic fibers are
nonsingular elliptic curves, and the latter implies that all
singular points of~$J\w$ are simple.)
Conversely, given a ruled surface $\Sigma=\Cp{}(1\oplus Y)$ with a
section~$s$ and
a trigonal curve~\eqref{trigeq},
a choice of a square root~$L$ of~$Y$ defines a unique double
covering of~$\Sigma$ ramified at~$s$ and the curve;
if the pair $(\p,\q)$ in~\eqref{trigeq}
satisfies~\loccit{WDelta} and~\loccit{Worder}
above, the double covering is the Weierstra{\ss} model
of a Jacobian elliptic surface. The $j$-invariant of the resulting
surface
is given
by~\eqref{jinv}.

\subsubsection{Real structures}\label{real-part}
The construction above is natural. Hence, if
the Jacobian surface $(J,s)$ is real, so are
the ruled surface~$\Sigma$, its section~$s$, and the branch
curve~$C$. Conversely, if the surface $(\Sigma,s)$, curve~$C$,
and square root~$L$ are real, then the Jacobian surface~$J$ resulting
from the construction above inherits two opposite real structures.
Under the assumptions, $\Sigma_\R$
splits into two \emph{halves} with common boundary
$C_\R \cup s_\R$, the halves being the projections of the real parts
of the two real structures on~$J$.
A choice of one of the two real structures
is equivalent to a choice of one of the two halves.

The real part~$J_\R$ is a double of the corresponding half; it
looks as shown in Figure~\ref{fig.topology}.
It splits into groups of components
$J_i=p_\R^{-1}(B_i)$.
Each group~$J_i$ has a distinguished component that contains the
section~$s$ over~$B_i$; we call this component \emph{principal}.
Its orientability is governed by the restriction~$L_i$ of the real
part~$L_\R$ of~$L$ to~$B_i$: the principal component is orientable
if and only if $L_i$ is topologically trivial.
Besides, there are a few \emph{extra} components disjoint from the
section~$s_\R$. In the
non-hyperbolic case all extra components are spheres.
In the hyperbolic case, there is exactly one extra component,
which is either a torus or a Klein bottle, depending on whether
$L_i$ is trivial or not. Thus, in the hyperbolic case the two
components are either both orientable or both not.

The following lemma states that a real Jacobian surface and its
branch curve have the same discrepancy.

\lemma\label{M-d}
An almost generic Jacobian surface~$J$
is an $(M-d)$-variety if and only if so is
the trigonal part of
the branch curve of the Weierstra{\ss} model of~$J$.
\endlemma

\proof
Indeed, the isomorphism
$\pi_1(J)\to\pi_1(B)$, see \cite{FM, Proposition 2.2.1}, gives
$\dim H_1(J;\Z_2)=\dim H_1(B;\Z_2)$; the other Betti numbers are
controlled using the Riemann-Hurwitz formula
$\chi(J)=2\chi(\Sigma)-\chi(C)=4\chi(B)-\chi(C)$
and Poincar\'{e} duality. The Betti numbers of the real part~$J_\R$
are found using the description of its topology given
in~\ref{real-part}, and the statement follows from a simple
comparison.
\endproof

\subsubsection{The homological invariant}\label{monodromy}
The Weierstra{\ss} model gives a clear geometric interpretation of
the $\PSL(2,\Z)$-bundle $j^*\Cal P$ defined
by~$j$, see~\ref{hominv}. Indeed,
the modular group
$\PSL(2,\Z)$
is naturally identified with the factorized braid group
$B_3/\Delta^2$, which, in turn, can be regarded as the mapping
class group of the triad $(F;s\cap F,C\cap F)$, where $F$ is a generic
fiber of the ruling. Thus, $j^*\Cal P$ is merely the monodromy
$\pi_1(B^\sharp)\to B_3/\Delta^2$ of the trigonal curve.

As explained in~\ref{counthominv}, the homological invariants
belonging to a given functional invariant~$j$ form an affine space over
$H^1(B^\sharp,\Z_2)$. The branch curve~$C$ narrows this choice
down to $H^1(B,\Z_2)$, as its singularities and vertical tangents
determine the singular
fibers of the covering elliptic surfaces. (Roughly speaking,
the singularities of~$C$ are encoded, in addition to~$j$,
in the presence and
multiplicities of the common roots
of the sections~$g_2$, $g_3$).
The rest of the homological invariant
is recovered \via\ the choice of the square root~$L$ of~$Y$.

In the real case, the choice is narrowed down to $\HH^1(B;\Z_2)$,
see~\ref{realhominv}, and partially it can be made canonical
using the
correspondence $L\mapsto\bigoplus_iw_1(L_i)$, which
is an affine map from the set of
homological invariants onto the subset
$$
\bigl\{\alpha\in H^1(B_\R;\Z_2)\bigm|\alpha[B_\R]=\deg L\bmod2\bigr\}.
$$
In other words, the real Jacobian elliptic surfaces with a given branch
curve are partially distinguished
by the orientability of their principal components.

\subsubsection{Existence of the roots}\label{tg->realinv}
Obviously, the necessary and sufficient condition for the
existence of a square root~$L$ of a line bundle~$Y$ on~$B$ (and,
hence, the existence of an elliptic surface over a given ruled
surface)
is that
$\deg Y$ should be even.

In the real case, if $B_\R\ne\varnothing$, for the existence of a
real square root~$L$ of a real line bundle~$Y$ one should require,
in addition, that the real part~$Y_\R$ is topologically trivial.
Indeed, the condition is obviously necessary. For the sufficiency
notice that, if $B_\R\ne\varnothing$, a line bundle is real if and
only if its class in $\Pic B$ is fixed by
the induced real structure. In particular, the property to have
real roots is invariant under equivariant deformations of~$Y$. On the
other hand,
the real part $\Pic^0_\R B$ has
$2^{b_0(B_\R)-1}$ connected components which are distinguished by
the restrictions to the components~$B_i$
(as the real structure on
$\Pic^0B$ is essentially the same as the induced involution in
$H^1(B;\Z)$\,). Hence, all bundles of a given degree whose real
part is trivial are deformation equivalent.

Combining the above statement with the beginning
of~\ref{topology}, one arrives at the following criteria.

\corollary\label{real-exist}
Assume that a real ruled surface $(\Sigma, s)$, $\Sigma_\R
\ne\varnothing$,
and a real trigonal curve $C \subset \Sigma$ do define
a Jacobian surface. The latter can be chosen real
if
and only if $\Sigma_\R$ is orientable.
\qed
\endcorollary

\subsubsection{An application: generic surfaces}\label{j.generic}
Any non-isotrivial Jacobian surface
can be deformed through
Jacobian surfaces
to an almost generic one.
(In general, that would change the base of
the fibration. The deformation can be
chosen
elementary and
arbitrary small.)
The $j$-invariant of an almost generic surface is similar to that
of an almost generic curve, the pull-backs of~$\infty$, $0$,
and~$1$
corresponding to the singular fibers (of type~$\I_1$), fibers with
complex multiplication of order~$6$, and those with complex
multiplication of order~$4$, respectively.
By another arbitrary small deformation one can
make the surface \emph{generic}, \ie,
achieve
that the $j$-invariant have
generic branching behavior.

If $(J,s)$ is real, the above deformations can also be chosen
real.
For proof one can use the same `cut-and-paste' arguments
as in the complex case, constructing local deformations and
patching them together.
To construct a local perturbation, one can use the Weierstra{\ss}
equation~\eqref{trigeq}, which contains a versal
deformation of the special point (\ie, singular point, vertical
flex, or multiple root of~$g_2$ or~$g_3$)
of the Weierstra{\ss} model
of~$J$. After a covering parameter change, any
perturbation of~$J$ admits a simultaneous resolution of
singularities to which one can extend any real structure and any
automorphism of the original perturbation. (This follows, \eg,
from the Grothendieck-Brieskorn model.) Due to the versality, the
covering deformation can be realized as a deformation of Weierstra{\ss}
models.

If an elliptic fibration
is deformed through almost generic ones,
its functional invariant does change
continuously. Conversely, if an analytic
family of degree~$12d$ functions having generic branching behavior
includes the $j$-invariant of a Jacobian elliptic fibration~$J$,
it results in a unique deformation of~$J$
through
generic Jacobian
fibrations. Observing that
Kodaira's proof
respects real structures (or
using the uniqueness of the deformation),
one obtains a real version of
the statement:
if the fibration~$J$ and the family of functions
are real, the resulting deformation is real.

\end


\ifx\defloaded\undefined\input elliptic.def \fi

\def\extra{^{\fam0 ext}}

\begin

\setcounter\section3

\section{Real Tate-Shafarevich group}\label{Tate-Shafarevich}

Fix a real Jacobian fibration $p\:J\to B$, not necessarily
compact, with a real section $s\:B\to J$. The \emph{real
\rom(analytic\rom) Tate-Shafarevich group} of~$J$ is the set
$\Rsha=\Rsha(J)$ of isomorphism classes of triples $(E,c,\Gf)$,
where $(E,c)$ is a real elliptic surface (over~$B$) without
multiple fibers and $\Gf\:J(E)\to J$ is a real isomorphism. (The
group structure on $\Rsha(J)$ is given by Theorem~\ref{Rsha}
below.) Our principal result in this section is the fact that the
`discrete part' $\Rsha\top=\Rsha/\Rsha^0$ (where $\Rsha^0$ is the
component of unity) is a topological invariant of the pair
$(p,c_J)$.

\subsection{Topological invariance}
Let $p\:J\to B$ be as above and let~$c$ be the canonical real structure
on the sheaf~$\tJac$ of germs of holomorphic sections of~$p$,
see~\ref{realsheafs}.

\theorem\label{Rsha}
There is a natural isomorphism $\Rsha(J)=\HH^1(\tJac,c)$.
\endtheorem

\proof
We mimic the standard proof of the similar result for the complex
(analytic) Tate-Shafarevich group~$\sha(J)$. Pick a triple
$(E,c,\Gf)\in\Rsha(J)$ and use~$\Gf$ to identify the
Jacobian of~$E$ and~$J$. Since $E$ has no multiple fibers, one can
cover~$B$ by $c_B$-invariant open sets~$U_i$ so that each
restriction $E_i=E|_{U_i}$ has a section (not
necessarily real), or, equivalently, there is an isomorphism
$\Gf_i\:E_i\to J_i=J|_{U_i}$.
Let $c_i=\Gf_i^{-1}\circ c_J\circ\Gf_i$ be the real structure
on~$E_i$ induced by~$\Gf_i$. Then the restriction of~$c_E$ to $E_i$ has
the form $c_i+s_i$ for some section
$s_i\in\Gamma(U_i;\tJac)$ satisfying $(1+c)s_i=0$, see~\ref{realsheafs}.
The restrictions of~$\Gf_i$ and~$\Gf^{-1}$ to the intersection
$U_{ij}=U_i\cap U_j$ have the same Jacobian and, hence, differ by
a section $t_{ij}\in\Gamma(U_{ij};\tJac)$: one has
$\Gf_j=\Gf_i+t_{ij}$ (see~\ref{realsheafs} again) and, as usual,
the $1$-cochain $(t_{ij})$ must be a cocycle in the \v{C}ech
complex $(\check C^*_{\CU}(\tJac),d_2)$.
Finally, the real structures on $c_i+s_i$
and~$c_j+s_j$ must coincide on $E|_{U_{ij}}$ and, since
$c_j=c_i-(1-c)t_{ij}$ over~$U_{ij}$, one has
$s_j-s_i=(1-c)t_{ij}$. Thus, the
sections $(s_i,t_{ij})$ form a $1$-cocycle in the \v{C}ech
bi-complex $\check C^{**}_{\CU}(\tJac)$ corresponding to the
covering $\CU=\{U_i\}$. Any other set of isomorphisms $\Gf'_i\:E_i\to J_i$
differs from~$\Gf_i$ by sections $r_i\in\Gamma(U_i;\tJac)$,
$\Gf_i'=\Gf_i+r_i$, and for the new sections~$s_i'$, $t_{ij}'$ one
has $t_{ij}'=t_{ij}+r_j-r_i$ and $s_i'=s_i-(1-c)r_i$. Thus, the
new cocycle
$(s'_i,t'_{ij})$ differs from $(s_i,t_{ij})$ by the coboundary of
the $0$-cochain $(r_i)$ (and, \emph{vice versa}, changing
the cocycle $(s_i,t_{ij})$ by the coboundary of $(r_i)$ can be realized by
replacing the isomorphisms~$\Gf_i$ with $\Gf_i+r_i$).

Conversely, let $(s_i, t_{ij})$ be a $1$-cocycle. Since $(t_{ij})$
is a $1$-cocycle in the ordinary \v{C}ech complex, it defines a
complex elliptic surface~$E$ (by gluing the pieces~$J_i$
along their intersections \emph{via} the translations~$t_{ij}$).
Since $(1+c)s_i=0$, the anti-automorphisms $c_J+s_i$ are real
structures on~$J_i$, and the cocycle condition guarantees
that these real structures agree on the intersections, thus
blending into a real structure on~$E$.
\endproof

In addition to the sheaves $\tJac$, $\Jac$, $\RR1p\CO_J$ and exact
sequences~\eqref{stJac}, \eqref{sJac}
consider the sheaves $\tJac\top$, $\Jac\top$, $(\RR1p\CO_J)\top$
of continuous sections of the corresponding bundles/fibrations
and exact sequences
$$
\gather
0@>>>\Jac\top@>>>\tJac\top@>>>\CS@>>>0,\eqtag\label{sttop}\\
0@>>>\RR1p\Z^-_J@>>>(\RR1p\CO_J)\top@>>>\Jac\top@>>>0.\eqtag\label{stop}
\endgather
$$
(For the discrete sheaves $\RR1p\Z^-_J$ and~$\CS$ one would have
$(\RR1p\Z^-_J)\top=\RR1p\Z^-_J$ and $\CS\top=\CS$.)

\lemma\label{vector.space}
If $\CA$ is a sheaf of $\R$-vector spaces and the involution
$c\:\CA\to c_B^*\CA$ is $\R$-linear, then
$\HH^*(\CA)=(H^*(B;\CA))^c$ \rom(the subspace of $c$-invariant
classes\rom).
In particular, $\HH^i(\RR1p\CO_J)=0$ for $i>1$, and
$\HH^i((\RR1p\CO_J)\top)=0$ for $i>0$.
\endlemma

\proof
For any involution~$c$ on a vector space~$U$ over a field of
characteristic~$0$ the sequence $U@>1-c>>U@>1+c>>U$ is exact, and
the first statement follows from~\eqref{GaloisII}. The rest is
immediate.
\endproof

\corollary\label{Jac} There is a natural exact sequence
\rom(induced from~\eqref{sJac}\rom)
$$
\HH^1(\RR1p\Z^-_J)@>>>(H^1(B;\RR1p\CO_J))^c@>>>
 \HH^1(\Jac)@>>>\HH^2(\RR1p\Z^-_J)@>>>0.
$$
In particular, the discrete part of $\HH^1(\Jac)$ is canonically
isomorphic to $\HH^2(\RR1p\Z^-_J)$.
Furthermore, there are natural isomorphisms
$\HH^i(\Jac)=\HH^{i+1}(\RR1p\Z^-_J)$, $i>1$, and
$\HH^i(\Jac\top)=\HH^{i+1}(\RR1p\Z^-_J)$, $i>0$.
\qed
\endcorollary

\theorem\label{Rshatop}
There is a canonical isomorphism $\Rsha\top(J)=\HH^1(\tJac\top)$
and a natural \rom(with respect to
real fiberwise
homeomorphisms\rom) exact sequence
$$
\HH^0(\CS)@>>>\HH^2(\RR1p\Z^-_J)@>>>\Rsha\top(J)@>>>
 \HH^1(\CS)@>>>\HH^3(\RR1p\Z^-_J).
$$
In particular, $\Rsha\top(J)$ is a topological invariant
of the pair $(p,c_J)$.
\endtheorem

\proof
In view of Theorem~\ref{Rsha}, it suffices to show that the group
$\HH^1(\tJac\top)$ is discrete, the inclusion homomorphism
$\HH^1(\tJac)\to\HH^1(\tJac\top)$ is onto, and its kernel is
connected. The last two statements follow from Corollary~\ref{Jac}
and (an obvious extension of) the $5$-lemma applied to the
commutative diagram
$$
\minCDarrowwidth{.25in}
\CD
\HH^0(\CS)@>>>\HH^1(\Jac)@>>>\HH^1(\tJac)@>>>\HH^1(\CS)@>>>\HH^2(\Jac)\\
@|@V{i}V\roman{onto}V@V{\tilde\imath}VV@|@V\cong VV\\
\HH^0(\CS)@>>>\HH^1(\Jac\top)@>>>\HH^1(\tJac\top)@>>>\HH^1(\CS)@>>>\HH^2(\Jac\top)
 \rlap{.}
\endCD
$$
One obtains that $\tilde\imath$ is an epimorphism and the induced
homomorphism $\Ker i\to\Ker\tilde\imath$ is onto. The exact sequence in the
statement of the theorem (and the discreteness of
$\HH^1(\tJac\top)$) follow then from the bottom row of the diagram
\via\ the identification $\HH^i(\Jac\top)=\HH^{i+1}(\RR1p\Z^-_J)$
given by Corollary~\ref{Jac}.
\endproof

We conclude this subsection with an explanation of the count of real
homological invariants corresponding to a fixed real functional
invariant, see~\ref{realhominv}.

\theorem\label{HHgen}
If not empty, the set of real local systems
belonging to and concordant with a given real map $j\:B\to\Cp1$ is
an affine space over $\HH^1(B^\sharp;\Z_2)$, where
$B^\sharp=j^{-1}(\Cp1\sminus\{0,1,\infty\})$.
\endtheorem

\proof
The proof repeats literally that of Theorem~\ref{Rsha}, with~$E$ and~$J$
replaced with, respectively, any and one of the local systems in
question, and $\tJac$ replaced with the constant
sheaf~$(\Z_2)_{B^\sharp}$.
\endproof

Thus,
the number of real local systems
as in
Theorem~\ref{HHgen} is
either~$0$
or $2^r$, where $r=\dim\HH^1(B^\sharp;\Z_2)$. The latter dimension
is given by the following proposition.

\proposition\label{HHcount}
For a
connected
complex curve~$B$
and
for any
real structure $c\:B\to B$
one has
$$
\dim\HH^1(B;\Z_2)=\dim(H^1(B;\Z_2))^c+1-\Ge,
$$
where $\Ge=1$ if $B_{\R}=\varnothing$ and either $B$ is not
compact or the genus of~$B$ is odd, and $\Ge=0$ otherwise.
\endproposition

\proof
For~$\HH^1$, the spectral sequence~\eqref{GaloisII} reduces to the
exact sequence
$$
0@>>>H^1(\Z_2;H^0(B;\Z_2))@>>>\HH^1(B;\Z_2)@>>>
 H^0(\Z_2;H^1(B;\Z_2))@>\Rd2>>\ldots\,.
$$
Here, $H^1(\Z_2;H^0(B;\Z_2))=\Z_2$ and
$H^0(\Z_2;H^1(B;\Z_2))=(H^1(B;\Z_2))^c$.
Stabilizing
(\cf.~\ref{K.Z2}),
one observes that the differential
$$
\Rd2\:H^0(\Z_2;H^1(B;\Z_2))\to H^2(\Z_2;H^0(B;\Z_2))=\Z_2
$$
is nontrivial if and only if so is the differential
$\Rd2\:\RH2^2\to\RH2^0$ in Kalinin's spectral sequence of~$B$. If
$B_{\R}\ne\varnothing$, then $\Rd2=0$ as $\RH2^0$ must survive
to~$\iH^0$. If $B_{\R}=\varnothing$ and $B$ is not compact, \ie,
$H^2(B;\Z_2)=0$, then $\Rd2\ne0$ as this is the only chance to
kill~$\RH2^0=\Z_2$. Finally, if $B_{\R}=\varnothing$ and $B$ is
compact, the involution is standard and a direct calculation shows
that $\Rd2=0$ if and only if the genus of~$B$ is even.
\endproof

\subsection{The case of generic singular fibers}
Fix a real Jacobian elliptic fibration $p\:J\to B$ with a real
section $s\:B\to J$. Recall that, for any abelian group~$G$,
$\RR0pG_J$ is the constant sheaf~$G_B$ and $\RR2pG_J$ is an
extension of~$G_B$ by a skyscraper sheaf~$\CS'$ concentrated at
the points of~$B$ corresponding to reducible fibers of~$p$.

\lemma\label{direct.sum.1}
For any abelian group~$G$, the homomorphism
$p^*\:H^*(B;G)\to H^*(J;G)$ \rom(equivalently, the edge
homomorphism $H^*(B;\RR0pG_J)\to H^*(J;G)$ of the Leray spectral
sequence of~$p$\rom) is a monomorphism and its image is a direct
summand, the decomposition respecting~$c_J^*$. Furthermore,
$p^*$ embeds Kalinin's spectral sequence of $(B,c_B)$ \rom(both
over~$\Z_2$ and over~$\Z$\rom) as a direct summand into Kalinin's
spectral sequence of $(J,c_J)$.
\endlemma

\proof
All statements follow immediately from the existence of a section;
the complementary direct summand is $\Ker s^*$, where $s^*$ is the
appropriate induced homomorphism.
\endproof

\corollary\label{d2}
Starting from $\RH2$,
the only
potentially nontrivial
differential in the Leray spectral sequence of~$p$ is
$\Rd2\:H^0(B;\RR2pG_J)\to H^2(B;\RR1pG_J)$.
\qed
\endcorollary

\lemma\label{J.nontrivial}
If $J$ is compact and non-isotrivial, then $H^0(B;\RR1p\Z_J)=0$ and for each
$r\ge1$, $m\in\Z_2$,
and $q=0,1$ the induced homomorphism
$p^*\:\rH^{mq}(B;\Z)\to\rH^{mq}(J;\Z)$ is an isomorphism.
As a consequence, $p^*\:\CF^m_q(B;\Z)\to\CF^m_q(J;\Z)$ is an
isomorphism for all $m\in\Z_2$ and $q=0,1$.
\endlemma

\proof The isomorphism $p^*\:\RH1^{m,0}(B;\Z)=\RH1^{m,0}(J;\Z)$ is
obvious. Since $J$ is non-isotrivial, its Mordel-Weil group is
discrete,
and from the
exact sequence~\eqref{sJac} one concludes that
$H^0(B;\RR1p\CO_J)=0$ and, hence, $H^0(B;\RR1p\Z_J)=0$. Thus,
$p^*\:\RH1^{m,1}(B;\Z)\to\RH1^{m,1}(J;\Z)$ is also an isomorphism.
Since $p^*$ is a direct summand embedding, it remains an
isomorphism for all $r\ge1$.
\endproof

\lemma\label{J.compact}
If $J$ is compact and non-isotrivial, then
\roster
\item\local1
the groups $H_*(J;\Z)$ and $H^*(J;\Z)$ are torsion free\rom;
\item\local3
the `edge'
homomorphism $H^3(J;\Z)\to H^1(B;\RR2p\Z_J)$ is an
isomorphism\rom;
\item\local4
the differential~$\Rd2\:H^0(B;\RR2p\Z_J)\to H^2(B;\RR1p\Z_J)$
establishes an isomorphism
$\Coker[H^2(J;\Z)\to H^0(B;\RR2p\Z_J)]=H^2(B;\RR1p\Z_J)$.
\endroster
\endlemma

\proof
As is known (see, for example,
\cite{FM, Proposition 2.2.1}), the induced homomorphism
$p_*\:\pi_1(J)\to\pi_1(B)$ is an isomorphism. Hence,
$p_*\:H_1(J;\Z)\to H_1(B;\Z)$ is also an isomorphism, the group
$H_1(J;\Z)$ is torsion free, and so are $H_*(J;\Z)$ and
$H^*(J;\Z)$.
Statement~\loccit3 follows then from the isomorphism
$p^*\:H^1(B;\Z)\to H^1(J;\Z)$, see Lemma~\ref{J.nontrivial},
and Poincar\'{e} duality, and \loccit4~is immediate.
\endproof

\remark{Remark}
An alternative description of the group $H^2(B;\RR1p\Z_J)$, due to
Kodaira, is as the group of coinvariants of~$p$, \ie,
$H^2(B;\RR1p\Z_J)=H^1(F_b;\Z)/(1-m_\gamma)a$,
where $F_b$ is the fiber over a point $b\in B^\sharp$,
$a$ runs through $H^1(F_b;\Z)$, $\gamma$ runs through
$\pi_1(B^\sharp,b)$, and
$m_\gamma\:H^2(F_b;\Z)\to H^2(F_b;\Z)$ is the monodromy
along~$\gamma$.
\endremark

\corollary\label{Gmax}
Assume that $J$ is compact and non-isotrivial and that the real
part~$J_{\R}$ is nonempty \rom(equivalently, the real
part~$B_{\R}$ is nonempty\rom). Then $J$ is both $\Z_2$- and
$\Z$-Galois maximal.
\endcorollary

\proof
Since both $H^*(B;\Z)$ and $H^*(J;\Z)$ are torsion free, reduction
modulo~$2$ induces isomorphisms of their Kalinin's spectral
sequences with coefficients in~$\Z$ and~$\Z_2$, see~\ref{K.ZtoZ2}.
In particular, $\Z$- and $\Z_2$-Galois maximality are equivalent.
Furthermore, Lemma~\ref{J.nontrivial} implies that
$p^*\:\rH^q(B;\Z_2)\to\rH^q(J;\Z_2)$ is an isomorphism for all
$r\ge1$ and $q=0,1$. The curve~$B$ with nonempty real part is
$\Z_2$-Galois maximal. Hence, the differentials~$\rd$, $r\ge2$,
landing in $\rH^q(J;\Z_2)$, $q=0$,~$1$, are trivial, and using
Poincar\'{e} duality~\ref{K.PoincareZ2} one concludes that so are all
differentials~$\rd$, $r\ge2$, \ie, $J$ is $\Z_2$-Galois maximal.
\endproof

\lemma\label{J.simple}
Assume that $J$ is compact and non-isotrivial and that all its
fibers are irreducible. Then
\roster
\item\local1
one has $H^2(B;\RR1p\Z_J)=0$\rom;
\item\local2
there are isomorphisms
$\HH^i(\RR1p\Z^-_J)=H^{i-1}(\Z_2;H^1(B;\RR1p\Z^-_J))$, $i\ge0$,
\endroster
\rm(as usual $H^{-1}=0$\rm).
If, in addition, $B_{\R}\ne\varnothing$, then
\roster
\item[3]\local3
there is an isomorphism
$\HH^2(\RR1p\Z^-_J)=H^2(\Z_2;H^2(J;\Z))$.
\endroster
\endlemma

\proof
Pick a generic fiber~$F$ of~$p$ and let $\inc\:F\hookrightarrow J$
be the inclusion. Denote by $f,b\in H_2(J;\Z)$ the
fundamental classes of~$F$ and $s(B)$, respectively. One has
$f^2=0$ and $f\circ b=1$. Hence, $f$ and~$b$ span a unimodular
sublattice in $H_2(J;\Z)$ and there is an orthogonal (with respect
to the intersection index form) direct sum decomposition
$H_2(J;\Z)=\<f,b\>\oplus\<f,b\>\ort$. The Poincar\'{e} duality
yields then an orthogonal decomposition
$H^2(J;\Z)=\Hom(\<f,b\>,\Z)\oplus\Ker(\mathord{\inc^*}\oplus s^*)$.
In particular, the induced homomorphism
$\inc^*\:H^2(J;\Z)\to H^2(F;\Z)$ is onto.

Since all fibers are irreducible,
one has $\RR2p\Z_J=\Z_B$ and $H^0(B;\RR2p\Z_J)=\Z=H^2(F;\Z)$, and
the edge homomorphism $H^2(J;\Z)\to H^0(B;\RR2p\Z_J)$ coincides
with~$\inc^*$. Since the latter is onto, \iref{J.compact}4
implies~\loccit1. Then $H^1(B;\RR1p\Z_J)$ is the only nontrivial
cohomology of $\RR1p\Z_J$, and~\loccit2 follows from the spectral
sequence~\eqref{GaloisII}. In particular, using the definition
of $H^*(\Z_2;\,\cdot\,)$ given in \ref{Z2cohomology}, one has
$\HH^2(\RR1p\Z^-_J)=H^1(\Z_2;H^1(B;\RR1p\Z^-_J))=H^2(\Z_2;H^1(B;\RR1p\Z_J))$.

Now, assume that $B_{\R}\ne\varnothing$. Then $F$ can be chosen
real, and the orthogonal decomposition above is
$c_J^*$-equivariant. In view of Lemma~\ref{direct.sum.1}, one can
identify $H^1(B;\RR1p\Z_J)$ with
$\Ker(\mathord{\inc^*}\oplus s^*)\subset H^2(J;\Z)$. Hence, one has
$$
H^2(\Z_2;H^2(J;\Z))=H^2(\Z_2;\Hom(\<f,b\>,\Z))\oplus
 H^2(\Z_2;H^1(B;\RR1p\Z_J)),
$$
and it remains to notice that, $F$ and~$s(B)$ being analytic
curves, $c_*$ acts via minus identity on $\<f,b\>$ and, hence,
$H^2(\Z_2;\Hom(\<f,b\>,\Z))=0$.
\endproof

We are ready to prove the principal result of this section.

\theorem\label{RSha.BR}
Let $p\:J\to B$ be a compact non-isotrivial real Jacobian elliptic
surface with irreducible fibers and nonempty real part. Then there
are canonical isomorphisms
$\Rsha\top(J)=H^2(\Z_2;H^2(J;\Z))=\CF^0_2(J;\Z)/\CF^0_1(J;\Z)$,
where $\{\CF^0_q\}$ is Kalinin's filtration, see~\ref{K.Z}.

Alternatively, there is a canonical isomorphism
$\Rsha\top(J)=(\Im p^*)\ort/\Im p^*$,
where $p^*$ is the induced
homomorphism $H\even(B_{\R};\Z_2)\to H\even(J_{\R};\Z_2)$
and the
complement is with respect to a certain perfect pairing
$H\even\otimes H\even\to\Z_2$.
\endtheorem

\proof
Under the assumptions one has $\CS=0$. Hence, Theorem~\ref{Rshatop}
and Lemma \iref{J.simple}3 imply
$\Rsha\top(J)=\HH^2(\RR1p\Z^-_J)=H^2(\Z_2;H^2(J;\Z))=\RH2^{0,2}(J;\Z)$.
In view of Corollary~\ref{Gmax}, Kalinin's spectral sequence
degenerates at~$\RH2$ and the latter group equals
$\iH^{0,2}(J;\Z)=\CF^0_2/\CF^0_1$.

Since the group $H^*(J;\Z)$ is torsion free,
see~\iref{J.compact}1, the map $1+\Sq^1$ establishes
isomorphisms
$\CF^0_q(J;\Z)=(1+\Sq^1)\CF_q(J;\Z_2)\cap H\even(J_\R;\Z_2)$,
see~\ref{K.ZtoZ2}. On the other hand, with respect to the pairing
described in~\ref{K.PoincareZ2},
one has $\CF_2=\CF_1\ort$,
and, in view of Lemma~\ref{J.nontrivial}, the map
$p^*\:\CF_1(B;\Z_2)\to\CF_1(J;\Z_2)$
is an isomorphism. Finally,
since~$B$ is a compact curve,
the intersection
$(1+\Sq^1)\CF_1(B;\Z_2)\cap H\even(B_\R;\Z_2)$ is merely the group
$H^0(B_\R;\Z_2)=H\even(B_\R;\Z_2)$.
\endproof

\subsection{The geometric interpretation}\label{S.geometry}
Numerically, Theorem~\ref{RSha.BR} states that, under the
assumptions, the discrete part $\Rsha\top(J)$ is a $\Z_2$-vector
space of dimension twice the number of extra components
of~$J_\R$.
Indeed, each
component of~$J_\R$ contributes~$2$ to $\dim H\even(J_{\R};\Z_2)$,
and the passage to $(\Im p^*)\ort/\Im p^*$ kills the contribution
of the
principal components. Below we
show that, in fact, each extra component of~$J_\R$ does contribute
a pair of $\Z_2$ summands in a natural way.

As in Theorem~\ref{RSha.BR}, fix a compact non-isotrivial real
Jacobian elliptic surface $p\:J\to B$
with irreducible fibers and nonempty real part. Let
$J'=p^{-1}(B_\R)$, and regard $p\:J'\to B_\R$ as a real `Jacobian
fibration.'
In particular, one can consider the discrete
Tate-Shafarevich group $\Rsha\top(J')$. Clearly,
Theorem~\ref{Rshatop} still applies and yields a natural
isomorphism $\Rsha\top(J')=\HH^2(\RR1p\Z^-_{J'})$. Note that the
sheaf $\RR1p\Z^-_{J'}$ is the restriction $i^*\RR1p\Z^-_J$, where
$i\:B_\R\to B$ is the inclusion.

\lemma\label{onto.JR}
Let $J$ and~$J'$ be as above.
Then the inclusion homomorphism
$i^*\:\HH^2(\RR1p\Z^-_{J})\to\HH^2(\RR1p\Z^-_{J'})$ is onto.
\endlemma

\proof
Let $\Cal K$ be the kernel of the epimorphism
$\RR1p\Z^-_{J}\to i_*\RR1p\Z^-_{J'}$. It suffices to show that
$\HH^3(\Cal K)=0$. Since $\Cal K$ is trivial
over~$B_\R$, one has $\CH^p(\Z_2;\Cal K)=0$ for all $p>0$, and the
spectral sequence~\eqref{GaloisI} reduces to the isomorphisms
$\HH^q(\Cal K)=H^q(B/c;\CH^0(\Z_2;\Cal K))$.
As $3>\dim B/c=2$, the statement follows.
\endproof

Similarly, one can speak about the discrete Tate-Shafarevich group
$\Rsha\top(F)$ of a single fiber $F=p^{-1}(b)$, $b\in B_\R$.
As above, for an irreducible fiber~$F$ one has
$\Rsha\top(F)=\HH^2(\RR1p\Z^-_F)=H^1(\Z_2;H^1(F;\Z))$. (As usual, the
dimension shift
is due to the non-standard action of the
involution on~$\Z^-_F$.)
This group is either~$\Z_2$ or~$0$, see
Table~\ref{table.H1}. In the case $\Rsha\top(F)=\Z_2$, the
non-trivial element corresponds to the real structure on~$F$ with
respect to which the normalization of~$F$ has empty real part.

\topinsert
\table\label{table.H1}
The groups $H^i(\Z_2;H^1(F;\Z))$
\endtable
\centerline{\vbox{%
\offinterlineskip
\halign{\vrule depth3pt height11pt\ #\hfil\ &\vrule\ \hfil$#$\hfil\ &
 \vrule\ \hfil$#$\hfil\ \vrule\cr
\noalign{\hrule}
\hfil Fiber~$F$&H^1(\Z_2;H^1(F;\Z))&H^2(\Z_2;H^1(F;\Z))\cr
\noalign{\hrule}
a nonsingular $M$-curve&\Z_2&\Z_2\cr
a nonsingular $(M-1)$-curve&0&0\cr
$\Cp1/\{0\sim\infty\}$ with $\conj\:z\mapsto\bar z$&\Z_2&0\cr
$\Cp1/\{0\sim\infty\}$ with $\conj\:z\mapsto1/\bar z$&0&\Z_2\cr
rational curve with a cusp&0&0\cr
\noalign{\hrule}}}}
\endinsert

Denote by $J_\R\extra$ the union of the extra components
of~$J_\R$, and let $B_\R\extra=p_*(J_\R\extra)$.

\theorem\label{Rsha(JR)}
Let $p\:J\to B$ be a compact non-isotrivial real Jacobian elliptic
surface with irreducible fibers and nonempty real part, and let
$J'=p^{-1}(B_\R)$. Then there is a natural exact sequence
$$
0@>>>
 H^1(B_\R;\CH^1(\Z_2;\RR1p\Z^-_{J'}))
 @>\alpha>>
 \Rsha\top(J')
 @>\beta>>
 H^0(B_\R;\CH^2(\Z_2;\RR1p\Z^-_{J'}))
 @>>>0
$$
and natural isomorphisms
$$
\align
H^1(B_\R;\CH^1(\Z_2;\RR1p\Z^-_{J'}))&=
 H^1(B_\R\extra,\partial B_\R\extra;\Z_2),\\
H^0(B_\R;\CH^2(\Z_2;\RR1p\Z^-_{J'}))&=
 H^0(B_\R\extra;\Z_2)=
 \bigoplus\Rsha\top(F_i),
\endalign
$$
where $F_i=p^{-1}(b_i)$ are the fibers over some points $b_i\in
B_\R$, one in the interior of each connected component
of~$B_\R\extra$. The composition of the last isomorphism
and~$\beta$ coincides with the homomorphism induced by the
inclusion $\bigcup\{b_i\}\hookrightarrow B_\R$.
\endtheorem

\proof
Since $\dim B_\R=1$, one has $H^q(B_\R;\,\cdot\,)=0$ for $q>1$
and, hence,
the spectral sequence~\eqref{GaloisI} for
$\Rsha\top(J')=\HH^2(\RR1p\Z^-_{J'})$ collapses at~$E_2$ and
results in the exact sequence in the statement.

The stalks of the sheaves $\CH^p(\Z_2;\RR1p\Z^-_{J'})$, $p=1,2$,
are given by Table~\ref{table.H1}. The stalks are at most~$\Z_2$,
and the sheaves are supported by the closure $\bar B_\R^M$ of the
subset $B_\R^M$ of the points of~$B_\R$ whose pull-back is a
nonsingular $M$-curve. Hence, the cohomology groups are
$H^q(B_\R;\,\cdot\,)=H^q(\bar B_\R^M,D;\Z_2)$, where $D$ is the
part of the boundary $\partial\bar B_\R^M$ over which
the sheaf in question
has trivial stalks. Considering the possibilities for a component
of~$\bar B^M_\R$, one easily concludes that only the components of
$B_\R\extra\subset\bar B^M_\R$ make nontrivial contributions to the
cohomology; this observation gives the isomorphisms in the statement.
\endproof

\corollary\label{iso.JR}
Under the hypothesis of Theorem \ref{Rsha(JR)},
the inclusion homomorphism
$i^*\:\Rsha\top(J)\to\Rsha\top(J')$
is an isomorphism.
\endcorollary

\proof
According to Theorem~\ref{Rsha(JR)}, each extra component
of~$J_\R$ contributes~$4$ to the order of $\Rsha\top(J')$. Hence,
the two groups are of the same order, and the statement follows from
Lemma~\ref{onto.JR}.
\endproof

\paragraph\label{modifications}
Combining Theorem~\ref{Rsha(JR)} and the proof of
Theorem~\ref{Rsha} (namely, the part explaining how a real elliptic
fibration can be modified), one can easily describe the real parts
of all, not necessarily Jacobian, elliptic surfaces. Each extra
component $X\subset J_\R$ contributes two $\Z_2$ summands to
$\Rsha(J)$: one to the subgroup
$H^1=H^1(B_\R;\CH^1(\Z_2;\RR1p\Z^-_{J'}))$, and one to the
quotient $H^0=H^0(B_\R;\CH^2(\Z_2;\RR1p\Z^-_{J'}))$. The
non-trivial element of~$H^1$ represents a real modification of
the fibration, given by a real $1$-cocycle, \cf. the proof of
Theorem~\ref{Rsha}.
If $X$ is a sphere,
Figure~\ref{fig.MM}(a), the resulting modification of the real
part is shown in Figure~\ref{fig.MM}(b), the
cocycle being the
partial section shown by a gray dotted line in Figure~\ref{fig.MM}(a).
If $X$ is a torus or a Klein bottle,
the two components over the corresponding hyperbolic component of~$B_\R$
are intertwined into one. Note that, since the restriction
of~$H^1$ to each fiber is trivial, the real structures of all
fibers remain intact.

\midinsert
\eps{S4Fig1}
\figure\label{fig.MM}
Modifications of the real part
\endfigure
\endinsert

The non-trivial element of~$H^0$ contributed by~$X$ restricts
non-trivially to each fiber~$F$ over the projection
$p(X)\subset B_\R$, thus resulting in
rebuilding the real structure of all fibers. If $X$ is a torus or
a Klein bottle,
the result has empty real part (over the corresponding hyperbolic
component of~$B_\R$). If $X$ is a sphere, the result is shown in
Figure~\ref{fig.MM}(c) and~(d). The new real structure is obtained
via the shift by a real section of the
opposite
Jacobian fibration. One such section,
which is
homotopically
non-trivial
in each fiber over $p(X)$,
is shown in gray in Figure~\ref{fig.MM}(a). Note that the real
parts shown in Figure~\ref{fig.MM}(c) and~(d) are homeomorphic;
the two figures are intended to indicate the fact that one can
just change the real structure or change both the real structure
and the fibration.

\corollary\label{section->Jacobian}
A compact non-isotrivial real elliptic surface $p\:E\to B$ with
irreducible fibers is equivariantly deformation equivalent
to a real Jacobian
surface
if and only if the restriction
$p\:E_\R\to B_\R$ admits a continuous section.
\endcorollary

\proof
The `only if' part is obvious. For the `if' part, it suffices to
notice that none of the fibrations shown in
Figure~\ref{fig.MM}(b)--(d) (neither the results of modifying a
toroidal extra component) has a section.
\endproof

\proposition\label{Betti-4}
Let $E$ be a compact non-isotrivial real elliptic
surface that is not equivariantly deformation equivalent
to its Jacobian $J=J(E)$. Then
$$
\dim H^*(E_\R;\Z_2)\le\dim H^*(J_\R;\Z_2)-4.
$$
\endproposition

\proof
Each modification shown on Figure~\ref{fig.MM}, from~(a) to any
of~(b), (c), (d), either leaves the total Betti number
$\dim H^*(\,\cdot\,;\Z_2)$
unchanged or reduces it by~$4$, depending on whether a new
component is created or not. As the first modification of this
kind, starting from $J=J(E)$, does not create a new
component, it does reduce the total Betti number.
Similarly, each nontrivial
modification over a hyperbolic component of~$B_\R$
turns a pair of tori (or Klein bottles) into a single
torus (respectively, Klein bottle) or the empty set, thus reducing
the total Betti number by~$4$ or~$8$,
respectively.
\endproof

\corollary\label{M->Jacobian}
A compact non-isotrivial real elliptic surface that is an $M$- or
$(M-1)$-variety is equivariantly deformation equivalent to its Jacobian.
\qed
\endcorollary

\remark{Remark}
The same arguments as above show that under the same hypotheses
the inequality
$$
\dim H^k(E_\R;\Z_2)\le\dim H^k(J_\R;\Z_2),
$$
where $J=J(E)$ is the Jacobian, holds for any $k$.
\endremark

\subsection{Deformations}
Let $p\:X\to S$, $\pi\:S\to D$ be an elementary real deformation
of non-isotrivial compact real Jacobian elliptic surfaces.
The projection $p\:X\to S$ can be regarded as an
elliptic fibration; hence, one can define the corresponding real
Tate-Shafarevich group $\Rsha(X)$, the sheaves $\tJac$ and
$\Jac=\RR1p\CO_X/\RR1p\Z^-_X$, and their topological counterparts
$\Rsha(X)\top$, $\tJac\top$, and~$\Jac\top$. As in
Theorems~\ref{Rsha} and~\ref{Rshatop}, one has
$\Rsha(X)=\HH^1(\tJac)$ and $\Rsha\top(X)=\HH^1(\tJac\top)$.

\theorem\label{Rsha.onto}
Assume
that the members $p_t\:X_t\to S_t$, $t\in D$, of the family, except
possibly~$X_0$, have no singular fibers other than those of
type~$\I_1$. Then  the restriction homomorphism
$\Rsha(X)\to\Rsha(X_0)$ is onto, and the restriction homomorphism
$\Rsha\top(X)\to\Rsha\top(X_0)$ is an isomorphism.
\endtheorem

\corollary\label{deformation}
Let $E$ be a compact real elliptic surface, $J$ its Jacobian, and
$p\:X\to D$, $\pi\:S\to D$ a real deformation of~$J$ \rom(\ie, $X_0=J$\rom) as in
Theorem~\ref{Rsha.onto}. Then there is a real deformation
$\tilde p\:\tilde X\to D$ of $E$ \rom(with the same base $\pi\:S\to D$\rom)
such that $X_t=J(\tilde X_t)$
for each $t\in D$.
\qed
\endcorollary

\proof[Proof of Theorem~\ref{Rsha.onto}]
Consider the sheaves $\Jac$, $\tJac$, and~$\CS$ on~$S$
and denote by $\Jac_0$, $\tJac_0$,
and~$\CS_0$, respectively, their restrictions to~$S_0$.
There are exact sequences
$$
0@>>>\Jac@>>>\tilde\Jac@>>>\CS@>>>0,\qquad
0@>>>\Jac_0@>>>\tilde\Jac_0@>>>\CS_0@>>>0,
$$
and, since the restriction homomorphisms
$\HH^*(\CS)\to\HH^*(\CS_0)$ are isomorphisms (as $\CS$ is
concentrated at points of~$S_0$), it suffices to show
that the homomorphism $\HH^i(\Jac)\to\HH^i(\Jac_0)$ is onto
for $i=1$ and one-to-one for $i=2$.

To this end,
we
compare the exact sequences
$$
\CD
0@>>>\RR1p\Z^-_X@>>>\RR1p\CO_X@>>>\Jac@>>>0\rlap,\\
0@>>>\RR1p\Z^-_{X_0}@>>>\RR1p\CO_{X_0}@>>>\Jac_0@>>>0.
\endCD
$$
For the coherent sheaves $\RR1p\CO_X$ and $\RR1p\CO_{X_0}$ one has
$\HH^i(\,\cdot\,)=H^0(\Z_2;H^i(\,\cdot\,))$ and
$H^j(\Z_2;H^i(\,\cdot\,))=0$ for $j>0$ (\cf\.
Lemma~\ref{vector.space}). Furthermore, for each fiber~$S_t$
of~$\pi$ one has $H^j(S_t;\RR1p\CO_{X_t})=0$ unless
$j=1$
(see \ref{sheafs}),
and the Leray
spectral sequence of~$\pi$ implies that
$H^i(X;\RR1p\CO_X)=H^{i-1}(D;\RR1\pi\RR1p\CO_X)$. Since $D$ is a
Stein manifold, one concludes that
$\HH^2(\RR1p\CO_X)=\HH^2(\RR1p\CO_{X_0})=0$ and the
homomorphism $\HH^1(\RR1p\CO_X)\to\HH^1(\RR1p\CO_{X_0})$ is onto.
Thus, it remains to prove that the restriction homomorphism
$\HH^1(\RR1p\Z^-_X)\to\HH^1(\RR1p\Z^-_{X_0})$ is also onto. We will show
that it is, in fact, an isomorphism.

Informally, the last assertion follows from the isomorphism
$H^*(X)=H^*(X_0)$. More precisely,
observe that $\RR0p\Z_X$ and $\RR0p\Z_{X_0}$ are both constant sheaves
(with fiber~$\Z$) and $\RR2p\Z_X$ and $\RR2p\Z_E$ are extensions
of constant sheaves by a sheaf concentrated at points of~$X_0$.
Hence, both $H^*(\RR0p\Z_X)\to H^*(\RR0p\Z_{X_0})$ and
$H^*(\RR2p\Z_X)\to H^*(\RR2p\Z_{X_0})$ are isomorphisms. From
comparing the Leray spectral sequences for the projections
$X\to S$ and $X_0\to S_0$ it follows that the maps
$H^i(\RR1p\Z_X)\to H^i(\RR1p\Z_E)$ are also isomorphisms (for
all~$i$), and the statement follows from the
spectral sequence~\eqref{GaloisII}.

The proof of the topological statement is similar, except that
for fine sheaves one
has $H^j(S_t;(\RR1p\CO_{X_t})\top)=0$ unless $j=0$
(for each $t\in D$) and, hence,
$\HH^i(\RR1p\CO_X)=\HH^i(\RR1p\CO_{X_0})=0$ for $i=1$,~$2$.
Alternatively, the assertion can be derived from
Theorem~\ref{Rshatop} and a similar exact sequence for
$\Rsha\top(X)$.
\endproof

\end


\ifx\defloaded\undefined\input elliptic.def \fi

\begin

\setcounter\section4


\section{Real trigonal curves and dessins d'enfants}\label{rtc}

We start this section
by
introducing
the
notion of trichotomic
graph. It is
a real version of Grothendieck's \emph{dessins d'enfants},
which is adjusted
for
dealing
with
real meromorphic
functions defined on
a real curve and
having a certain
preset ramification over the three real points
$0,1,\infty\in\Cp1$.
More precisely, a trichotomic graph is the quotient by the complex
conjugation of
a properly decorated pull-back of $(\Rp1;0,1,\infty)$; the
pull-backs of~$0$, $1$, and~$\infty$ being marked with \black--,
\white--, and \cross-- respectively.
Note that the function may (and usually does) have other ramification points,
which are ignored unless they are real.

\subsection{Trichotomic graphs}\label{tg}
Let $\base$ be a (topological) compact connected surface, possibly
with boundary. (Unless specified otherwise, in the topological
part of this section we are working in the {\sl PL\/}-category.)
We use the term \emph{real} for points, segments, \etc.
situated at the boundary~$\partial\base$. For a graph
$\Gamma\subset\base$, we denote by~$\Base\Gamma$ the closed cut
of~$D$ along~$\Gamma$. The connected components of~$\Base\Gamma$
are called \emph{regions} of~$\Gamma$.

A \emph{trichotomic graph} on~$\base$
is an embedded oriented
graph $\Gamma\subset\base$
decorated with the
following additional structures (referred to as \emph{colorings}
of the edges and vertices of~$\Gamma$, respectively):
\roster
\item"--"
each edge of~$\Gamma$ is of one of the three kinds: \solid-, \bold-,
or \dotted-;
\item"--"
each vertex of~$\Gamma$ is of one of the four kinds: \black-, \white-,
\cross-, or monochrome (the vertices of the first three kinds being called
\emph{essential});
\endroster
and satisfying the following conditions:
\roster
\item\local{tg-boundary}
the boundary $\partial\base$ is a union of
edges
and vertices
of~$\Gamma$;
\item\local{tg-valency}
the valency of each essential vertex of~$\Gamma$ is at
least~$2$,
and the valency of each monochrome vertex of~$\Gamma$ is at
least~$3$;
\item\local{tg-oriented}
the orientations of the edges of~$\Gamma$ form an orientation
of the boundary $\partial\Base\Gamma$;
this orientation extends to an orientation of~$\Base\Gamma$;
\item\local{tg-monochrome}
all edges incident to a monochrome vertex are of the same kind;
\item\local{tg-cross}
\cross-vertices are incident to incoming \dotted-edges and
outgoing \solid-edges;
\item\local{tg-black}
\black-vertices are incident to incoming \solid-edges and
outgoing \bold-edges;
\item\local{tg-white}
\white-vertices are incident to incoming \bold-edges and
outgoing \dotted-edges;
\item\local{tg-triangle}
each triangle (i.e., region with three essential vertices in
the boundary) is a topological disk.
\endroster
In \loccit{tg-cross}--\loccit{tg-white} the lists are complete,
\ie, vertices cannot be incident to edges of other kinds or with
different orientation.

In view of~\loccit{tg-monochrome}, the monochrome vertices can
further be subdivided into \solid-, \bold-, and \dotted-,
according to their
incident edges.
The sets of \solid-, \bold-, and \dotted- monochrome vertices
of~$\Gamma$ will be denoted by $\Gamma_{\msolid}$,
$\Gamma_{\mbold}$, and $\Gamma_{\mdotted}$, respectively.
The \emph{monochrome part}
of~$\Gamma$ of a given kind
(\solid-, \bold-, or \dotted-) is the union of (open)
edges and monochrome vertices of the corresponding kind.
Thus,
essential vertices {\bf never} belong to a monochrome part.

Condition~\loccit{tg-oriented} implies, in particular, that the
orientations of the edges incident to a vertex alternate. (This
statement is equivalent to the first part
of~\loccit{tg-oriented}.) Thus, all inner vertices of~$\Gamma$
have even valencies.

A number of examples of trichotomic graphs
is found further in this section (see, \eg, Figure~\ref{fig.cubics},
where complete graphs are drawn).

\paragraph\label{tg-chessboard}
Let $\Gamma$ be a trichotomic graph on~$\base$.
If $\base$ is orientable, a choice of the orientation defines a
chessboard coloring of~$\Base\Gamma$: a region
$\base_i\subset\Base\Gamma$ is said to be \emph{positive}
(\emph{negative}) if its orientation induced from~$\base$ coincides
with (respectively, is opposite to) that defined by~$\Gamma$.
Conversely, a chessboard coloring of~$\Base\Gamma$
defines
an
orientation of~$\base$.

\paragraph\label{tg-admissible}
A path in a trichotomic graph~$\Gamma$ is called \emph{monochrome}
if it belongs to a monochrome part of~$\Gamma$.
Given
two monochrome vertices $u,v\in\Gamma$, we say that $u\prec v$ if
there is an oriented monochrome path from~$u$ to~$v$. (Clearly,
only vertices of the same kind can be compatible.) The graph is
called \emph{admissible} if $\prec$ is a partial order. Since
$\prec$ is obviously transitive,
this condition is equivalent to the requirement that $\Gamma$
should have no oriented monochrome cycles.

\remark{Remark}
Note
that the orientation of~$\Gamma$ is almost
superfluous. Indeed, $\Gamma$ may have at most two orientations
satisfying~\loccit{tg-oriented}, and if $\Gamma$ has at least one
essential vertex, its orientation is uniquely determined by
\loccit{tg-cross}--\loccit{tg-white}.
Note also that (each connected component of) an admissible
graph does have essential vertices, as otherwise any component
of~$\partial\Base\Gamma$ would be an oriented monochrome cycle.
\endremark

\remark{Remark}
In fact, all three decorations of an admissible graph~$\Gamma$
(orientation and the two colorings) can be recovered from any of
the colorings. However, for clarity we retain both
colorings in the diagrams.
\endremark

%

\paragraph\label{tg-symmetric}
Let $B$ be the orientable double of~$\base$,
\ie, the orientation double covering of~$\base$ with the two preimages
of each real point $d\in\partial\base$ identified. (We
exclude the case when $\base$ is closed and orientable, as then
$B$ would be disconnected.) Denote by
$p\:B\to\base$ the projection and by $c\:B\to B$
its deck
translation, which is an orientation reversing involution.
We will show that trichotomic graphs on~$D$ are merely a way of
describing $c$-invariant trichotomic graphs on~$B$.
For this purpose,
given a trichotomic graph $\Gamma$ on $D$,
consider its
pull-back $\Gamma'=p^{-1}(\Gamma)$
and equip it
with the decorations
induced by~$p$.
Clearly, the deck translation
$c$ {\bf preserves} the decorations
of~$\Gamma'$, including its orientation.

\lemma\label{j-double}
Given a trichotomic graph $\Gamma\subset\base$,
its pull-back $\Gamma'=p^{-1}(\Gamma)$, with the decorations
induced by~$p$, is a $c$-invariant trichotomic graph on~$B$.
Conversely, given a $c$-invariant trichotomic graph
$\Gamma'\subset B$, its quotient $\Gamma=p(\Gamma')$ is a
trichotomic graph on~$\base$. The graph~$\Gamma'$ is admissible if
and only if so is~$\Gamma$.
\endlemma

\proof
The
direct
statement is immediate.
For the converse,
assume that $\Gamma'\subset B$ is a $c$-invariant trichotomic graph.
Since $c$ is orientation reversing, the graph $\Gamma'$ has the following
{\bf separation property}:
each region~$B_i$ of~$\Gamma'$ is disjoint from its image
$c(B_i)$\rom. (In fact, $B_i$ and $c(B_i)$ have opposite signs in
the sense of~\ref{tg-chessboard}.) Hence, the restriction of~$p$
to~$B_i$ is a one-to-one map onto a region
$\base_i\subset\Base\Gamma$. Since, in addition, the restriction
$p\:\Gamma'\to\Gamma$ is orientation preserving,
property~\itemref{tg}{tg-oriented} for~$\Gamma$ follows
from~\itemref{tg}{tg-oriented} for~$\Gamma'$.
%
The separation property implies
also that $\Gamma'$ contains the fixed point set $\Fix c$; this
yields~\itemref{tg}{tg-boundary} and~\ditto{tg-valency}
for~$\Gamma$.

Since $p$ preserves the decorations of $\Gamma'$ and $\Gamma$,
the admissibility of one
of the graphs
implies the admissibility of the other.
\endproof


The \emph{full valency} of a vertex of~$\Gamma$ is the valency of
any of its pull-backs in~$\Gamma'$. The full valency of an inner
vertex coincides with its valency; the full valency of a real
vertex equals $2\cdot\text{valency}-2$. The full valency of any
vertex is even.

In what follows, we denote by $\#_{\text{\white-}}(\Gamma)$,
$\#_{\text{\black-}}(\Gamma)$, and $\#_{\text{\cross-}}(\Gamma)$
the numbers of, respectively, \white--, \black--, and \cross-vertices
of $\Gamma'$. These numbers can be regarded as weighted
numbers of respective vertices in $\Gamma$, each inner vertex
being counted twice.

\paragraph\label{tg-maps}
A typical example of a trichotomic graph is the following. Let $B$ be
a connected closed surface with involution~$c$, and let
$j\:(B,c)\to(\Cp1,\spbar)$
be an equivariant ramified covering. (Thus, $B$ is necessarily
orientable, $c$ is orientation reversing, and one can assume~$j$
orientation preserving.)
Then $j$ defines a trichotomic graph $\Gamma'(j)\subset B$.
As a set, $\Gamma'(j)$ is the pull-back $j^{-1}(\Rp1)$.
The trichotomic graph structure on~$\Gamma'(j)$ is introduced
as follows: the \black--, \white--, and
\cross-vertices are  the pull-backs of~$0$, $1$,
and~$\infty$, respectively (monochrome vertices being the
branch points with other real critical values), the edges
are \solid-, \bold-, or \dotted- provided that their images belong to
$[\infty,0]$, $[0,1]$, or $[1,\infty]$, respectively, and
the
orientation of~$\Gamma'(j)$ is that induced from the positive
orientation of~$\Rp1$ (\ie, order of~$\R$).

\lemma\label{j->tg}
The graph~$\Gamma'(j)\subset B$ constructed above is an admissible
$c$-in\-vari\-ant trichotomic graph. Hence, its image
$\Gamma(j)=\Gamma'(j)/c\subset B/c=\base$ is an admissible trichotomic
graph.
\endlemma

\proof
By Lemma~\ref{j-double}, the second statement follows from the first one.
Axiom~\itemref{tg}{tg-oriented} for $\Gamma'(j)$
follows from the fact
that a region
$B_i\subset B_{\Gamma'(j)}$ is positive (negative) in the sense
of~\ref{tg-chessboard} if its image is the disk
$\{\Impart z\ge0\}$ (respectively, $\{\Impart z\le0\}$).
Other axioms are straightforward. The admissibility follows
from the fact that, since $j\:\Gamma'(j)\to\Rp1$ is orientation
preserving, $\prec$ is a subset of the partial order
induced by the linear orders on the intervals
$(1,\infty)$, $(\infty,0)$, and $(0,1)$.
\endproof

\theorem\label{tg->j}
Let $D$ be a compact connected surface and
$(B,c)$ its orientable double.
Exclude the case of oriented $D$ without boundary, and
equip $B$ with its canonical orientation. Then a
trichotomic graph $\Gamma\subset D$ is admissible if and only if it
has the form $\Gamma(j)$ for some orientation preserving
equivariant ramified covering $j\:(B,c)\to(\Cp1,\spbar)$.
Furthermore, $j$ is determined by~$\Gamma$ up to
homotopy in the class of equivariant ramified covering having
a fixed trichotomic graph.
\endtheorem

\proof
The `if' part is given by Lemma~\ref{j->tg}. For the `only if'
part, we will construct a map~$j$ and, at each step, check that the
construction is unique up to homotopy.

Any map~$j$ in question must have
an orientation
preserving descent $\bj\:\Gamma\to\Rp1$. The
images of the essential vertices are predefined, and the
extension of~$\bj$ to, say, the \solid- part of the graph is
determined, up to homotopy, by a
monotonous map
$\bj\:(\Gamma_{\msolid},\prec)\to((1,\infty),<)$.
The set of
such maps is defined by linear inequalities $\bj(u)<\bj(v)$ whenever
$u\prec v$, $u,v\in\Gamma_{\msolid}$. Hence, as a convex subset
of a Cartesian power of $(1,\infty)$, it is connected. For the
existence, one can, \eg, extend~$\prec$ to a linear order
(any maximal order) on~$\Gamma_{\msolid}$ and map the vertices to
consecutive integers.

Let $\Gamma'\subset B$ be the pull-back of~$\Gamma$,
see~\ref{tg-symmetric} for the notation.
The composition $\bj\circ p$ is an
equivariant orientation preserving map
$j\:(\Gamma',c)\to(\Rp1,\spbar)$. For each positive (in the sense
of~\ref{tg-chessboard}) region~$B_i$ of~$\Gamma'$,
the restriction $j\:\partial B_i\to\Rp1$ is a covering;
since the orientations on $B_i$ and $\partial B_i$ agree, $j$
extends to an orientation preserving ramified covering
$j\:B_i\to\{\Impart z\ge0\}$. Then
$\spbar\circ j\circ c\:c(B_i)\to\{\Impart z\le0\}$
extends~$j$ to the negative components.
The separation property of~$\Gamma'$ (see~\ref{j-double})
assures that the
extension $j\:B\to\Cp1$ is well defined and equivariant.
Each inner
point of an edge of~$\Gamma'$ is regular (as adjacent components
of~$B_{\Gamma'}$ have opposite signs); hence, $j$ has isolated
critical points and thus is a ramified covering.

The only ambiguity in the last step of the construction is in
extending a covering $\partial B_i\to\Rp1$ of the circle
to a ramified covering
$B_i\to\{\Impart z\ge0\}$ of the disk.
Any such extension can be perturbed to
a generic one (with all branch points double and all critical
values distinct), and the latter is unique up to homotopy due to
an analog of the Hurwitz
theorem
(see,
\eg,~\cite{Protopopov} or~\cite{Natanzon}; a very transparent
proof is indicated in~\cite{BE}).
\endproof

Theorem~\ref{tg->j} and the Riemann existence theorem result in the
following corollary.

\corollary\label{tg->an}
Given an admissible trichotomic graph $\Gamma\subset\base$, there
is
a complex structure on~$B$
and a holomorphic map
$j\:B\to\Cp1$ such that
the canonical orientation of $B$ coincides with its complex orientation,
$c$ is a real structure on $B$,
$j\:(B,c)\to(\Cp1,\spbar)$ is equivariant,
and $\Gamma=\Gamma(j)$.
Both the complex structure and the map are unique up to deformation.
\qed
\endcorollary

\remark{Remark}
As it
follows from the proof, a slightly stronger statement holds. On each of the
sets $\Gamma_{\msolid}$, $\Gamma_{\mbold}$, $\Gamma_{\mdotted}$
one can fix in advance a partial order extending~$\prec$. Then
$j$ can be chosen compatible with the given partial orders, and $j$
is unique up to homotopy in the class of such maps.
\endremark

\subsection{Deformations}\label{dfm}
Let us fix an oriented closed connected
surface~$B$
with an orientation reversing involution $c\:B\to B$. Let
$\base=B/c$ and let $p\:B\to\base$ be the projection. We are
interested in orientation preserving equivariant ramified
coverings $j\:(B,c)\to(\Cp1,\spbar)$.
A \emph{deformation} of coverings is a homotopy $B\times I\to\Cp1$
in the class of equivariant ramified coverings.
A deformation is called \emph{simple} if it preserves the
multiplicities of all the points with values $0$, $1$, and
$\infty$ and the multiplicities of
all branch points with real critical values.
Clearly, any deformation is locally simple with the exception of
finitely many
isolated values of the parameter $t\in I$. (As in
Section~\ref{tg}, we are working in the {\sl PL\/}-category; in
particular, this implies the finiteness.)
The
following statement
is an immediate consequence of Theorem~\ref{tg->j} and
the definition of~$\Gamma(j)$.

\proposition\label{dfm-simple}
Two equivariant ramified coverings $j_0,j_1\:B\to\Cp1$ can be
connected by a simple deformation if and only if their
graphs~$\Gamma({j_0})$, $\Gamma({j_1})$ are isotopic.
\qed
\endproposition

Let $\Gamma_0\subset\base$ be a trichotomic graph. Pick some disjoint
regular
neighborhoods~$U_v$ of all (or some) vertices~$v$
of~$\Gamma_0$ (we assume that $U_v\cap\partial\base=\varnothing$
unless $v$ is real)
and replace each intersection $\Gamma_0\cap U_v$ with another
decorated graph, so that the result~$\Gamma_1$ is again a trichotomic
graph. If each intersection $\Gamma_1\cap U_v$ contains essential
vertices of at most one kind, $\Gamma_1$ is called a
\emph{perturbation} of~$\Gamma_0$ (and $\Gamma_0$ is called a
\emph{degeneration} of~$\Gamma_1$).
A perturbation~$\Gamma_1$ of an
admissible trichotomic graph~$\Gamma_0$ is admissible if and only if
none of the intersections $\Gamma_1\cap U_v$ contains an oriented
monochrome cycle. (Note that there
are
no simple local criteria
for the admissibility of a degeneration.)

\remark{Remark}
Assume that $\Gamma_1$ is a perturbation of~$\Gamma_0$,
and $\Gamma_1\cap U_v$ contains no oriented monochrome
cycles. Since the intersection $\Gamma_1\cap\partial U_v$ is
fixed, the assumption on $\Gamma_1\cap U_v$ implies that
$\Gamma_1\cap U_v$ either is monochrome (if $v$ is monochrome) or
consists of monochrome vertices, essential vertices of the same
kind as~$v$, and edges of the two kinds incident to~$v$.
\endremark

Any deformation~$j_t$ of ramified coverings whose restriction to
$B\times(0,1]$ is simple results in a perturbation of the
graph $\Gamma_0=\Gamma({j_0})$. (The requirement that each
intersection $\Gamma_1\cap U_v$ should contain essential vertices
of at most one kind is due to the fact that essential vertices
have predefined distinct images in~$\Cp1$.) Our goal is to prove
the converse.

\proposition\label{dfm-perturbation}
Given an admissible graph~$\Gamma_0$ and its
admissible perturbation~$\Gamma_1$, there is a deformation
$j_t\:B\to\Cp1$,
$t\in[0,1]$, with the following properties\rom:
\roster
\item
one has $\Gamma_0=\Gamma({j_0})$ and $\Gamma_1=\Gamma({j_1})$\rom;
\item
the restrictions of all maps~$j_t$ to
$B\sminus\bigcup_v p^{-1}(U_v)$ coincide\rom;
\item
the restriction of the deformation to $B\times(0,1]$ is
simple.
\endroster
\endproposition

\proof
Let~$j_0$ be any ramified covering given by Theorem~\ref{tg->j}.
We can assume that the restriction of~$j_0$ to each
pull-back $U'_v=p^{-1}(U_v)$ has no branch points other than the pull-backs
of~$v$ itself. Then it suffices to construct a desired homotopy
(fixed on the boundary) on each pull-back $U'_v$.

Assume that $v$ is a \black-vertex, so that $j(v)=0$. (In the
other cases the proof is literally the same after reordering the
colors and a coordinate change in~$\Cp1$.) First, assume that $v$
is real. Let~$d$ be the full valency of~$v$. Regard~$U'_v$ as a
hemisphere in a sphere $\bar U'_v\cong S^2$
and extend both
$\Gamma'_0\cap U'_v$ and $\Gamma'_1\cap U'_v$ to
symmetric
trichotomic graphs $\bar\Gamma_0'$, $\bar\Gamma_1'$ on~$\bar U'_v$
by adding a real \cross-vertex~$\bar v$ of valency~$d$, $d$
\white-vertices of valency~$2$, and appropriate edges.
The graphs are admissible, and
Corollary~\ref{tg->an} gives real regular analytic maps
$f_0,f_1\:\bar U'_v=\Cp1\to\Cp1$ corresponding to~$\bar\Gamma_0'$,
$\bar\Gamma_1'$, respectively. Clearly, $f_0(z)=z^d$ and
$f_1(z)$ is a real polynomial of degree~$d$,
so that the family
$f_t(z)=t^df_1(z/t)$ is a desired homotopy. More precisely, we can
assume that all critical points of~$f_1$ other than~$\bar v$ are
mapped, say, to the disk $\{\mathopen|z\mathclose|<1/2\}$
(otherwise, replace~$f_1$ with some
$\varepsilon^df_1(\,\cdot\,/\varepsilon)$\,); then,
$f_t^{-1}\{\mathopen|z\mathclose|\le1/2\}$, $t\in I$,
is a disk bundle
over~$I$, and it can be identified with $U'_v\times I$ so that the
restriction of the homotopy to the boundary
$\partial U'_v\times I$ is constant.

If $v$ is not real, the same construction applies to one of the
two disks constituting~$U'_v$ (with~$c$ ignored) and
extends to the other disk by symmetry.
\endproof

Fix a set~$\Cal G$ of admissible trichotomic graphs closed
under isotopies, and let
$\Cal J$ be the set of equivariant ramified coverings $j\:B\to\Cp1$
defined via $j\in\Cal J$ if and only if $\Gamma(j)\in\Cal G$.

\corollary\label{deformations}
Let $j\:(B,c)\to(\Cp1,\spbar)$ be
an equivariant holomorphic map, $j\in\Cal J$. Assume that there is
a chain $\Gamma(j)=\Gamma^0,\Gamma^1,\ldots,\Gamma^n$ so that
$\Gamma^i\in\Cal G$, $i=0,\ldots,n$, and each $\Gamma^i$,
$i=1,\ldots n$, is a perturbation of, a degeneration of, or
isotopic to~$\Gamma_{i-1}$. Then there is a piecewise analytic
equivariant deformation~$j_t$, $t\in I$, of~$j=j_0$ such that all
$j_t\in\Cal J$, $t\in I$, and $\Gamma({j_1})=\Gamma_n$. Moreover,
each piece can be chosen as a closed real subinterval of an
equivariant deformation in the sense of Kodaira-Spencer over an
open complex disc. \rom(In general, the complex structure of~$B$
changes.\rom)
\endcorollary

\proof
Using Propositions~\ref{dfm-simple} and~\ref{dfm-perturbation},
one can construct a topological deformation $B\times
I\to\Cp1\times I$ as in the statement. By construction, the branch
set in $\Cp1\times I$ can be made piecewise analytic. Moreover, by
the choice made in the construction, for each (real closed) piece
the equivariant ramified covering extends to an open
complexification of the piece and the Grauert-Remmert
theorem applies to produce a complex structure.
\endproof

\subsection{Dessins}\label{sketches}
From now on, we
will only consider trichotomic graphs arising from the
$j$-invariants of almost generic elliptic surfaces
(see~\ref{counthominv})
or, more generally,
almost generic trigonal curves (see~\ref{topology}).
In view of~\ref{j.generic}, this is the case if and only if
\roster
\item"$(*)$"\local{star}
the full valency of each \cross--
(respectively, \white-- or \black--) vertex is~$2$ (respectively,
$0\bmod4$ or $0\bmod6$).
\endroster

\proposition\label{tg->curve}
Any admissible trichotomic graph satisfying~$(*)$ above is
of the form~$\Gamma(j)$, where $j\:B\to\Cp1$ is the $j$-invariant
of an
almost generic
real trigonal curve.
The latter
is determined uniquely up to deformation equivalence.
\endproposition

\proof
The
deformation
uniqueness
of an equivariant holomorphic map
$j\:(B,c)\to(\Cp1,\spbar)$ such that $\Gamma=\Gamma(j)$ is given
by Corollary~\ref{tg->an}. Let $G_3$, $G_2$, and~$I$ be,
respectively, the sum of all \white--, \black--, and
\cross-vertices considered as divisors on~$B$. By construction,
$2G_3$ is the zero divisor of~$j$, $3G_2$ is the zero divisor of
$j-1$, and $I$ is the pole divisor of both~$j$ and $j-1$. In
particular, $2G_3\sim3G_2$ and, hence, $G_2\sim2(G_3-G_2)$ and
$G_3\sim3(G_3-G_2)$. Thus, one can take for the bundle~$Y$
generating the ruled surface (see~\ref{trigonal-curves}) the line
bundle defined by the real divisor $G_3-G_2$.

Now, pick a real section $\tilde g_2\in\Gamma(B;\CO_B(Y\pos2))$ whose
zero divisor is~$G_2$ and a real section
$\tilde g_3\in\Gamma(B;\CO_B(Y\pos3))$ whose zero divisor is~$G_3$.
For $\alpha,\beta\in\R$ let $g_2=\alpha\tilde g_2$ and
$g_3=\beta\tilde g_3$. The sections
$4g_2^3j^{-1}$ and $27g_3^2(j-1)^{-1}$ of $\CO_B(Y\pos6)$
are regular and have the same zero divisor~$I$. Hence,
$\alpha$ and $\beta$ can be chosen so that $j$ is given
by~\eqref{jinv}.
They are defined up to the
transformation $(\alpha,\beta)\mapsto(t^2\alpha,t^3\beta)$,
$t\in\R$; the corresponding sections $g_2$, $g_3$ define
deformation equivalent
trigonal curve.
\endproof

\paragraph
Any graph satisfying~$(*)$ can be perturbed to a graph~$\Gamma$
such that
\roster
\item\local1
the full valency of each \cross--, \white--, or
\black-- vertex of~$\Gamma$ is, respectively,~$2$, $4$, or $6$;
\item\local2
the valency of any real monochrome vertex of~$\Gamma$ is~$3$;
\item\local3
$\Gamma$ has no inner monochrome vertices.
\endroster
An admissible graph
satisfying conditions~\loccit1--\loccit3
is called a \emph{dessin}; such a graph corresponds to a generic
trigonal curve. We always assume that the boundary of the
underlying surface is nonempty.
We freely extend to dessins
all
terminology that
applies to almost generic trigonal curves. Thus, we speak about
\emph{$(M-d)$-dessins}, \emph{\hbox{\rom(non-\rom)}hyperbolic}
(components of) dessins, \emph{ovals} and \emph{zigzags}
(see~\ref{tg->real}
for more details and a reinterpretation of these notions in
terms of the dessins).

The ramified covering defined by a dessin~$\Gamma$
has generic branching behavior; its degree is of the form~$6k$,
$k\in\Z$, and one has
$\#_{\text{\black-}}(\Gamma)=2k$, $\#_{\text{\white-}}(\Gamma)=3k$,
and $\#_{\text{\cross-}}(\Gamma)=6k$.
The number $3k$ is called the \emph{degree}
of~$\Gamma$. By definition, it is positive and divisible by~$3$.

A dessin~$\Gamma$ of degree~$3$ on a disk is
called a \emph{cubic}. Such a dessin~$\Gamma$ is indeed the dessin
of a nonsingular cubic curve in the projective plane blown-up at
one point.

Two dessins are called
\emph{equivalent} if, after a homeomorphism
of the underlying surfaces, they can be connected by a finite sequence of
isotopies and the following \emph{elementary moves}:
\roster
\item"--"
\emph{monochrome modification}, see
Figure~\ref{fig.moves}(a);
\item"--"
\emph{creating \rom(destroying\rom) a bridge}, see
Figure~\ref{fig.moves}(b),
a \emph{bridge} being a pair of
monochrome vertices connected by a real monochrome edge;
\item"--"
\emph{\white-in} and its inverse \emph{\white-out}, see
Figure~\ref{fig.moves}(c) and~(d);
\item"--"
\emph{\black-in} and its inverse \emph{\black-out}, see
Figure~\ref{fig.moves}(e) and~(f);
\endroster
(In the first two cases, a move is valid if and only if the result
is again a dessin, \ie, one needs to check its
admissibility.)

\midinsert
\eps{Fig1all}
\figure\label{fig.moves}
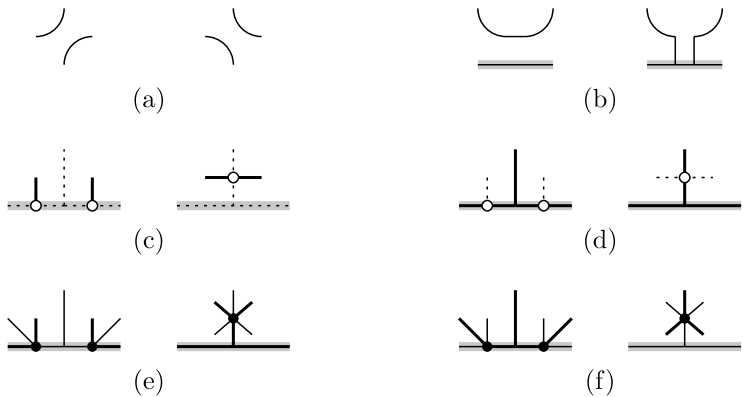
Elementary moves of dessins. Wide gray lines indicate
real points.
\endfigure
\endinsert

Clearly,
the elementary moves are exactly the results of passing
through
codimension~$1$ degenerations still satisfying~$(*)$. Hence,
in view of
Proposition~\ref{tg->curve} and Corollary~\ref{deformations}, the
following statement holds.

\proposition\label{equiv.curves}
Two generic real trigonal curves are deformation equivalent in the
class of
almost generic
real trigonal curves
if and only if their dessins are equivalent.
\qed
\endproposition

\paragraph\label{tg->topology}
The definition of the $j$-invariant gives an easy way to
recover the topology of a generic real trigonal curve
$C\subset\Sigma=\Cp{}(1 \oplus Y)$ from its dessin~$\Gamma$. Let
$q\:\Sigma\to B$ be the projection and $q_C$ its restriction
to~$C$. The pull-back
$q_C^{-1}(b)\subset q^{-1}(b)$ of each point
$b\in B\sminus\{\text{\cross-vertices}\}$ consists of three
points.
\roster
\item\local{b.region}
If $b$ is an inner point of a region of~$\Gamma$, the three points
of the pull-back $q_C^{-1}(b)$ form a triangle with all three
edges distinct. As a consequence, the restriction of~$q_C$ to the
interior of each region of~$\Gamma$ is a trivial covering.
\item\local{b.dotted}
If $b$ belongs to a \dotted- edge of~$\Gamma$, the three points of
the pull-back $q_C^{-1}(b)$ are collinear. The ratio
$(\text{smallest distance})/(\text{largest distance})$
is in $(0,1/2)$; it
tends to~$0$ ($1/2$) when $b$ approaches a \cross--
(respectively,~\white--) vertex.
\item\local{b.edge}
If $b$ belongs to a \solid- (\bold-) edge of~$\Gamma$, the three
points of the pull-back $q_C^{-1}(b)$ form an isosceles triangle
with the angle at the vertex less than (respectively, greater than)
$\pi/3$. The angle
tends to~$0$, $\pi/3$, or~$\pi$ when $b$ approaches, respectively,
a \cross--, \black--, or \white-vertex.
\endroster

\paragraph\label{tg->monodromy}
In particular, statements~\loccit{b.region}--\loccit{b.edge} above
give a very simple description of the $\Braids_3/\Delta^2$-valued monodromy,
see~\ref{monodromy},
along any loop~$\gamma$ in~$B^\sharp$. As a consequence, the
following
statements
hold:
\roster
\item"--"
if $\gamma$ does not intersect the closure of~$\Gamma_{\mdotted}$,
then the monodromy along~$\gamma$ is
determined by the corresponding permutation, which must be even
(as in this case the three points
in the fiber never become collinear);
\item"--"
in particular,
if $\gamma$ does not intersect the closure of
$\Gamma_{\mdotted}\cup\Gamma_{\msolid}$ (or the closure
of $\Gamma_{\mdotted}\cup\Gamma_{\mbold}$),
then the monodromy along~$\gamma$ is
trivial;
\item"--"
if $\gamma$ belongs to the closure of~$\Gamma_{\mdotted}$, then
the monodromy along~$\gamma$ is $\Delta^{\epsilon\bmod2}$,
where $\epsilon$ is the number of \white-vertices on~$\gamma$.
\endroster

\paragraph\label{tg->real}
Let $\Gamma\subset\base$ be a dessin.
The collection of all vertices and edges of~$\Gamma$
contained in a given connected component
of~$\partial\base$ is called a {\it real component}
of~$\Gamma$. In the drawings, (portions of) the real components
of~$\Gamma$ are indicated by wide grey lines.

Every
maximal \dotted-
segment on a non-hyperbolic real component
(respectively, every maximal real \bold- segment)
is bounded by
two
\cross--
(respectively, \black--) vertices.
(Here, segments are allowed to contain
monochrome vertices and \white-vertices.)
In particular, the numbers of
\cross-- and \black-vertices in each real component of~$\Gamma$
are even.

A real component of~$\Gamma$
(and the corresponding
component
of $\partial\base$) is
called
\roster
\item"--" \emph{even}/\emph{odd}, if it contains
an even/odd number of \white-vertices of~$\Gamma$,
\item"--" \emph{hyperbolic}, if all
edges of this component are dotted.
\endroster
In addition, define the \emph{parity} of each maximal \dotted-
segment of~$\Gamma$ and each complementary segment as the parity
of the number of \white-vertices contained in the segment.
Equivalence of dessins preserves their even, odd, and hyperbolic
components, as well as the parity of the segments.
A dessin
is called \emph{hyperbolic} if all its real components are
hyperbolic.

Now, let~$\Gamma$ be the dessin of a generic real trigonal
curve $C\subset\Sigma$ (see~\ref{tg->topology} for the notation).
Then, the real components $\Gamma_i$ of~$\Gamma$
are identified with the connected components~$B_i$
of~$B_\R$.
The pull-back $q_C^{-1}(b)$ of a real point $b\in\partial\base$
has three real points if $b$ is a \dotted- point or a
\white-vertex
adjacent to two real \dotted- edges; it has two real points, if $b$ is a
\cross-vertex, and a single real point otherwise.
A component $\Sigma_i$ of $\Sigma_\R$ is orientable
(equivalently, the restriction $Y_i$ of $Y_\R$ is topologically trivial,
see~\ref{topology}) if and only if the corresponding real
component $\Gamma_i$ is even.
(Indeed, recall that $Y$ is defined by the real divisor $G_3-G_2$,
see the proof of
Proposition~\ref{tg->curve}, and the restriction of~$G_2$ to $B_i$ is even.)

A component~$B_i$ is hyperbolic (in the sense of~\ref{topology})
if and only if so is~$\Gamma_i$.
If $B_i$ is non-hyperbolic,
its ovals and zigzags are represented
by the maximal \dotted- segments
of~$\Gamma_i$, even and odd, respectively.
The latter are also called \emph{ovals} and \emph{zigzags}.
Two consecutive ovals of~$\Gamma$ belong to a single chain,
see~\ref{topology}, if and
only if they are separated by an even number of \white-vertices.

\subsection{The oval count}\label{tg->ovals}
Let $\Gamma\subset\base$ be a dessin of degree $\deg\Gamma=3k$,
and let $C\subset\Sigma$ be the corresponding trigonal curve. Its
genus is $g(C)=3k-3\chi(\base)+1$.
Introduce the following notation:
\roster
\item"--"
$\Leven$, $\Lodd$: the numbers of even/odd
hyperbolic real components;
\item"--"
$\Lnh$: the number of non-hyperbolic real components;
\item"--"
$\Noval$, $\Nzigzag$, $\Ninner$: the numbers of ovals,
zigzags, and inner \cross-vertices, respectively;
\item"--"
$\delta=2-(\Leven+\Lodd+\Lnh)-\chi(\base)$:
the `excessive' Euler characteristic.
\endroster
Note that $2(\Noval+\Nzigzag+\Ninner)=6k$ is the weighted number of
\cross-vertices. Note also that all quantities introduced are
nonnegative and that $\Lnh>0$ unless $\Gamma$ is hyperbolic.
The following statement is an immediate
consequence of the discussion in~\ref{tg->real}.

\proposition\label{oval.count}
If $\Gamma$ is an $(M-d)$-dessin, one has
$$
2\Lnh+\Lodd+\Nzigzag+\Ninner+3\delta=d+4.\qed
$$
\endproposition

If $\Gamma$ is hyperbolic, one has $\Lnh=\Nzigzag=\Noval=0$
and $\Ninner = \deg\Gamma$, and
the identity in Proposition~\ref{oval.count} takes the form
$$
\Lodd+\deg\Gamma+3\delta=d+4.\eqtag\label{oval.count-hyp}
$$
As in this case one also has $\Lodd=\deg\Gamma\bmod2$,
the following statement holds.

\corollary\label{hyp.cong}
For a hyperbolic dessin, one has $d=\delta\bmod2$.
\qed
\endcorollary


\subsection{Inner \white-- and \black-vertices}\label{without-bw}
A dessin~$\Gamma$
is called \emph{bridge free} if
any bridge of~$\Gamma$ belongs to a monochrome
real
component, the
latter containing exactly two vertices.
A non-hyperbolic dessin~$\Gamma$ is called
\emph{almost connected} if
each connected component of~$\Gamma$ contains a non-hyperbolic
real component.

\lemma\label{connect}
Any dessin~$\Gamma$ is equivalent to a
bridge free
dessin~$\Gamma'$
with the same
numbers of essential inner vertices.
If, in addition, $\Gamma$ is
hyperbolic \rom(respectively, non-hyperbolic\rom),
then $\Gamma'$ can be
chosen connected
\rom(respectively, almost connected\rom).
\endlemma

\proof
Assume that $\Gamma$ has a bridge, and denote by~$\gamma$ the
intersection of the corresponding monochrome part of~$\Gamma$ and
the real component containing the bridge. If $\gamma$ is a whole
(monochrome) component containing more than~$2$ vertices, pick a
minimal (in the sense of~$\prec$) vertex~$v_0$, a vertex~$v_1$
adjacent to~$v_0$, and the other vertex $v_2\ne v_0$ adjacent
to~$v_1$; then destroying the bridge $[v_1,v_2]$ is an admissible
operation. Otherwise, $\gamma$ has a bridge $[v_1,v_2]$ adjacent
to an essential vertex of~$\Gamma$, and destroying $[v_1,v_2]$ is
also admissible.

Assume that the resulting dessin~$\Gamma$ is disconnected.
Consider a region~$R$ whose boundary
contains two circles~$\alpha_1$,
$\alpha_2$
in two different connected components $\Gamma_1$, $\Gamma_2$
of~$\Gamma$.
We need to show that~$\Gamma_1$ and~$\Gamma_2$
can be joined together provided that one of them,
say $\Gamma_2$, is hyperbolic.
Each of the circles $\alpha_1$, $\alpha_2$,
has edges of all three colors.
Furthermore, $\alpha_1$ has a \black-vertex and, hence, an inner
\solid- or \bold-edge $e_1$.
On the other hand, all real edges
of $\Gamma_2$ are \dotted-; hence, $\alpha_2$ has an inner edge~$e_2$
of the same color as~$e_1$.
The inner modification involving~$e_1$ and~$e_2$ is
admissible, it does not create bridges, and it reduces the number
of connected components of~$\Gamma$.
\endproof

The \emph{reduction} (a \emph{partial reduction}) of a trichotomic
graph~$\Gamma$ is the image $\tilde\Gamma\subset\tilde\base$
of~$\Gamma$ in the surface~$\tilde\base$ obtained from~$\base$ by
contracting all
(respectively, some)
monochrome real
components of~$\Gamma$.
The original graph~$\Gamma$ is called an \emph{inflation}
of~$\tilde\Gamma$.

The reduction carries a natural structure of a trichotomic graph.
(The image of a monochrome real component of~$\Gamma$ is a
monochrome vertex of~$\tilde\Gamma$ unless the resulting valency
is~$2$; in the latter case the image is ignored and considered
part of an edge. In some instances, the image of valency~$2$ is
retained as a marked point in~$\tilde\Gamma$.)
The reduction of a
bridge free
dessin~$\Gamma$ is a dessin
unless all real components of~$\Gamma$ are monochrome.
Furthermore,
if~$\Gamma$ is bridge free,
so is its reduction.
The reduction preserves the counts of inner/real essential
vertices of each type.
A dessin~$\Gamma$ is called \emph{reduced} if it
has no monochrome real components.
In this (and only this) case
$\Gamma$ coincides with its reduction.
A dessin is called \emph{totally reduced} if it
has no even real components without \cross-vertices.
A dessin~$\Gamma$ is totally reduced if and only if any dessin
equivalent to~$\Gamma$ is reduced. A dessin is equivalent to an
inflation of a totally reduced one if and only if it has an odd
component or a real \cross-vertex.

The following lemma is obvious (as
\white-- and \black-vertices can freely be `dragged' through the
marked points).

\lemma\label{partial}
Let $\Gamma$ be a bridge free dessin,
and let~$\tilde\Gamma$ be its \rom(partial\rom)
reduction. Then any dessin
equivalent
to~$\tilde\Gamma$ is a partial reduction of a dessin
equivalent to~$\Gamma$.
\qed
\endlemma

A dessin is called \emph{peripheral} if it has no inner vertices
other than \cross-vertices.

\proposition\label{bw}
Any non-hyperbolic dessin is equivalent to
a peripheral one.
\endproposition

\proof
Suppose that there exists a non-hyperbolic dessin not equivalent
to a dessin without inner \white-- and \black-vertices. Among such
dessins choose a dessin~$\Gamma$ with the smallest number of
essential inner vertices. According to Lemma~\ref{connect}, one
can assume~$\Gamma$ bridge free and almost connected, and, in view
of Lemma~\ref{partial}, it suffices to show that either~$\Gamma$
or its reduction~$\tilde\Gamma$ is equivalent to a dessin with
fewer inner vertices.

If all non-hyperbolic
real
components
of~$\Gamma$ are monochrome, then at least one such component is
adjacent to a \black-vertex, which must be inner, and a \black-out
move reduces the number of essential inner vertices. Otherwise,
the reduction~$\tilde\Gamma$ is a non-hyperbolic dessin
and we can replace~$\Gamma$
with~$\tilde\Gamma$, \ie, assume~$\Gamma$ reduced. Since $\Gamma$
is also bridge free,
any
nontrivial monochrome modification of~$\Gamma$ is
admissible.


Define an \emph{inner chain} (of length~$k$) in~$\Gamma$ as
a path $v_0$, \dots, $v_k$ in~$\Gamma$ such that
all edges $[v_i,v_{i+1}]$, $0\le i<k$, and all vertices~$v_i$,
$0<i<k$, are inner.

First, suppose that $\Gamma$ has an inner chain connecting an
inner \white-- or \black-vertex with a non-hyperbolic real
component. Let $v_0$, $v_1$, \dots, $v_k$ be a shortest inner
chain with this property, and assume that either $v_k$ is
monochrome or else no inner chain of length~$k$ connects an inner
\white-- or \black-vertex with a monochrome vertex at a
non-hyperbolic real component. In particular,
this assumption guarantees that the
creating a bridge modifications used below in the proof are
admissible.\footnote{In this proof, we are mainly interested in $\Gamma$ as an
abstract graph (i.e., regions do not matter), and the
modifications can be performed so as to keep condition
\iref{tg}{tg-triangle}
in the definition}

\Case0:
$v_k$ is monochrome and $v_{k-1}$ is a
\black-- or \white-vertex. Then the number of
inner
vertices is reduced
by a single \black-out (respectively, \white-out).

\Case1.1: $v_k$ is a \black-vertex and $v_{k-1}$ is a
\cross-vertex. Then $k\ge2$ and $v_{k-2}$ is a \white-vertex.
This case reduces to Case~0 by creating a \bold-bridge, see
Figure~\ref{fig5}.


\midinsert
\line{\vbox{\hsize=.5\hsize
\eps{Fig5}
\figure\label{fig5}
\endfigure}%
\vbox{\hsize=.5\hsize
\eps{Fig6}
\figure\label{fig6}
\endfigure}}
\endinsert

\Case1.2: $v_k$ is a \black-vertex and $v_{k-1}$ is a
\white-vertex. Consider the region~$R$ whose boundary includes
$[v_{k-1},v_k]$ and the inner \solid-edge incident to $v_k$, see
Figure~\ref{fig6}.
The vertex~$u$ following $v_k$, $v_{k-1}$ in
the boundary of~$R$ is a \cross-vertex.
If necessary, reduce~$R$ to a triangle by a monochrome
modification.
Then creating a bold bridge reduces this case to Case~0.

\Case2.1: $v_k$ is a \white-vertex and $v_{k-1}$ is a
\cross-vertex. Then $k\ge2$ and $v_{k-2}$ is a \black-vertex, and
creating a \bold-bridge reduces this case to Case~0, see
Figure~\ref{fig7}.


\midinsert
\line{\vbox{\hsize=.5\hsize
\eps{Fig7}
\figure\label{fig7}
\endfigure}%
\vbox{\hsize=.5\hsize
\eps{Fig8}
\figure\label{fig8}
\endfigure}}
\endinsert

\Case2.2: $v_k$ is a \white-vertex and $v_{k-1}$ is a
\black-vertex. If among the real neighbors of~$v_k$
(\ie, real vertices connected to~$v_k$ by a real edge)
there is a \cross-vertex,
creating a \solid-bridge reduces this case to Case~0, see
Figure~\ref{fig8}.
Otherwise, the real neighbors of~$v_k$ are
monochrome. Let~$a$ be one of them, and
let~$w$
be the
\white-vertex following~$a$ in the real component. Since the real
component is non-hyperbolic,
$w$
is distinct from~$v_k$.
Consider the region~$R$ whose boundary includes $[v_k,v_{k-1}]$
and $[v_k,a]$, see Figure~\ref{fig9}.
The vertex~$u$
following $v_k$, $v_{k-1}$ in the boundary of~$R$ is a
\cross-vertex. If necessary, reduce~$R$ to a
triangle
by a monochrome modification
and, if
$v_{k-1}$
and~$w$
are not adjacent, perform a monochrome modification to
create a \bold-edge
$[v_{k-1},w]$, see Figure~\ref{fig9}.
Now,
replace the original chain with $v_0$, \dots,
$v_{k-1}$,
$v_k'=w$.
Since the real component in question is
non-hyperbolic, iterating this procedure (in the same direction)
will produce a chain \dots, $v_{k-1}$, $v_k''$ with $v_k''$ having
a \cross-vertex as a real neighbor. This reduces the situation to
that considered at the beginning of this paragraph (Figure~\ref{fig8}).


\midinsert
\line{\vbox{\hsize=.5\hsize
\eps{Fig9}
\figure\label{fig9}
\endfigure}%
\vbox{\hsize=.5\hsize
\eps{Fig10}
\figure\label{fig10}
\endfigure}}
\endinsert

\Case3:
$v_k$ is monochrome and $v_{k-1}$ is a \cross-vertex.
Then $k\ge2$ and $v_{k-2}$ is a \white-- or
\black-vertex. By a monochrome modification one can create a
\bold-edge connecting $v_{k-2}$ with one of the real neighbors
of~$v_k$ and thus reduce this case to Case~1.2 (see Figure~\ref{fig10})
or~2.2 (see Figure~\ref{fig11}).



\midinsert
\line{%
\vbox{\hsize=.5\hsize
\eps{Fig11}
\figure\label{fig11}
\endfigure}%
\vbox{\hsize=.5\hsize
\eps{Fig18a}
\figure\label{fig18'}
\endfigure}%
}
\endinsert

Now, suppose that $\Gamma$ has no inner chain connecting an inner
\white-- or \black-vertex to a non-hyperbolic real
component. Note that any inner chain
connecting two hyperbolic real
components has a \black-vertex.
Since
$\Gamma$ is almost connected, one can find two inner chains
$C=(v_0,\ldots,v_k)$ and $C'=(v'_0,v'_1,\ldots)$
so
that $v_k$
belongs to a non-hyperbolic
real
component, $v_0$ and~$v_0'$ are
connected by a real edge in a hyperbolic
real
component, and $C'$
contains an inner \white-- or \black-vertex.
Observe that $k=1$ or~$2$, in the latter case $v_1$ being a
\cross-vertex. Denote by~$R$ the region incident to $[v_0,v_0']$.

\Case4: $v_0'$ is a \white-vertex. Then $v_0$ is monochrome and
$v_1'$ is an inner \black-vertex. If $k=1$, then $v_1$ is
monochrome, the vertex following $v_0$, $v_1$ in the boundary
of~$R$ is a real \cross-vertex, and
creating a \solid-bridge reduces this case to Case~0, see
Figure~\ref{fig18'}.
If $k=2$,
the reduction to Case~0 is obtained by creating a \solid-bridge
as in Figure~\ref{fig18} (if $v_2$ is monochrome) or
by creating a \bold-bridge as in
Figures~\ref{fig19} and~\ref{fig20}
(if $v_2$ is a \black-vertex and the \bold-edge
following~$v_2$ in the boundary of~$R$ is, respectively, inner or
real; in the former case, a \solid-inner modification is performed
first).


%

\midinsert
\line{%
\vbox{\hsize=.5\hsize
\eps{Fig18}
\figure\label{fig18}
\endfigure}%
\vbox{\hsize=.5\hsize
\eps{Fig19}
\figure\label{fig19}
\endfigure}%
}
\endinsert

\midinsert
\line{%
\vbox{\hsize=.5\hsize
\eps{Fig20}
\figure\label{fig20}
\endfigure}%
\vbox{\hsize=.5\hsize
\eps{Fig21}
\figure\label{fig21}
\endfigure}%
}
\endinsert


\Case5: $v_0'$ is monochrome. Then
$v'_1$ is a \cross-vertex, $v'_2$ is a \black-vertex,
$k=1$, and $v_0$ is a \white-vertex.
This case is reduced to Case~0 by creating a \bold-bridge
as in Figure~\ref{fig21} (if $v_1$ is monochrome) or
by creating a \solid-bridge as in
Figures~\ref{fig22}
and~\ref{fig23} (if $v_1$ is a \black-vertex and the \solid-edge
following~$v_1$ in the boundary of~$R$ is, respectively, real or
inner; in the latter case, a \bold-inner modification is required).
\endproof

\midinsert
\line{%
\vbox{\hsize=.5\hsize
\eps{Fig22}
\figure\label{fig22}
\endfigure}%
\vbox{\hsize=.5\hsize
\eps{Fig23}
\figure\label{fig23}
\endfigure}%
}
\endinsert

Next statement is an analogue of Proposition~\ref{bw} for
hyperbolic dessins.

\proposition\label{bw-hyp}
Any hyperbolic dessin is equivalent to a dessin whose all
\white-vertices are real.
\endproposition

\proof
As in Proposition~\ref{bw}, one can assume the dessin~$\Gamma$ in
question bridge free, connected, and reduced. Consider a shortest
inner chain $v_0$, \dots, $v_k$ connecting an inner
\white-vertex~$v_0$ with a real vertex~$v_k$. It is easy to see that
$k\le3$ and, since $\Gamma$ is bridge free, $k>1$.

If $k=2$, then $v_1$ is a
\black-vertex and $v_2$ is a \white-vertex, see Figure~\ref{fig23-1}.
Consider the region~$R$ as in the figure
and, if necessary, reduce it to a triangle by a monochrome
modification. Now, the number of inner \white-vertices is reduced
by creating a \dotted- bridge followed by a \white-out, see
Figure~\ref{fig23-1}.

\midinsert
\line{%
\vbox{\hsize=.5\hsize
\eps{Fig23-1}
\figure\label{fig23-1}
\endfigure}%
\vbox{\hsize=.5\hsize
\eps{Fig23-2}
\figure\label{fig23-2}
\endfigure}%
}
\endinsert

\midinsert
\eps{Fig23-3}
\figure\label{fig23-3}
\endfigure
\endinsert

If $k=3$, then either $v_1$ is a \cross-vertex,
$v_2$ is
a \black-vertex, and $v_3$ is a \white-vertex (see
Figure~\ref{fig23-2}), or $v_1$ is a \black-vertex,
$v_2$ is a \cross-vertex, and $v_3$ is monochrome
(see Figure~\ref{fig23-3}).
In the former case, all three \white-vertices
adjacent to $v_2$ are real (as otherwise
the chain $v_0$, \dots, $v_k$ would not be shortest),
and the number of inner \white-vertices is reduced
by creating a \dotted- bridge followed by a \white-out, see
Figure~\ref{fig23-2}.
In the latter case,
all three \white-vertices
adjacent to $v_1$ are inner,
and at least one of them (not necessarily $v_0$)
can be pushed out
by creating a \dotted- bridge followed by a \white-out, see
Figure~\ref{fig23-3}.
\endproof




\subsection{Indecomposable dessins}\label{normal-forms}
In this section, we allow dessins on disconnected
surfaces (which are merely unions of dessins on
the components of the surface).

Consider a dessin~$\Gamma\subset\base$. Let
$I_1,I_2\subset\partial D$ be a pair of segments whose endpoints are
not vertices of~$\Gamma$, and let
$\Gf\:I_1\to I_2$ be an isomorphism, \ie, a diffeomorphism of the
segments establishing a graph isomorphism
$\Gamma\cap I_1\to\Gamma\cap I_2$
and preserving the kinds of the vertices and edges. (Note that, if
$I_1$ contains at least one essential vertex of~$\Gamma$, then $\Gf$
necessarily preserves
the orientations of the edges
given by the trichotomic graph structure.)
Consider the
quotient $D_\Gf=D/\{x\sim\Gf(x)\}$ and the image
$\Gamma'_\Gf\subset D_\Gf$ of~$\Gamma$, and denote by~$\Gamma_\Gf$
the graph obtained from~$\Gamma'_\Gf$ by erasing the image
of~$I_1$, if $\Gf$ is orientation reversing, or converting the
images of the endpoints of~$I_1$ to monochrome vertices otherwise.

In what follows we always
assume that either $I_1$ is part of an edge of~$\Gamma$ or
$I_1$ contains a single \white-- or \cross-vertex. In the latter case, $\Gf$
is unique up to isotopy; in the former case, $\Gf$ is determined
by whether it is orientation preserving or orientation reversing.
If $\Gamma_\Gf$ is a dessin, it is called the result of
\emph{gluing} $\Gamma$ along~$\Gf$.
(Sometimes we speak about gluing several dessins, meaning gluing
their disjoint union.)
The image of~$I_1$ is called a \emph{cut}
in~$\Gamma_\Gf$, and $\Gamma$ is called the result of a cut. The cut is
called \emph{genuine} (\emph{artificial}) if $\Gf$ is orientation
preserving (respectively, reversing); it is called a \solid-,
\dotted-, \bold-, or \cross-cut according to the structure of
$\Gamma\cap I_1$.
(The terms \dotted- and \bold- still apply
to cuts containing a \white-vertex.)

A dessin that is not equivalent to the result of gluing another
dessin is called \emph{indecomposable}.
A \emph{generalized cubic}
is
a dessin whose reduction is a cubic.

\theorem\label{reduction}
Any indecomposable dessin is a
disjoint
union of generalized cubics.
\endtheorem

\corollary\label{corollary}
Any
dessin can be obtained
from a disjoint union of generalized cubics by a sequence of gluing
operations and equivalences.
\qed
\endcorollary

\remark{Remark}
At present, we do not know whether a given graph is equivalent to
the result of gluing of a union of cubics.
As shown below, this is
true for $M$- and $(M-1)$-dessins.
\endremark

In view of Propositions~\ref{bw} and~\ref{bw-hyp},
Theorem~\ref{reduction} is an immediate consequence of
Propositions~\ref{reduction-hyp}
(the hyperbolic case)
and~\ref{reduction-details}
(the non-hyperbolic case).


\proposition\label{reduction-hyp}
Let $\Gamma$ be a connected reduced hyperbolic dessin
whose all \white-vertices
are real. Then
$\Gamma$ either is a cubic, or
has a cut; in the former case, $\Gamma$
is isotopic to the dessin shown in Figure~\ref{fig-star-like}.
\endproposition

\proof
Consider a \black-vertex~$v$ of~$\Gamma$.
Under the hypothesis, $v$ has a neighborhood shown
in Figure~\ref{fig-star}. If this neighborhood
does not close up to a cubic
({\it i.e.}, at least one of the regions adjacent
to~$v$ is not a triangle),
then $\Gamma$
has an artificial dotted cut
(located in the above region).
\endproof

\midinsert
\line{%
\vbox{\hsize=.5\hsize
\eps{Fig-star}
\figure\label{fig-star}
\endfigure}%
\vbox{\hsize=.5\hsize
\eps{Fig-starlike}
\figure\label{fig-star-like}
\endfigure}%
}
\endinsert

\remark{Remark}
One can show that, on the disc, any two hyperbolic
dessins of the same degree
are equivalent.
If all \white-vertices are real,
such a dessin~$\Gamma$ is a perturbation of
a star-like trichotomic graph as in Figure~\ref{fig-star-like},
with $2\deg\Gamma$ alternating rays radiating from
a single multiple \black-vertex.
(Note that the latter graph does satisfy~\ref{sketches}$(*)$,
and thus represents
the $j$-invariant of an almost generic curve,
see Proposition~\ref{tg->curve}.)
\endremark

\proposition\label{reduction-details}
Let $\Gamma$ be a reduced peripheral dessin on
a connected surface. Then
either $\Gamma$ is a cubic, or $\Gamma$ is
equivalent to a peripheral dessin with a cut.
\endproposition

Proposition~\ref{reduction-details}
is a mere combination of Lemmas~\ref{regions-list}
and~\ref{regions-cubic} proved at the end of this section.

Given a region~$R$, a component of the boundary $\partial R$ is
called a \emph{$3m$-gonal component} if it contains $3m$
essential vertices (equivalently, $m$ vertices of any given kind).
If $\partial R$ consists of a single $3m$-gonal
component, then $R$ itself is called a \emph{$3m$-gon}.

Recall that the real \white-vertices of a dessin can be subdivided
into two types, depending on the type of the real edges incident
to the vertex. Similarly, the real \black-vertices in the boundary
of a given region~$R$ can be subdivided into three types,
depending on which of the three angles at the
vertex belongs to~$R$.

\lemma\label{regions}
Let~$R$ be a region in a reduced peripheral indecomposable dessin.
Then the following holds\rom:
\roster
\item\local1
the boundary $\partial R$ cannot contain two distinct real edges of
the same kind\rom;
\item\local2
the boundary $\partial R$ cannot contain two distinct \white-vertices
of the same type\rom;
\item\local3
the boundary $\partial R$ cannot contain two distinct \black-vertices
of the same type\rom;
\item\local4
the boundary $\partial R$ consists of either one or two triangles or
a hexagon\rom;
\item\local5
unless $R$ is a triangle, the boundary $\partial R$ cannot contain
an inner \cross-vertex adjacent to a \solid-
monochrome vertex\rom;
\item\local6
if $\partial R$ is disconnected, it cannot contain
a real \cross-vertex.
\endroster
\endlemma

\proof
If $\partial R$ contains two real edges of the same kind, they either
are connected by an inner edge of the same kind or can
be connected by an artificial
cut; in both cases the graph is decomposable. This proves~\loccit1.
Statement~\loccit2 follows directly from~\loccit1, and \loccit3
follows from~\loccit1 unless $R$ has no real edges at the two
vertices in question. In the latter case, a \bold- inner
modification results in a region with two distinct \solid- real
edges, which contradicts~\loccit1. (Alternatively, a \solid- inner
modification results in a region with two distinct \bold- real
edges.)

In view of~\loccit2, $\partial R$ contains at most two
\white-vertices. This implies~\loccit4.

Let~$u$ be a \cross-vertex as in~\loccit5. In view of~\loccit2,
since $\partial R$ is not a triangle, it contains a \white-vertex
incident to \dotted- real edges.
Then, creating a \dotted-bridge produces a cut
(containing~$u$).

Let~$u$ be a real \cross-vertex in~$\partial R$ and let~$v$ be a
\cross-vertex in another component of~$\partial R$. Due
to~\loccit1, $v$ is an inner vertex, and one can create a
\solid-bridge (close to~$u$), converting~$R$ to a hexagon and~$v$,
to a \cross-vertex as in~\loccit5.
\endproof

In Lemma~\ref{regions-list} below
we list all regions appearing in
an indecomposable dessin
(see Figures~\ref{fig24} and~\ref{fig25}).
Various brackets in the notation indicate
the `ends' of a region, \ie, the components of the inner parts of
its boundary. (Clearly, it is these components that govern the
adjacencies of the regions.) The symbols $\lbold$, $\lfloor$, and
$\lceil$ (and the corresponding right delimiters) stand,
respectively, for a \bold-, \solid-, and \dotted- edge, and
$\lbrack$ stands for a pair of edges separated by an inner
\cross-vertex. The brace $\{$ indicates several `ends' that are
not of particular interest, and $($ indicates no `end' at all.

\midinsert
\labeldelta=.125cm
\def\Y{17}
\eps{Fig-trg}
\centerline{%
\labelat[-36,\Y]{$(1\rbold$}%
\labelat[-18,\Y]{$\lbrack2_1\rbold$}%
\labelat[0,\Y]{$\lceil2_2\rbold$}%
\labelat[18,\Y]{$\lbrack3_1\rbold$}%
\labelat[36,\Y]{$\lfloor3_2\rbold$}%
\labelat[-36,-1]{$\lbrack4_1\rbold$}%
\labelat[-18,-1]{$\lbrack5_1\rbold$}%
\labelat[0,-1]{$\lfloor5_2\rbold$}%
\labelat[18,-1]{$\lfloor6\rceil$}%
\labelat[36,-1]{$\lbrack7)$}%
}
\figure\label{fig24}
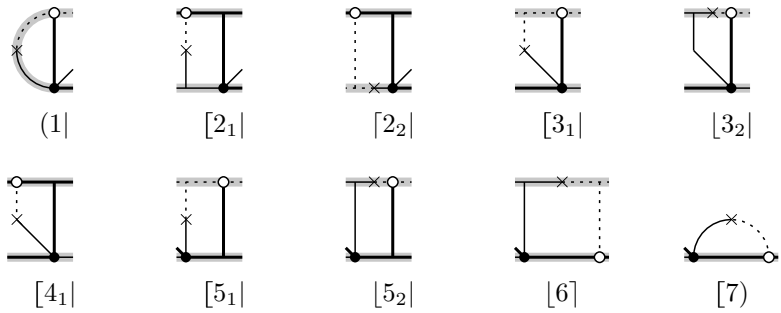
Triangular regions of indecomposable dessins
\endfigure
\bigskip
\bigskip
\eps{Fig-xg}
\centerline{%
\labelat[-28,-1]{$\{A_1\rbrack$}%
\labelat[0,-1]{$\{A_2\rceil$}%
\labelat[28,-1]{$\{B\rbold$}%
}
\figure\label{fig25}
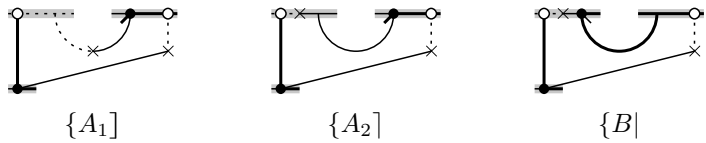
Hexagonal regions of indecomposable dessins
\endfigure
\bigskip
\bigskip
\eps{Fig-ex}
\figure\label{fig26}
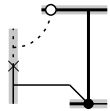
The exceptional triangle
\endfigure
\endinsert

\lemma\label{regions-list}
Any region~$R$
in a reduced peripheral indecomposable dessin
is either one of the
triangles
in Figure~\ref{fig24}
or one of the hexagons
in Figure~\ref{fig25}.
\endlemma

\proof
Lemma~\iref{regions}1
restricts all possible
triangle components of the boundary of~$R$ to those
listed in Figures~\ref{fig24} and~\ref{fig26},
and \ditto1--\ditto3 and~\ditto5
restrict the hexagons to those listed in Figure~\ref{fig25}.
Furthermore, a hexagon bounds a region, see~\iref{regions}4, and
if the latter is not a disk, it can be modified to a region with
disconnected boundary, see below.

Assume that $R$ is the triangle in Figure~\ref{fig26}.
Its \bold- edge can
only be adjacent to
a triangle of type~$2_1$ or~$2_2$ or a
hexagon of type~$B$. In the former case, a \white-in modification
followed by a \white-out along any \dotted-edge produces a
\dotted- cut. In the latter case, a \bold-inner modification
within the hexagon results in a region with two \solid- (as well
as two \dotted-) real edges.

Finally, assume that $\partial R$ consists of two triangles.
Lemma~\iref{regions}6 reduces the list of triangles to $2_1$,
$3_1$, $4_1$, $5_1$, and~$7$, and \iref{regions}5
eliminates~$2_1$. Thus, in view of~\iref{regions}1, the
boundary~$\partial R$ must be formed by one of the pairs $3_1$,
$4_1$ or $3_1$, $7$. The former is eliminated by~\iref{regions}3,
and in the latter case, a \solid- (or \dotted-) inner modification
results in a region with two \bold- real edges.
\endproof

\lemma\label{regions-cubic}
Any reduced dessin \rom(on a connected surface\rom)
whose regions are those listed in Lemma~\ref{regions-list}
is a cubic. Conversely, all regions of a
peripheral
cubic
are among those listed in Lemma~\ref{regions-list},
and they
are attached to one another according to one of the following
adjacency schemes \rom(see Figure~\ref{fig.cubics}\rom)\rom:
\smallskip
{\openup3pt
\halign{\kern2\parindent\hbox to2em{\hss$#$\rom:\ }&&\hss$\quad#$\hss\cr
\II_1&(1\rbold&@---&\lbold3_1\rbrack&@---&\lbrack5_1\rbold
 &@---&\lbold5_1\rbrack&@---&\lbrack3_1\rbold&@---&\lbold1)\cr
\I_1&(1\rbold&@---&\lbold3_1\rbrack&@---&\lbrack5_1\rbold
 &@---&\lbold5_2\rfloor&@---&\lfloor3_2\rbold&@---&\lbold1)\cr
\I_1&(1\rbold&@---&\lbold3_2\rfloor&@---&\lfloor5_2\rbold
 &@---&\lbold5_1\rbrack&@---&\lbrack3_1\rbold&@---&\lbold1)\cr
\II_3&(1\rbold&@---&\lbold3_2\rfloor&@---&\lfloor5_2\rbold
 &@---&\lbold5_2\rfloor&@---&\lfloor3_2\rbold&@---&\lbold1)\cr
\I_2&(1\rbold&@---&\lbold3_2\rfloor&@---&\lfloor6\rceil
 &@---&\lceil6\rfloor&@---&\lfloor3_2\rbold&@---&\lbold1)\cr
\I_1&(1\rbold&@---&\lbold3_2\rfloor&@---&\lfloor6\rceil
 &@---&\lceil2_2\rbold&@---&\lbold4_1\rbrack&@---&\lbrack7)\cr
\noalign{\smallskip}
\II_0&(7\rbrack&@---&\lbrack4_1\rbold&@---&\lbold2_1\rbrack&
 @---&\lbrack2_1\rbold&@---&\lbold4_1\rbrack&@---&\lbrack7)\cr
\I_0&(7\rbrack&@---&\lbrack4_1\rbold&@---&\lbold2_2\rceil&
 @---&\lceil2_2\rbold&@---&\lbold4_1\rbrack&@---&\lbrack7)\cr
}\smallskip\noindent
\halign{\kern2\parindent\hbox to2em{\hss$#$\rom:\ }&&\hss$\quad#$\hss\cr
\I_0&\{A_1\rbrack&@---&\lbrack3_1\rbold&@---&\lbold1)\cr
\II_2&\{A_2\rfloor&@---&\lfloor3_2\rbold&@---&\lbold1)\cr
\II_1&\{B\rbold&@---&\lbold4_1\rbrack&@---&\lbrack7)\cr
}}\smallskip\noindent
\rom(in the last three cases each hexagon being also adjacent to a
triangle of type~$1$ and a triangle of type~$7$\rom).
\endlemma

\remark{Remark}
Some pairs of dessins listed in Lemma~\ref{regions-cubic}
and Figure~\ref{fig.cubics} are equivalent.
It is easy to see that, in fact, there are seven equivalence
classes of cubics. They differ by the type ($\I$ or~$\II$;
equivalently, cubics with an oval are of type~$\I$, and those
without ovals are of type~$\II$) and the
number of zigzags (shown as a subscript in the notation). The
equivalence
class represented by each dessin is also listed in
Lemma~\ref{regions-cubic}
and Figure~\ref{fig.cubics}.
\endremark

\midinsert
\labeldelta=.125cm
\def\Y{\the\count0}
\def\X{-\XX}
\def\XX{18}
\eps{Fig-cubics}
\centerline{%
\count2=18
\count0=\count2
\multiply\count0 by5
\advance\count0 -1
\advance\count0-\count2
\labelat[\X,\Y]{$\II_1$}%
\labelat[\XX,\Y]{$\I_1$}%
\advance\count0-\count2
\labelat[\X,\Y]{$\II_3$}%
\labelat[\XX,\Y]{$\I_2$}%
\advance\count0-\count2
\labelat[\X,\Y]{$\I_1$}%
\labelat[\XX,\Y]{$\II_0$}%
\advance\count0-\count2
\labelat[\X,\Y]{$\I_0$}%
\labelat[\XX,\Y]{$\I_0$}%
\advance\count0-\count2
\labelat[\X,\Y]{$\II_2$}%
\labelat[\XX,\Y]{$\II_1$}%
}
\figure\label{fig.cubics}
Peripheral cubic dessins
\endfigure
\endinsert

\proof
It suffices to consider a dessin whose all regions are among those
listed in Lemma~\ref{regions-list}. (Any reduced
non-hyperbolic
cubic has
this property since it is indecomposable.)
Comparing the `ends' of the regions,
one arrives at the following list of
adjacencies:
\smallskip
{\openup3pt
\halign{\kern2\parindent\hss$#$&\hss$\quad#\quad$\hss&\hss$#$\hss&\hss$\quad#\quad$\hss&$#$\hss\cr
&&(1\rbold&@---&\rbold3_1\rbrack,\rbold3_2\rfloor,\{A_1\rbrack,\{A_2\rfloor,\{B\rbold\cr
\lbold2_1\rbrack&@---&\lbrack2_1\rbold&@---&\lbold4_1\rbrack\cr
\lbold2_2\rceil,\lfloor6\rceil&@---&\lceil2_2\rbold&@---&\lbold4_1\rbrack\cr
\lbold5_1\rbrack,\{A_1\rbrack&@---&\lbrack3_1\rbold&@---&\lbold1)\cr
\lbold5_2\rfloor,\lceil6\rfloor,\{A_2\rfloor&@---&\lfloor3_2\rbold&@---&\lbold1)\cr
(7\rbrack&@---&\lbrack4_1\rbold&@---&\lbold2_1\rbrack,\lbold2_2\rceil,\lbold B\}\cr
\lbold3_1\rbrack&@---&\lbrack5_1\rbold&@---&\lbold5_1\rbrack,\lbold5_2\rfloor\cr
\lbold3_2\rfloor&@---&\lfloor5_2\rbold&@---&\lbold5_1\rbrack,\lbold5_2\rfloor\cr
\lbold3_2\rfloor&@---&\lfloor6\rceil&@---&\lceil2_2\rbold,\lceil6\rfloor\cr
&&(7\rbrack&@---&\lbrack4_1\rbold,\{A_1\rbrack,\{A_2\rfloor,\{B\rbold\cr
}}
\smallskip\noindent
It remains to list all
chains of regions joined according to these rules, terminating
a chain whenever there are no free `ends' left.

Assume that all regions of~$\Gamma$ are triangles. If $\Gamma$ has
a triangle of type~$1$, starting from it one obtains one of the
first six schemes in the statement. Otherwise, $\Gamma$ has no
triangle of types~$3_1$, $3_2$ and, hence, no triangle of
types~$5_1$, $5_2$, or~$6$. Assuming that $\Gamma$ has a triangle
of type~$7$, one arrives at the last two schemes with triangles
only. Otherwise, $\Gamma$ has no triangle of type~$4_1$ and, hence, no
triangle of type~$2_1$ or~$2_2$, \ie, such a dessin does not exist.
Finally, any hexagon that $\Gamma$ may have extends uniquely to
one of the last three schemes in the statement. It is
straightforward to observe that all eleven schemes do represent
cubics.
\endproof

\subsection{Scraps}\label{scraps}
Given a dessin and one or
several
of
its inner edges, each connecting a real \white-vertex and
a real monochrome vertex,
one can
cut the dessin along these edges;
the connected components of the result (which, in general, is not
a dessin anymore) are called \emph{scraps}.
The edges used in the cut are called \emph{breaks}; they
can be \dotted- or \bold-. (In the sequel we need \dotted-
breaks only.) Note
that a scrap with breaks is not a dessin; it can be regarded as a
`dessin with boundary.' Two scraps can be glued along a break of the
same kind. The result is a dessin if and only if it is admissible
and has no breaks.

We extend to scraps the weighted numbers
$\#_{\text{\white-}}$, $\#_{\text{\black-}}$, and $\#_{\text{\cross-}}$.
Given a scrap~$\sigma$ on a surface~$D$, denote by
$\beta(\sigma)$ the number of breaks in
the boundary of~$\sigma$. Let further
$\kappa(\sigma)=\chi(D)-\frac12\beta(\sigma)$.
The latter quantity is additive; one has $\kappa(\sigma)>0$
if and only if $D$ is
a disk and
$\beta(\sigma)\le1$,
and $\kappa(\sigma)=0$ if and only if $D$
is a disk and $\beta(\sigma)=2$, or $D$ is an annulus or a
M\"{o}bius band and $\beta(\sigma)=0$.
Another additive quantity associated to a scrap~$\sigma$
is the degree
$\deg(\sigma) = \#_{\text{\white-}}(\sigma)
-\frac12\beta(\sigma)$.
The degree of any scrap is positive.

\lemma\label{scrap-degree}
For a scrap~$\sigma$, one has
$\#_{\text{\cross-}}(\sigma)=2\deg(\sigma)$ and
$\#_{\text{\black-}}(\sigma)=\frac{2}{3}\deg(\sigma)$.
Furthermore, $\deg(\sigma)+\frac{3}{2}\beta(\sigma)=0 \bmod 3\Z$.
\endlemma

\proof
It
suffices
to complete~$\sigma$ to a true dessin by patching
each break with a half of a cubic (say,
$\lbold5_1\rbrack@---\lbrack3_1\rbold@---\lbold1)$ or
$\lceil6\rfloor@---\lfloor3_2\rbold@---\lbold1)$ in the notation
of section~\ref{normal-forms}) and to use the known
identities and congruences for dessins.
\endproof

\corollary\label{kappa>=0}
A scrap~$\sigma$ with $\beta(\sigma)=1$
\rom(respectively,~$2$\rom)
has $\deg(\sigma)\ge\frac{3}{2}$
\rom(respectively, $\deg(\sigma)\ge3$\rom).
\qed
\endcorollary

\paragraph\label{breaks}
The importance of scraps is in the following construction.
Let~$\Gamma$ be a dessin (or a scrap). Each oval and each odd
real
hyperbolic component of~$\Gamma$ has
at least one \dotted- monochrome vertex~$u$. Let~$e$ be the inner
edge incident to~$u$ (and extended through any inner \white-vertex),
and let~$v$ be the other end of~$e$. Then either $v$ is an inner
\cross-vertex, or $v$ is a monochrome vertex and, hence, $e$ is a
\dotted- cut, or else $v$ is a real \white-vertex.
In the last case, $e$ has no inner vertices, and, thus,
breaks~$\Gamma$ into smaller scrap(s).

As an immediate consequence, since
a
monochrome vertex~$u$ in an
odd
real
hyperbolic
bridge free
component cannot be adjacent to a \white-vertex not
in the component, we obtain the following statement.

\lemma\label{Lodd<Ninner}
If a dessin~$\Gamma$ has no genuine \dotted- cuts,
then
$\Lodd\le\Ninner$.
\qed
\endlemma

\lemma\label{kappa>0}
A scrap~$\sigma$ with $\kappa(\sigma)>0$ contains a
zigzag or
an
inner \cross-vertex.
\endlemma

\proof
If $\sigma$ has no inner \cross-vertices, one can use~\ref{breaks}
to subdivide it into smaller scraps so that none of them has
ovals. At least one of the pieces still has $\kappa>0$.
Such a piece $\sigma'$ can only be a scrap
on a disk with $\beta(\sigma')\le1$.
Due to
Lemmas~\ref{kappa>=0}
and~\ref{scrap-degree}
it has at least three \cross-vertices and,
hence, at least one zigzag.
\endproof

\theorem\label{6k<=}
If an $(M-d)$-dessin~$\Gamma$ has no genuine \dotted-
cut, then
$$
2\deg\Gamma\le3(\Nzigzag+\Ninner)+3d-3\delta.
$$
\endtheorem

\proof
Let $\deg\Gamma=3k$. Using the construction of~\ref{breaks}, one
can break~$\Gamma$ into scraps, the total number of breaks
being~$2b$, where
$$b\ge b_0=\Lodd+\Noval-\Ninner=\Lodd+3k-\Nzigzag-2\Ninner
\eqtag\label{temp1}$$
(we count each break twice,
once in each of the two scraps incident to it).
Let $m_+$ be the nubmer of scraps~$\sigma$ with $\kappa(\sigma)>0$.
Using Lemma~\ref{scrap-degree} one can split $m_+=m_+'+m_+''$,
where $m_+'$ is the number of scraps with $\deg(\sigma)=\frac{3}{2}$
and $m_+''$ is the number of scraps with $\deg(\sigma)\ge\frac{9}{2}$.
According to Lemma~\ref{kappa>0},
at least $m_+' - \Nzigzag$ inner \cross-vertices
are separated by breaks from the ovals,
and the inequality~\eqref{temp1}
can be sharpened to $b\ge b_0+m_+'-\Nzigzag$.

Let $b_-$ be the
total number of breaks in the scraps
with $\beta\ge3$.
Then, according to Corollary~\ref{kappa>=0} and the definition
of $m_+''$, the number $3k$ of
\white-vertices of~$\Gamma$ is at least
$b+(2b-b_-)+3m_+''\ge3\Lodd+9k-6(\Nzigzag+\Ninner)-b_-+3m_+$. Hence, one must
have $6k\le6(\Nzigzag+\Ninner)-3\Lodd+b_--3m_+$.
On the other hand,
since $\kappa/\beta\le-\frac{1}{6}$
for a scrap with $\beta \ge 3$,
the additivity of~$\kappa$ yields
$\frac16b_-\le\frac12m_+-\chi(D)$.
Hence,
$6k\le 6(n_z+n_i)-3\Lodd-6\chi(D)$,
and it remains to substitute $\chi(\base)=2-(\Lodd+\Lnh)-\delta$ and
use Proposition~\ref{oval.count}.
\endproof

\end


\ifx\defloaded\undefined\input elliptic.def \fi

\begin

\setcounter\section5


\section{Applications: $M$- and $(M-1)$-cases}\label{applications}

Recall that we only consider generic curves and
surfaces and classify them
up to deformation in the class of almost generic ones.
The base of the fibration is never assumed fixed; it is also
subject to a deformation. Unless stated otherwise,
trigonal curves never intersect
the exceptional section.

\subsection{Junctions}\label{junctions}
Define a \emph{\rom(self-\rom)junction} as a genuine gluing of
a dessin along
isomorphic parts of two zigzags (respectively, one zigzag)
so that the resulting cut connects two ovals
(respectively, an oval and an odd hyperbolic component)
of the dessin obtained, see~\ref{normal-forms} for the
terminology concerning cuts.
Note that a (self-)junction consumes the zigzags involved.
In particular, any two junctions commute.

Below we state several structure theorems that deal
with junctions of cubics. From the point of view of the junction
operation, a cubic can be regarded as a `black box'
with a certain extra decoration of its boundary.
More precisely, we define a \emph{ribbon box}
as a disc with
a few disjoint
segments (called \emph{dotted})
marked in the boundary,
and a parity assigned to each dotted and each complementary
segment;
a box is required to be one of those listed
in Figure~\ref{blocks}.

\midinsert
\labeldelta=.125cm
\def\Y{20}
\eps{S6block}
\centerline{%
\labelat[-32,\Y]{$\I_0$}%
\labelat[-11,\Y]{$\I_1$}%
\labelat[11,\Y]{$\I_2$}%
\labelat[-32,-1]{$\II_0$}%
\labelat[-11,-1]{$\II_1$}%
\labelat[11,-1]{$\II_2$}%
\labelat[32,-1]{$\II_3$}%
}

\figure\label{blocks}
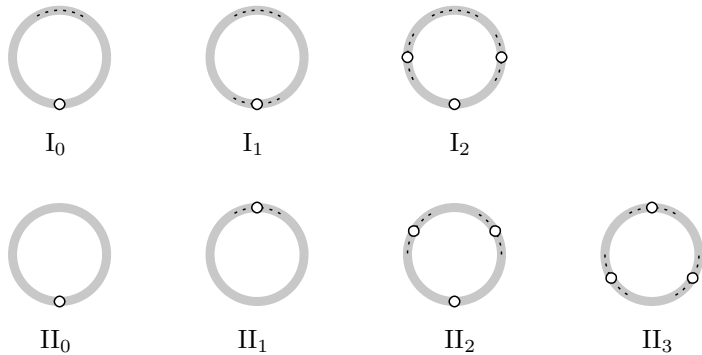
Ribbon boxes; odd segments are marked with \white-points
\endfigure
\endinsert

Each non-hyperbolic cubic dessin $\Gamma$ gives rise to a ribbon box:
the disk is the underlying surface of~$\Gamma$,
the dotted segments of the box are the maximal real dotted
segments of~$\Gamma$, and the parity is given by the number
of \white-vertices, as in the cubic, \cf.~\ref{tg->real}. Conversely,
in view of the classification given by Lemma~\ref{regions-cubic},
each box is obtained in this way from a cubic
dessin which is unique up to equivalence
leaving real \cross-vertices fixed.
In view of this correspondence, we will refer
to even (odd) dotted segments of a ribbon box
as oval (respectively, zigzag) segments.

A \emph{ribbon curve structure} is a collection
of boxes in which
some of the boxes are glued via identifying certain pairs of
zigzag segments, so that the result is a connected surface,
and
some of the remaining zigzag segments are selected for future
self-junction.
The selected segments are called {\it vanishing}.
An \emph{isomorphism} of two ribbon curve structures
is a homeomorphism of the underlying surfaces
preserving the decorations, {\it i.e.},
taking boxes to boxes
and
dotted segments to dotted segments,
preserving all parities,
and
taking
vanishing segments to vanishing segments.

Each of the two zigzag segments of a box of type~$\I_2$ or~$\II_2$
can be given a preferred orientation, say, towards
its odd complementary neighbor.
Thus, each adjacency of two such boxes
has a sign: it is said to be \emph{positive}
or \emph{negative} depending on whether the orientations
of two segments involved do or do not coincide.

The \emph{junction graph} of a ribbon curve structure is the graph
obtained by replacing each box with a vertex and connecting each
pair of glued boxes by an edge (one for each pair of segments
identified). The junction graphs of isomorphic ribbon
curve structures are isomorphic. Clearly,
the valency of a vertex of a junction graph
does not exceed the number of zigzag segments of the box
represented by the vertex.
In particular, the valency is at most~$3$,
and any vertex of valency~$3$ represents
a box of type~$\II_3$.

A ribbon curve structure defines an equivalence class
of
totally
reduced dessins: one replaces each box
with
a corresponding cubic,
performs a junction on each pair of identified zigzag
segments,
and performs a self-junction on each vanishing segment.
A \emph{ribbon curve} is a trigonal curve
whose dessin is equivalent to one
obtained in this way.
Note that,
since cubics are indecomposable,
each particular dessin
admits at most one ribbon curve structure.

An \emph{enhanced ribbon curve structure} is a ribbon curve
structure equipped with a collection of nonnegative integers, one
for each box other than $\II_3$, one for each adjacency,
and one for each vanishing segment.
The notion of isomorphism extends naturally:
one requires that the isomorphism should
preserve the enhancement. An enhanced ribbon curve structure
defines an equivalence class of dessins: one takes the
totally
reduced
dessin constructed in the previous paragraph and
inflates it by placing
the
indicated number of
\dotted- monochrome
components to each
(self-)junction and to an inner \dotted- edge within each box.
Considering the types of boxes one by one, see Lemma~\ref{regions-cubic},
one can easily show that the equivalence class is indeed well
defined:
for each cubic dessin, any two distinct inner
\dotted- edges in it can be connected by a sequence of
elementary moves; hence, the \dotted- monochrome components introduced
inside the box can be placed to any preselected inner \dotted-
edge.

\subsection{Classification of trigonal $M$-curves}\label{class-M}

\theorem\label{M-curves}
The collection of all \dotted- cuts of any totally reduced
non-hyperbolic $M$-dessin~$\Gamma$
represents~$\Gamma$ as
an iterated \rom(self-\rom)junction
of a union of $M$-cubics.
Furthermore,
any elementary
move of~$\Gamma$ is either a simple modification of the
junction or an elementary move in one of the cubics
\rom(not involving the cuts\rom).
\endtheorem

\proof
We will show that a totally reduced
non-hyperbolic $M$-dessin $\Gamma \subset \base$
that is not a cubic is a (self-)junction of another
$M$-dessin.
One has $d=0$; hence, $\delta=0$ and
$\Nzigzag+\Ninner\le 2$ if $\base$ is a disk,
and $\Nzigzag+\Ninner \le 1$ otherwise,
see~\ref{tg->ovals}.
Thus,
according to Theorem~\ref{6k<=}, the dessin~$\Gamma$ has a dotted cut.
From the oval
count given by Proposition~\ref{oval.count} it follows that any such cut
is a (self-)junction, the result being an $M$-dessin.
(Roughly, to keep the maximal number of components, each
gluing must create at least two ovals, and each self-gluing
must create at least three components.
The possibility to form an even hyperbolic component
is ruled out by the assumption that the dessin is totally
reduced.)
In particular, the result of the cut has no hyperbolic components
and, hence, is still totally reduced.

The second statement follows from the first one and the fact that
a cubic is indecomposable. Indeed, the only elementary move not as
in the theorem is an inner modification joining two cuts (within one
cubic) and
producing an alternative pair of cuts. However, from the point of
view of the cubic, that would
imply the existence of
an artificial \dotted- cut, which
would contradict to
the fact that the cubic is indecomposable.
\endproof

\corollary\label{M-structure}
Any totally reduced non-hyperbolic $M$-dessin
admits a ribbon curve structure,
which has the following properties\rom:
\roster
\item"--"
each box represents an $M$-cubic\rom;
\item"--"
the underlying surface is orientable.
\endroster
Conversely, any ribbon curve structure with the properties above
defines a totally reduced $M$-dessin.
Furthermore, two such dessins
are equivalent if and only if
their ribbon curve structures are isomorphic.
\endcorollary

\remark{Remark}
As each box representing an $M$-cubic has valency at most~$2$, the
junction graph of an $M$-curve is either a linear tree or a single
cycle, see Figure~\ref{figure-graphs1}.
In the former case, the underlying surface is a disk, in
the latter case it is an annulus.
\endremark

\midinsert
\eps{Fig-M}
\figure\label{figure-graphs1}
Junction graphs of $M$-dessins
\endfigure
\endinsert

\remark{Remark}
The statement of Corollary~\ref{M-structure} for the case of rational base
was first obtained by
S.~Orev\-kov~\cite{Orevkov-private}.
\endremark

\proof[Proof of Corollary~\ref{M-structure}]
The only statement that needs proof is the fact that, if the
underlying graph is a single cycle, the underlying surface cannot
be a M\"{o}bius band. This possibility is eliminated
by Proposition~\ref{oval.count}.
\endproof

\theorem\label{M-all}
The deformation classes
of
almost generic non-hyperbolic
trigonal $M$-curves are in a canonical
one-to-one correspondence with the isomorphism classes
of enhanced ribbon curve structures as in
Corollary~\ref{M-structure}.
\endtheorem

\proof
Due to Proposition~\ref{equiv.curves},
it suffices to enumerate the equivalence
classes of non-hyperbolic $M$-dessins.
Any such dessin~$\Gamma$ has ovals, and hence real \cross-vertices.
Thus, $\Gamma$ is equivalent to an inflation
of a totally reduced dessin, and, in view
of Lemma~\ref{partial}, the statement
of the theorem follows from Corollary~\ref{M-structure}.
\endproof

\theorem\label{M-hyp}
For each integer $g \ge 0$, there is a unique deformation class
of almost generic hyperbolic trigonal $M$-curves over a base
of genus~$g$.
\endtheorem

\proof
Proposition~\ref{equiv.curves}
reduces the problem to dessins.
Consider a hyperbolic
$M$-dessin. The oval count~\eqref{oval.count-hyp}
implies that $\delta = 0$, $\Lodd = 1$, and $\deg\Gamma = 3$.
Hence, $\Gamma$ is equivalent to an inflation
of a totally reduced dessin, which is a
hyperbolic cubic; the latter is equivalent to
the dessin shown in Figure~\ref{fig-star-like},
see Lemma~\ref{bw-hyp}.
As above, the equivalence class of the inflation is determined
by the number~$g$ of the components inserted.
\endproof

\subsection{Classification of elliptic $M$-surfaces}\label{M-surfaces}
Define a \emph{ribbon surface structure} as a ribbon curve
structure satisfying the following additional requirements:
\roster
\item"--"
there
are
no vanishing segments;
\item"--"
the combined parity of the segments within
each boundary component of the underlying
surface
is
even;
\endroster
and\mnote{It: they modified the rosters, and now this line
looks very strange}
enriched with the following decorations:
\roster
\item"--"
each junction (\ie, common segment of two boxes glued together)
is subdivided into segments, and each segment and
each inner vertex of the subdivision are given a sign;
\item"--"
each box other than $\II_3$ is given a sequence of signs
of even length.
\endroster
The decorations must be subject to the following condition:
the product of the signs
given to all segments in the junctions
adjacent to one box is $-1$.

From the definition it follows that the total number of boxes in
a ribbon surface structure must be even. In particular, there
are
no boxes of types~$\I_0$ or~$\II_0$.

An \emph{isomorphism} of ribbon surface structures is an
isomorphism of the corresponding ribbon curve structures
preserving the additional decorations.

\theorem\label{th.M-surfaces}
Each
almost generic elliptic $M$-surface defines an isomorphism
class of ribbon surface structures
with the following properties\rom:
\roster
\item"--"
each box represents an $M$-cubic\rom;
\item"--"
the underlying surface is orientable
\endroster
\rom(\cf. Corollary~\ref{M-structure}\rom).
Conversely, each
ribbon surface structure
with the properties above
defines one or two \rom(depending on
whether the underlying surface is a disk or an annulus,
respectively\rom)
deformation
classes of pairs of opposite elliptic $M$-surfaces.
\endtheorem

\proof
Essentially, the statement follows from
Theorem~\ref{M-all}; we will just explain the relation between
elliptic surfaces, trigonal curves,
ribbon surface structures, and enhanced ribbon
curve structures. In view of Corollary~\ref{M->Jacobian},
up to deformation each
$M$-surface is Jacobian; hence, it is described by its
Weierstra{\ss} model, see~\ref{wmodels},
and the branch curve is $M$-, see Lemma~\ref{M-d}.
Thus, an $M$-surface determines and is determined by the
following data: a trigonal $M$-curve~$C$ over the same base~$B$
(hence, an enhanced ribbon
curve structure), a lift of the monodromy
$\pi_1(B^\#)\to\Braids_3/\Delta^2$ to
$\SL(2,\Z)$, see~\ref{monodromy}, and a
choice of one of the two opposite real structures (which is the
reason why the theorem is stated about pairs of opposite
surfaces).

If the underlying surface of the ribbon curve structure is a disk,
then the group $\pi_1(B^\#)$ is generated by small
loops~$\alpha_i$ around the singular fibers, the classes~$\beta_j$
of the hyperbolic components of~$B_\R$,
and the doubles~$\gamma_k$
of the \dotted-segments connecting the components of~$B_\R$. The
monodromy along each loop~$\alpha_i$ and its lift to $\SL(2,\Z)$
are determined by the requirement that the singular fibers of the
surface must be of type~$\I_1$. The $\Braids_3/\Delta^2$-valued
monodromy along each loop~$\beta_j$ or~$\gamma_k$ is trivial,
see~\ref{tg->monodromy}, and, hence, its lift to $\SL(2,\Z)$ is
$\pm\id$. The sign in front of~$\id$ is the sign assigned to the
corresponding \dotted-component/segment, and the collection of
signs thus obtained (as well as the number of hyperbolic
components)
is encoded by the additional decoration in the
definition of ribbon surface structure: the vertices subdividing
the junctions into segments represent even \dotted-components,
and a chain of $2m$ signs assigned to a box
represents $m$ even \dotted-components on the \dotted-segment
inside the box attached to the oval. (Lemma~\ref{regions-cubic}
asserts that such a segment always exists.) In the latter case,
the even numbered signs in the chain are those assigned to the
\dotted-components (counted starting from the oval), and the odd
numbered ones are the signs assigned to the segments connecting
two consecutive components (or the first component and the oval).

The relation between the signs required in the definition of the
ribbon surface structure is a manifestation of a relation between
equivariant characteristic classes of a line bundle. In simple
terms, it can be
obtained as follows. As explained
in~\ref{realhominv} and~\ref{monodromy},
the homological invariants form an affine
space over $\HH^1(B;\Z_2)$, whence the signs involved are subject
to one relation for each ribbon box. Pick a box, cut it out of the
underlying surface, and cut out a disk containing the inner
monochrome components. The result can be regarded as a cubic with
a few half disks at the boundary removed.
By an affine shift the signs can be
chosen so that the monodromy along each cut but one is trivial.
If it were also trivial on the remaining cut, it would extend
to the whole cubic, thus producing an elliptic surface over a
trigonal curve of odd degree, which is a contradiction.

If the underlying surface is an annulus, there is an additional
pair of complex conjugate cycles, on which a lift should be chosen.
This
accounts for the fact that, in this case, the ribbon surface
structure defines two pairs of deformation classes.
\endproof

\remark{Remark}
If the underlying surface is an annulus,
the two lifts of the monodromy along an additional cycle
are also topologically distinct, see
Lemma~\ref{Aconj-A}.
Instead of choosing a lift, which depends on a particular choice
of the cycle, one can distinguish the two surfaces by the
orientability of the principal component over one of the two
boundary components of the annulus,
see~\ref{real-part}.
\endremark

\corollary\label{M.g=0}
The deformation classes of pairs of opposite
almost generic elliptic $M$-surfaces
over a rational base are in a canonical
one-to-one correspondence with the isomorphism classes of ribbon
curve structures with the following properties\rom:
\roster
\item"--"
each box
represents an $M$-cubic,
and the number
of boxes is even\rom;
\item"--"
there
are
no vanishing segments\rom;
\item"--"
the junction graph is a linear tree.\mnote{It: they removed
systematically \qed after corollaries without proofs.
I would prefer to put \qed back}
\qed
\endroster
\endcorollary

\subsection{$(M-1)$-curves and surfaces}\label{S.M-1}

\theorem\label{M-1-curves}
The collection of all \dotted- cuts of any
totally reduced non-hyperbolic $(M-1)$-dessin~$\Gamma$
represents~$\Gamma$ as
an iterated \rom(self-\rom)junction\mnote{It: `self' should be in italic}
of a union
of $M$-cubics and at most one
$(M-1)$-block, the latter being either a sextic
\rom(\ie, a dessin of degree six
on a disk\rom) or a block of degree~$3$.
Furthermore, any elementary
move in~$\Gamma$ is either a simple modification of the
junction, or an elementary move in one of the blocks
\rom(not involving the cuts separating distinct blocks\rom),
or an elementary move in a sextic that is the junction
of the $(M-1)$-cubic and an $M$-cubic.
\endtheorem

\proof
The
proof is almost literally the same as in the case of
$M$-curves.
Let~$\Gamma\subset\base$ be a totally reduced non-hyperbolic
$(M-1)$-dessin
that is neither of degree~$3$ nor a sextic.
One has $d = 1$; hence, either
\roster
\item"" $\delta=1$, $\Nzigzag+\Ninner = 0$, and $\base$
is a M\"{o}bius band, or
\item"" $\delta=0$ and
$\Nzigzag+\Ninner\le 3$ if $\base$ is a disk,
and $\Nzigzag+\Ninner \le 2$ otherwise,
\endroster
see~\ref{tg->ovals}.
Thus,
according to Theorem~\ref{6k<=}, the dessin~$\Gamma$ has a \dotted- cut.
The oval
count~\ref{oval.count} implies that a \dotted- cut
in~$\Gamma$ must be a (self-)junction,
the result being
an $M$- or $(M-1)$-dessin.

For the second statement,
an elementary move not as in the theorem would destroy a junction
and create a new one (due to the first statement). Hence, at the
very moment of the modification the graph would have a genuine
\dotted- cut that is not a junction.
\endproof

\theorem\label{M-1-structure}
Any totally reduced non-hyperbolic $(M-1)$-dessin
is equivalent to a dessin
that admits a ribbon curve structure such that either
\roster
\item"--"
each box
represents an $M$-cubic and has valency~$2$ \rom(so that
the junction graph is a single cycle\rom),
and the underlying surface is a M\"{o}bius band, or
\item"--"
exactly one
box
represents an $(M-1)$-cubic, while all other
boxes
represent $M$-cubics, and the underlying surface is orientable.
\endroster
Conversely, any ribbon curve structure as above defines a totally reduced
non-hyperbolic \hbox{$(M-1)$}-\penalty0dessin. Furthermore, two such dessins are
equivalent if and only if their ribbon curve structures are
isomorphic or connected by
one or several of the following moves:
$\II_2+\I_2\leftrightarrow\I_2+\II_2$
provided that the adjacency is positive,
$\II_2+\I_1\leftrightarrow\I_2+\II_1$,
$\II_1+\I_2\leftrightarrow\I_1+\II_2$,
or $\II_1+\I_1\leftrightarrow\I_1+\II_1$.
\endtheorem

\proof
First, show the decomposability. In view of
Theorem \ref{M-1-curves}, it suffices to consider a
dessin~$\Gamma$ that either is a sextic on a disk or has
degree~$3$. Theorem~\ref{reduction} implies that $\Gamma$~is
equivalent to a dessin with a cut. It remains to consider all
gluings of one or two cubics, select those that are
$(M-1)$-curves, and, in each case, see that the dessin is
equivalent to a junction, making sure that the equivalence leaves
intact the zigzags.
This is straightforward. In fact, from
Proposition~\ref{oval.count} it follows that, since $\Gamma$~is an
$(M-1)$-dessin, the cut is either a \cross-cut (the result of the
cut being an $M$-curve) or a \dotted- cut. A genuine \dotted- cut
is necessarily a junction. An artificial \dotted- cut cannot join
two real \dotted- segments adjacent to \cross-vertices (as it must
destroy two ovals); hence, creating a bridge, one can replace the
cut with a genuine one.

\midinsert
\labeldelta=.125cm
\def\Y{\the\count0}
\def\X{-\XX}
\def\XX{18}
\eps{Fig-decomp}
\centerline{%
\count2=21
\count0=\count2
\multiply\count0 by4
\advance\count0 -1
\advance\count0-\count2
\labelat[\X,\Y]{$\I_2+\II_3$}%
\labelat[\XX,\Y]{$\I_1+\II_3$}%
\advance\count0-\count2
\labelat[\X,\Y]{$\I_2+\II_2\to\II_2+\I_2$}%
\labelat[\XX,\Y]{$\I_1+\II_2\to\II_1+\I_2$}%
\advance\count0-\count2
\labelat[\X,\Y]{$\I_2+\II_2$}%
\labelat[\XX,\Y]{$\I_1+\II_1\to\II_1+\I_1$}%
\advance\count0-\count2
\labelat[\X,\Y]{$\I_2+\II_1\to\II_2+\I_1$}%
}
\figure\label{figure}
(Re-)decompositions of $(M-1)$-sextics
\endfigure

\bigskip\bigskip
\eps{Fig-M-1}
\figure\label{figure-graphs2}
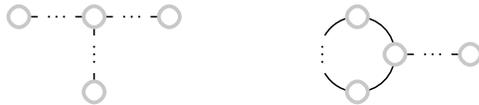
Extra junction graphs of $(M-1)$-dessins
\endfigure
\endinsert

The forms of the ribbon curve structures are easily enumerated
using
Proposition~\ref{oval.count}. According to
Theorem~\ref{M-1-curves}, in order to study the equivalences, it
suffices to consider an \hbox{$(M-1)$}-sextic
decomposed into a junction of two cubics,
and study its re-decompositions.
The sextics and their re-decompositions
are shown schematically in Figure~\ref{figure}.
(Recall that a box
of type~$\I_2$ and a box of type~$\II_2$ can be joined in two
different ways, forming a positive or negative
junction.)
Each re-decomposition shown (the double lines in the figure)
can easily be realized by a sequence
of elementary moves.
A way to prohibit the other re-decompositions is to consider the
distribution of the maximal real \dotted-/complementary segments
and their parities.
\endproof


\remark{Remark}
The junction graph of an $(M - 1)$-dessin
is either a linear tree,
or a single cycle, or one of the two graphs shown in
Figure~\ref{figure-graphs2}, the vertex of valency~$3$
representing the $(M-1)$-cubic, which is
of type~$\II_3$.
In the second graph in
Figure~\ref{figure-graphs2}, the cycle may as well consist of one
or two vertices, \cf. Figure~\ref{figure-graphs1}.
\endremark

\corollary\label{M-1-all}
The deformation classes
of almost generic trigonal $(M-1)$-curves are in a canonical
one-to-one correspondence with the isomorphism classes
of enhanced ribbon curve structures as in
Theorem~\ref{M-1-structure} modulo the following
additional equivalence relation\rom:
\roster
\item"--"
the moves as in Theorem~\ref{M-1-structure}\rom;
the three integers assigned to the two boxes
and the adjacency involved can be chosen arbitrarily
provided that the sum of the integers is left intact\rom;
\item"--"
the integers within a box of type~$\II_1$ or~$\II_2$
\rom(\ie, those assigned to the box itself,
its adjacencies, and its vanishing segments\rom) can be changed
arbitrarily provided that their sum is left intact.
\endroster
\endcorollary

\proof
First, notice that an $(M-1)$-curve cannot be hyperbolic.
Indeed, with $d=1$
the oval count~\eqref{oval.count-hyp} implies
$\delta=0$,
which
contradicts Corollary~\ref{hyp.cong}. Thus,
the curve is non-hyperbolic, and, similar to Theorem~\ref{M-all}, the
problem can be reduced to Theorem~\ref{M-1-structure}.

The realizability of both moves can be deduced
from Lemma~\ref{regions-cubic}. As in the proof of
Theorem~\ref{M-1-structure}, to show that there
are
no
others, it suffices to consider sextic dessins. This can be done
on a case by case basis, using Figure~\ref{figure} and a
careful analysis of
real \dotted- monochrome vertices.
\endproof

\paragraph\label{M-1.surface}
In view of Corollary~\ref{M->Jacobian}, the classification of
almost generic elliptic $(M-1)$-surfaces also reduces to the classification of
almost generic trigonal $(M-1)$-curves enhanced with a lift of the monodromy
$\pi_1(B^\#)\to\Braids_3/\Delta^2$ to $\SL(2,\Z)$. As in the case
of $M$-surfaces, this procedure could be expressed in terms of
ribbon surface structures. An additional complication is the fact
that $(M-1)$-curves enjoy much more freedom (the equivalence
relations described in Corollary~\ref{M-1-all}) resulting in a vast
number of moves for the monodromy. For this reason, we confine
ourselves to the case of genus zero, where the monodromy is
uniquely determined by the dessin and, hence, the result is an
immediate consequence of Corollary~\ref{M-1-all}.

\proposition\label{M-1.g=0}
The deformation classes of pairs of opposite almost generic elliptic
$(M-1)$-surfaces
over a rational base
are in a canonical one-to-one correspondence with the
isomorphism classes of ribbon curve structures with the following
properties\rom:
\roster
\item"--"
the number of boxes is even\rom;
\item"--"
one box represents an $(M-1)$-cubic,
the others representing $M$-cubics\rom;
\item"--"
there
are
no vanishing segments\rom;
\item"--"
the junction graph is a tree.
\qed
\endroster
\endproposition

\subsection{Oval chains}
In this section we derive a few simple consequences of the
classification results obtained above.

\theorem\label{odd.chains}
Let $C$ be a nonsingular trigonal $M$-curve of degree~$3k$ on a real ruled
surface over a base~$B$. Then
the following holds\rom:
\roster
\item
each non-complete maximal chain of ovals of~$C$ is of odd
length\rom;
\item
if $C$ has no complete chains, then it has
$k-2+\Lodd$ \rom(respectively,~$k$\rom) maximal chains if $B_\R$ has one
\rom(respectively, two\rom) non-hyperbolic components.
\endroster
\endtheorem

\proof
Corollary~\ref{M-structure} lists all dessins of
non-hyperbolic
$M$-curves, and
the maximal chains of ovals are easily seen: the ovals are
described in~\ref{tg->real}, and chain breaks are
the
maximal sequences consisting of an odd number of
odd segments (and any number of even segments other than ovals).
\endproof

\corollary\label{odd.chains.rational}
The ovals of a nonsingular trigonal $M$-curve on a real rational ruled
surface $\Sigma_k$, $k\ge3$, form $k-2$ maximal chains,
each maximal chain being of odd length.\mnote{It: \qed back}
\qed
\endcorollary

Let $C'\subset\Sigma'\to B'$ and $C''\subset\Sigma''\to B''$ be two real
trigonal curves on real ruled surfaces.
The curves~$C'$ and~$C''$ are said to have the
same \emph{fibered real scheme} if there is a fiberwise
homeomorphism $\Gf\:\Sigma''_\R\to\Sigma'_\R$
such that $C'_\R$ and the
image
$\Gf(C''_\R)$ can be connected by an isotopy $C^t\subset\Sigma'_\R$
during which the intersection of~$C^t$ with any fiber of the
projection $\Sigma'_\R\to B'_\R$ consists of at most three points.
In other words, to make the result slightly more general,
we allow passing through vertical flexes, \ie, straightening
zigzags,
\cf.~\ref{straightening}.

\theorem
There
is a nonsingular trigonal $(M-1)$-curve
on a rational ruled surface
such that the
fibered real scheme of the curve cannot
be obtained by a single Morse modification from the fibered real
scheme of a nonsingular $M$-curve.
\endtheorem

\midinsert
\eps{S6M-1}
\figure\label{fig.M-1}
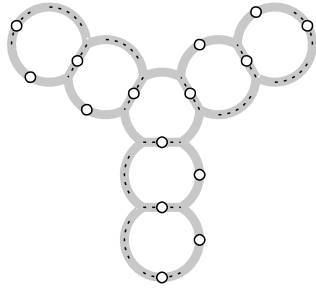
$(M-1)$-curve with three maximal chains of length six
\endfigure
\endinsert

\proof
In fact, any $(M-1)$-curve whose
junction graph has a vertex of valency three and three branches
of length at least two each has the desired property. An
example is shown in Figure~\ref{fig.M-1}. The reason is that
the curve in question has three maximal chains of ovals of even
length; hence, in view of Corollary~\ref{odd.chains.rational}, its
fibered real scheme
cannot be obtained by erasing one oval from the fibered real
scheme of an $M$-curve.
\endproof

\subsection{Further generalizations and open questions}

\subsubsection{Singular curves}\label{singular-curves}
According to the definition given in~\ref{trigonal-curves}, a
trigonal curve in a ruled surface $\Sigma$ is not supposed to
intersect the distinguished section~$s$ of~$\Sigma$. Relax this
requirement and consider a curve $C\subset\Sigma$ that does
intersect~$s$ at a point~$P$. The elementary transformation
of~$\Sigma$ at~$P$ (\ie, blowing up~$P$ and blowing down the fiber
through~$P$) produces a new surface~$\Sigma'$, section~$s'$, and
curve $C'\subset\Sigma'$ which intersects~$s'$ with a
smaller multiplicity; the image of the fiber blown down is a
singular point of~$C'$ (a node or a cusp if $C$ was nonsingular).
Iterating this procedure, one arrives at a
surface~$\Sigma''$, section~$s''$, and curve $C''\subset\Sigma''$
disjoint from~$s''$, \ie, a trigonal curve in the sense
of~\ref{trigonal-curves}. The curve is singular: it has one
type~$\bA_{2m-1}$ or~$\bA_{2m}$ singular point for each
$m$-fold intersection point
of the original curve~$C$ and section~$s$. The inverse elementary
transformations (blowing up the singular points and contracting
the corresponding fibers) convert~$C''$ back to~$C$. Thus, the
deformation classification of trigonal (in the wide sense) curves
intersecting~$s$ at several points with prescribed multiplicities
can be reduced
to that of trigonal (in the sense
of~\ref{trigonal-curves}) curves with several type~$\bA$ singular points.
(Note that the degenerations $\bA_{2m-1}\to\bA_{2m}$ should be
allowed during the deformations; these degenerations correspond to
the confluences of vertical tangents and intersections with the
exceptional section.)
If the multiplicities are not prescribed, one should consider
curves with a certain number of nodes
and allow
deeper
confluence of the nodes during the
deformations.

The classification of non-hyperbolic
singular curves (with type~$\bA$ singularities only)
that perturb to nonsingular
$M$- or $(M-1)$-curves is essentially contained in Theorem~\ref{M-all}
and Corollary~\ref{M-1-all}.
Indeed, type~$\bA$ singularities are obtained
by bringing together some of the
vertical tangents of a nonsingular curve.
Hence, the graph (of the $j$-invariant) of the singular curve
is obtained from that of a nonsingular one
by bringing together some of the
\cross-vertices. Obviously,
several
consecutive
\cross-vertices can be brought together if and only if, after a
sequence of \white-ins and \black-ins, they are not separated by
\white-- or \black-vertices. This observation gives a clear
description of
the singular curves in question,
by either
referring to the ribbon curve structures of nonsingular curves and
indicating the sequences of \cross-vertices to be brought
together, or else constructing singular curves directly from boxes of
more general form (including the graphs of singular cubic)
\via\ a more general
junction operation (allowing forming singular points instead of
ovals).

\subsubsection{Straightening zigzags}\label{straightening}
On the account of the principal tool used in this paper
(constructing deformations of curves \via\ deformations of $j$-invariants),
we state our results in the language of equivariant fiberwise deformations
without confluence of singular fibers.
However, from the point of view of geometry of nonsingular curves,
vertical flex should not be considered
a singularity.
Passing through a vertical flex during a deformation
results in the removing (straightening) or creating a zigzag.
In spite of its apparent simplicity, this operation
does not lead to a deformation of the $j$-invariant:
at the very moment of the modification the degree of~$j$
drops by~$2$. The corresponding modification
of dessins is shown in Figure~\ref{s-zigzag}
(see also~\cite{Zvonilov}):
two adjacent triangles are removed, and two new triangles
are inserted. Note that forming a vertical flex is not
as local as forming a type $\bA$ singular point
(\cf.~\ref{singular-curves}); the possibility to bring
together a pair of \cross-vertices bounding a zigzag
cannot be deduced solely from the real part.
An example of a nonsingular trigonal curve with a zigzag
that cannot be straightened in a single step
was found by Orevkov~\cite{Orevkov-private}.

\midinsert
\eps{Fig-zigzag}
\figure\label{s-zigzag}
Straightening a zigzag
\endfigure
\endinsert

Due to Theorems~\ref{M-curves} and~\ref{M-1-curves}, whenever an $M$-
or an $(M-1)$-curve is decomposed in ribbon boxes,
each zigzag is localized within
a single cubic, and using Lemma~\ref{regions-cubic} one can conclude
that each zigzag can be straightened
without destroying the junctions.
Hence, the classification of nonsingular $M$- and $(M-1)$-curves
up to the new relaxed equivalence relation
can be deduced from Theorems~\ref{M-curves} and~\ref{M-1-curves}.
We refrain from attempting to formulate a precise statement for
$(M-1)$-curves. Just note that straightening a zigzag
increases the number of moves as in Theorem~\ref{M-1-curves}
that can be applied to the dessin.
Effectively, this increased flexibility means that, after its
disappearance, a zigzag can freely slide along a chain of ovals and
reappear at a new place.
In the $M$-case, the ribbon curve structure is always rigid.
Thus, the only modification is
the disappearance of
a zigzag, possibly followed by its reappearance next to the same
oval,
pointing to it from the other side.
At the level
of the ribbon curve structure,
this modification is either a change of type ($\I_1$ to $\I_2$
or {\it vice versa}) of the corresponding cubic (which is
necessarily located at one of the two ends of the junction graph)
or a change of the sign of its junction.

\subsubsection{Ribbon vs. unstructured curves}
The ribbon curve construction produces an interesting class of
trigonal curves with a clearly defined structure. Under various
mild assumptions (\eg, if all \cross-vertices are real, \ie,
assuming that all boxes are of type~$\I_2$ or~$\II_3$) the
structure is rigid: the junctions are present in any dessin
equivalent to a given one, and they cannot be destroyed or
modified by elementary transformations. Thus, the correspondence
between curves and ribbon curve structures gives a deformation
classification of such curves. On the other hand, there obviously
are large `unstructured' curves whose graphs do not contain a
junction or even a cut. At present, it is unclear whether and how
the property of being a ribbon curve can be characterized in
topological terms, or what a general classification theorem would
look like. Probably, the dessin of a general curve would be a
union of a ribbon part and a few unstructured pieces.

The simplest example of an `unstructured' curve is found among
$(M-2)$-sextic.
Indeed, one can glue two
cubics of type~$\I_2$ each along a pair of \solid- segments so
that the resulting sextic has two ovals not separated by zigzags
(and hence four zigzags not separated by ovals).
On the other hand, any sextic that is a junction
of two cubics may have at most two zigzags within each real segment
connecting the two ovals produced by the junction.

This example
shows that, starting from the $(M-2)$-case,
the structure theorems should have a form different from
the classification results of the paper.
(That is why we confine ourselves to $M$- and
$(M-1)$-curves only.)
Both of the key ingredients used in our approach fail:
first, Theorem~\ref{6k<=} does not break
a dessin into sufficiently small pieces;
second, it is no longer true that any dotted cut
is a junction (note that the sextic above is equivalent
to a dessin with a dotted cut).

\subsubsection{Quasi-simplicity: still open}
Let us
briefly discuss the relation between
the ribbon box decomposition of a trigonal $M$-curve (without
hyperbolic components) and its real part.
Certainly, the former does determine the latter,
and there is a number of situations in which the converse
also holds, \ie, the ribbon curve structure is recovered
from the sequence formed by the zigzags and maximal chains of ovals.
Among these situations are:
\roster
\item\local1 $M$-curves over a base~$B$ of genus one;
\item\local2 $M$-curves with at least one zigzag;
\item\local3 $M$-curves with a sufficiently generic (in the sense
described below) distribution of ovals.
\endroster
In case~\loccit1 the sequence of maximal chains of ovals in one
of the two components
of~$B_\R$
clearly determines the box decomposition.
In the case of rational base
the junction graph is a linear tree
and
there are two distinguished ovals which are located
in the cubics corresponding to
its
two extreme vertices;
we call them \emph{extreme} ovals.
If at least one extreme oval is known, the rest of the ribbon curve
structure is found uniquely. This observation covers
case~\loccit2, as zigzags of an $M$-curve always point
at its extreme ovals.

In fact, the ribbon curve structure can still be recovered,
at worst up to the `horizontal' symmetry,
starting from one extreme chain,
\ie, chain containing an extreme oval.
If the base is rational
and there
are
no zigzags, each maximal chain
of length $2k + 1$ containing $l = 0$, $1$, or $2$ extreme ovals
is opposed by a maximal sequence of $k - l$ solitary ovals
(\ie, those forming maximal chains of length~$1$).
Hence, the lengths of the extreme chains
are recovered from the sequence of maximal chains.
By case~\loccit3 we mean the situation when
these lengths determine the extreme chains.

However, there are sequences of maximal chains that can be obtained
from two non-isomorphic ribbon curve structures.
The simplest example that we know is the following sequence
of $24$ chains:
$$
\underline5\ 1\ 3\ 3\ 1\ 3\ 3\ 5\
\underline{\underline5}\ 3\ 1\ 3\
\underline5\ 3\ 5\ 3\ 1\ 3\ 1\ 3\ 3\
\underline{\underline5}\ 3\ 5 \ .
$$
(Of course, only the lengths of the chains are listed.)
The chains containing extreme ovals
are either those underlined, or those double underlined.
As the sequence have no symmetries, the two ribbon curve structures
are not isomorphic.
\end


\ifx\defloaded\undefined\input elliptic.def \fi

\begin

\end

\endgroup


\Refs \widestnumber\key{DIK1}

\ref{AMn}
\by M.~Akriche, F.~Mangolte
\paper Nombres de Betti des surfaces elliptiques r\'{e}elles
\jour Preprint ArXiv: math.AG/0604131, 2006
\endref\label{Mn-Akriche}

\ref{Ba}
\by E.~Ballico
\paper
Real analytic elliptic surfaces
\jour J. Pure Appl. Algebra
\vol 178
\yr 2003
\issue 3
\pages 225--234
\endref\label{Ballico}

\ref{BPV}
\by W.~Barth, C.~Peters, A.~van~de~Ven
\book Compact complex surfaces
\bookinfo Ergebnisse der Mathematik und ihrer
Grenzgebiete (3)
\publ Springer-Verlag
\publaddr Berlin-New York
\yr 1984
\endref\label{BPV}

\ref{BE}
\by I.~Bernstein, A.~L.~Edmonds
\paper On the classification of generic branched coverings of surfaces
\jour Illinois J. Math
\vol 28
\issue 1
\yr 1984
\endref\label{BE}

\ref{BMn}
\by F.~Bihan, F.~Mangolte
\paper Topological types of real regular jacobian elliptic surfaces
\jour
Geom. Dedicata
\vol 127
\yr 2007
\pages 57-73
\endref\label{Mn-Bihan}


\ref{Br1}
\by G.~E.~Bredon
\book Introduction to compact transformation groups
\publ Academic Press
\publaddr New York, London
\yr 1972
\endref\label{Bre1}

\ref{Br2}
\by G.~E.~Bredon
\book Sheaf Theory, Second Edition
\bookinfo Graduate Texts in Mathematics
\yr 1997
\publ Springer-Verlag
\publaddr New York, Berlin, Heidelberg
\endref\label{Bre2}

\ref{Bro}
\by K.S.~Brown
\book Cohomology of groups
\publ
Springer-Verlag
\publaddr New York-Heidelberg-Berlin
\yr 1982
\endref\label{Bro}

\ref{CF}
\by F.~Catanese, P.~Frediani
\paper Real hyperelliptic surfaces and the orbifold fundamental group
\jour J. Inst. Math. Jussieu
\vol 2
\issue 2
\pages  163-233
\yr 2003
\endref\label{CF}


\ref{Co1}
\by A.~Comessatti
\paper Fondamenti per la geometria sopra le superficie
razionali dal punto di vista reale
\jour  Math. Ann.
\vol 43
\yr 1912
\pages 1--72
\endref\label{Co-rat}

\ref{Co2}
\by A.~Comessatti
\paper Sulle variet\`{a} abeliane reali, \rom{I} e \rom{II}
\jour  Ann. Mat. Pura Appl.
\vol   2 and 4
\yr    1924 {\rom and} 1926
\pages 67--106 and 27--71
\endref\label{Co-abel}

\ref{DIK1}
\book Real Enriques surfaces
\by A.~Degtyarev,
I.~Itenberg, V.~Kharlamov
\bookinfo Lecture Notes in Math.
\vol 1746
\yr 2000
\publ Springer-Verlag
\endref\label{DIK1}

\ref{DIK2} \paper Finiteness and Quasi-Simplicity for Symmetric
K3-Surfaces \by A. Degtyarev, I. Itenberg, V. Kharlamov \jour Duke
Math.~J. \yr 2004 \vol 122 \pages 1--49
\endref\label{duke}

\ref{DK1}
\by A.~Degtyarev, V.~Kharlamov
\paper Real Enriques surfaces without real points and Enriques-Einstein-Hitchin 4-manifolds
\inbook
The Arnoldfest (Toronto, ON, 1997)
\pages  131--140
\bookinfo Fields Inst. Commun., 24
\publ Amer. Math. Soc.
\publaddr Providence, RI
\year 1999
\endref\label{DK-empty}

\ref{DK2}
\by A.~Degtyarev, V.~Kharlamov
\paper Real rational surfaces are quasi-simple
\jour Jour. reine angew. Matem
\vol 551
\yr 2002
\pages 87--99
\endref\label{DK-rat}

\ref{Fr}
\by P. Frediani
\paper Real Kodaira surfaces
\jour Collect. Math.
\vol 55
\yr 2004
\issue 1
\pages 61--96
\endref\label{Frediani}

\ref{FM}
\by R.~Friedman, J.~W.~Morgan
\book Smooth four-manifolds and complex surfaces
\bookinfo Ergebnisse der Mathematik und ihrer
Grenzgebiete (3)
\publ Springer-Verlag
\publaddr Berlin-New York
\yr 1994
\endref\label{FM}

\ref{GW}
\by M.~Gross, P. M. H. Wilson
\paper Mirror symmetry via 3-tori
for a class of Calabi-Yau threefolds
\jour Math. Ann.
\vol 309
\issue 3
\pages 505--531
\yr 1997
\endref\label{Gross-Wilson}

\ref{Ka}
\by A.~Kas
\paper On the deformation types of regular elliptic surfaces
\inbook Complex analysis and algebraic geometry,
pp. 107--111. Iwanami Shoten, Tokyo, 1977
\endref\label{Kas}

\ref{Kh}
\by V.~Kharlamov
\paper Towards the Viro, Ragsdale, and Petrovskii conjectures
\jour Proc. Internat. Topological Conference, Baku, 1987
\pages 412
\endref\label{Kharlamov}

\ref{KK}
\paper Surfaces with DIF $\ne$ DEF real structures
\by V.~Kharlamov, Vik.~Kulikov
\jour
Izv. Math.
\vol 70
\yr 2006
\issue 4
\pages 769-807
\endref\label{KK}

\ref{Kl} \by F.~Klein \paper Ueber Fl\"{a}chen dritter Ordnung
\jour Math. Ann. \vol   6 \yr    1873 \pages 551--581
\endref\label{Klein}

\ref{Ko}
\by K.~Kodaira
\paper On compact analytic surfaces, II--III
\jour Annals of Math.
vol 77--78
\yr 1963
\pages 563--626, 1--40
\endref\label{Kod}

\ref{Mn}
\by F.~Mangolte
\paper Surfaces elliptiques r\'{e}elles et in\'{e}galit\'{e} de Ragsdale-Viro
\jour Math. Z.
\vol 235
\yr 2000
\issue 2
\pages 213--226
\endref\label{Mangolte}

\ref{Ma}
\by Yu.~Manin
\paper Rational surfaces over perfect fields
\jour Mat. Sbornik
\yr 1967
\issue 1
\pages 141--168
\endref\label{Manin}

\ref{N1}
\by S.~Natanzon
\paper Klein surfaces
\jour Uspekhi Mat. Nauk
\vol 45
\yr 1990
\pages 47--90
\endref\label{Nat}

\ref{N2}
\by S.~Natanzon
\paper Topology of $2$-dimensional covering and meromorphic functions on
real and complex algebraic curves
\jour Selecta Math.
\vol 12
\issue 1
\yr 1993
\pages 251--291
\endref\label{Natanzon}

\ref{NSV}
\by S.~Natanzon, B.~Shapiro, A.~Vainshtein
\paper Topological classification of generic
real rational functions
\jour J. Knot Theory Ramifications
\vol 11
\yr 2002
\issue 7
\pages 1063--1075
\endref\label{NSV}

\ref{Ni}
\by V.~Nikulin
\paper Integer quadratic forms and some of their geometrical applications
\jour Math. USSR--Izv.
\vol 43
\yr 1979
\pages 103--167
\endref\label{Nikulin}

\ref{Or1}
\by S.~Orevkov
\paper Riemann existence theorem and construction of
real algebraic curves
\jour
Ann. Fac. Sci. Toulouse Math. (6)
\vol 12
\issue 4
\pages
517-531
\yr 2003
\endref\label{Orevkov}  

\ref{Or2}
\by S.~Orevkov
\jour Private communications
\endref\label{Orevkov-private}

\ref{P}
\by A.~Protopopov
\paper Topological classification of branched coverings
of the two-dimensional sphere
\inbook Zapiski Nauchnyh Seminarov LOMI
\vol 167
\pages 135--156
\yr 1988
\endref\label{Protopopov}

\ref{Se}
\by W.~Seiler
\paper
Moduln elliptischer Fl\"{a}chen mit Schnitt
\jour Karlsruhe, Thesis, 1982
\endref\label{Sei}

\ref{SV}
\by B.~Shapiro, A.~Vainshtein
\paper Counting real rational functions with all real critical
values
\jour Mosc. Math. J.
\vol 3
\yr 2003
\issue 2
\pages 647--659
\endref\label{SV}


\ref{Si}
\by R.~Silhol
\book Real algebraic surfaces
\bookinfo Lecture Notes in Math.
\vol 1392
\yr 1989
\publ Sprin\-ger-Verlag
\endref\label{Silhol}

\ref{Wa}
\by C. T. C. Wall
\paper Real forms of smooth del Pezzo surfaces.
\jour J. reine angew. Math.
\vol 375/376
\yr 1987
\pages 47--66
\endref\label{Wall}

\ref{We}
\by J.-Y.~Welschinger
\paper Real structures on minimal ruled surfaces
\jour Comment. Math. Helv.
\vol 78
\issue 2
\pages  418-446
\yr 2003
\endref\label{Welschinger}


\ref{Z}
\by V.~Zvonilov
\paper Rigid isotopies of threenomial curves with the maximal number of ovals
\jour Vestnik Syktyvkarskogo Universiteta
\issue 6
\pages  23 pages
\yr 2006
\endref\label{Zvonilov}

\endRefs
\enddocument